\documentclass[12pt]{amsart}

\usepackage{dsfont}
\usepackage{amsxtra}
\usepackage{amsmath}
\usepackage{amscd}
\usepackage{amssymb}
\usepackage{amsfonts}
\usepackage[all]{xy}
\usepackage{mathrsfs}
\usepackage{amsthm} 
\usepackage{color}
\usepackage{upgreek }
\usepackage{tikz}
\usepackage{enumerate}
\usepackage{geometry}
\usepackage{ae}
\usepackage{enumitem}

\numberwithin{itemcounter}{subsection}

\theoremstyle{plain}

\newtheorem*{thmA}{Theorem A}
\newtheorem*{thmB}{Theorem B}
\newtheorem*{thmC}{Theorem C}
\newtheorem{theorem}{Theorem}[section]
\newtheorem{lemma}[theorem]{Lemma}
\newtheorem{definition-lemma}[theorem]{Definition-Lemma}

\newtheorem{proposition}[theorem]{Proposition}
\newtheorem{conjecture}[theorem]{Conjecture}
\newtheorem{corollary}[theorem]{Corollary}
\theoremstyle{definition}
\newtheorem{definition}[theorem]{Definition}
\theoremstyle{remark}
\newtheorem{remark}[theorem]{Remark}
\newtheorem{example}[theorem]{Example}
\numberwithin{equation}{section}

\newenvironment{dedication}
        {\vspace{0ex}\begin{quotation}\begin{center}\begin{em}}
        {\par\end{em}\end{center}\end{quotation}}

\def\bbA{\mathbb{A}}

\def\bbC{\mathbb{C}}

\def\bbF{\mathbb{F}}
\def\bbG{\mathbb{G}}

\def\bbN{\mathbb{N}}

\def\bbQ{\mathbb{Q}}

\def\bbZ{\mathbb{Z}}

\def\scrC{\mathscr{C}}

\def\scrH{\mathscr{H}}
\def\scrI{{I}}

\def\scrL{\mathscr{L}}

\def\scrO{\mathscr{O}}
\def\scrP{\mathscr{P}}

\def\scrU{\mathscr{U}}

\def\scrX{\mathscr{X}}

\def\frakg{\mathfrak{g}}

\def\fraks{\mathfrak{s}}
\def\frakl{\mathfrak{l}}

\def\frakE{\mathfrak{E}}
\def\frakF{\mathfrak{F}}

\def\frakH{\mathfrak{H}}

\def\frakS{\mathfrak{S}}

\def\calG{\mathcal{G}}

\def\calL{\mathcal{L}}

\def\calO{\mathcal{O}}

\def\frakg{\mathfrak{g}}
\def\frakh{\mathfrak{h}}
\def\frakl{\mathfrak{l}}

\def\fraksl{\mathfrak{sl}}

\def\bfb{\mathbf{b}}

\def\bfe{\mathbf{e}}

\def\bfL{\mathbf{L}}

\def\bfB{\mathbf{B}}
\def\bfC{\mathbf{C}}

\def\bfF{\mathbf{F}}
\def\bfG{\mathbf{G}}
\def\bfH{\mathbf{H}}

\def\bfM{\mathbf{M}}
\def\bfN{\mathbf{N}}

\def\bfP{\mathbf{P}}

\def\bfS{\mathbf{S}}
\def\bfT{\mathbf{T}}
\def\bfU{\mathbf{U}}

\def\bfW{\mathbf{W}}

\def\bfZ{\mathbf{Z}}

\def\bfP{\mathbf{P}}
\def\bfT{\mathbf{T}}

\def\geqs{\geqslant}

\def\simto{\overset{\sim}\to}

\def\Unip{{\operatorname{Unip}\nolimits}}
\def\Irr{{\operatorname{Irr}\nolimits}}
\def\wIrr{{\operatorname{WIrr}\nolimits}}
\def\sp{{\operatorname{sp}\nolimits}}

\def\e{{\operatorname{e}\nolimits}}
\def\op{{\operatorname{op}\nolimits}}

\def\GL{\operatorname{GL}\nolimits}
\def\GU{\operatorname{GU}\nolimits}
\def\diag{\operatorname{diag}\nolimits}

\def\k{{\operatorname{k}\nolimits}}

\def\P{\operatorname{P}\nolimits}
\def\Q{\operatorname{Q}\nolimits}
\def\R{\operatorname{R}\nolimits}
\def\X{\operatorname{X}\nolimits}
\def\Y{\operatorname{Y}\nolimits}

\def\add{\operatorname{add}\nolimits}
\def\rem{\operatorname{rem}\nolimits}
\def\top{\operatorname{top}\nolimits}
\def\cl{{\operatorname{cl}\nolimits}}
\def\lw{{\operatorname{\,lw}\nolimits}}
\def\up{{\operatorname{\,up}\nolimits}}

\def\soc{\operatorname{soc}\nolimits}

\def\min{{\operatorname{min}\nolimits}}
\def\Ker{{\operatorname{Ker}\nolimits}}

\def\add{\operatorname{add}\nolimits}

\def\Hom{\operatorname{Hom}\nolimits}
\def\RHom{\operatorname{RHom}\nolimits}

\def\End{\operatorname{End}\nolimits}
\def\Res{\operatorname{Res}\nolimits}
\def\res{\operatorname{res}\nolimits}
\def\Ind{\operatorname{Ind}\nolimits}

\def\id{\operatorname{id}\nolimits}

\def\wt{\operatorname{wt}\nolimits}
\def\hw{{\operatorname{hw}\nolimits}}

\def\height{\operatorname{ht}\nolimits}

\def\umod{\operatorname{-umod}\nolimits}
\def\mod{\operatorname{-mod}\nolimits}

\def\proj{\operatorname{-proj}\nolimits}

\def\height{\operatorname{ht}\nolimits}

\def\ct{{{\text{ct}}}}

\author{O. Dudas, M. Varagnolo, E. Vasserot}

\title[]{Categorical actions on unipotent~representations~I. \\
Finite unitary groups.}

\begin{document}
\begin{abstract}
Using Harish-Chandra induction and restriction, we construct a categorical action 
of a Kac-Moody algebra on the category of unipotent representations of finite
unitary groups in non-defining characteristic. We show that the decategorified
representation is naturally isomorphic to a direct sum of level 2 Fock spaces. From our 
construction we deduce that the Harish-Chandra branching graph coincide
with the crystal graph of these Fock spaces, solving a recent conjecture of
Gerber-Hiss-Jacon. We also obtain derived equivalences between blocks, yielding
Brou\'e's abelian defect groups conjecture for unipotent $\ell$-blocks at
linear primes $\ell$. 
\end{abstract}

\thanks{This research was partially supported by the ANR grant number ANR-10-BLAN-0110 and
ANR-13-BS01-0001-01.}

\maketitle

\begin{dedication}
To Gerhard Hiss for his sixtieth birthday
\end{dedication}

\setcounter{tocdepth}{3}

\tableofcontents

\section*{Introduction}

Let $G =\bfG(\bbF_q)$ be a finite reductive group. The irreducible representations
of $G$ over fields of characteristic $\ell \nmid q$ fall into Harish-Chandra series, which 
are defined in terms of Harish-Chandra induction $R_L^G$ and restriction ${}^*R_L^G$
from proper Levi subgroups $L \subset G$. 
Therefore the classification of the irreducible representations (up to isomorphism)
can be reduced to the following two problems.
\begin{itemize}
  \item[(a)] Classification of the cuspidal irreducible representations. 
  \item[(b)] Determination of the endomorphism algebra of representations obtained by
   Harish-Chandra induction of cuspidal representations.
\end{itemize}
This was achieved by Lusztig in \cite{Lu84} when $\ell = 0$ but it remains open for
representations in positive characteristic for most of the finite reductive
groups. By results of Geck-Hiss-Malle \cite{GHM96}, we know however that the algebras 
$\mathrm{End}_G(R_L^G(V))$ occuring in (b) have a structure of (extended) Hecke algebras
of finite type, only the parameters of the deformation are unknown in general. 

\smallskip

When $G$ is a classical group, e.g., $\mathrm{GL}_n(q),$ $\mathrm{GU}_n(q)$, $\mathrm{Sp}_{2n}
(q)$... it turns out that most of the structure of $\mathrm{End}_G(R_L^G(V))$ does not
depend on the representation $V$. This suggests to rather study the endomorphism
algebra $\mathrm{End}(R_L^G)$ of the Harish-Chandra induction functor $R_L^G$ rather
than the endomorphism algebra of the induced representation. This was already
achieved in \cite{CR} for $G= \mathrm{GL}_n(q)$. In this series of papers we
aim at extending Chuang-Rouquier's framework to other classical groups. 

\smallskip

In this paper we will focus on the case of finite unitary groups $\mathrm{GU}_n(q)$.
We will work with both ordinary representations (characteristic zero) and
modular representations in non-defining characteristic (characteristic $\ell \nmid q$).
More precisely, the field of coefficients $R$ of the representations will be an extension
of either $\bbQ_\ell$ or $\bbF_\ell$. Using the tower of inclusion of groups
$\cdots \subset \mathrm{GU}_n(q) \subset  \mathrm{GU}_{n+2}(q) \subset \cdots$
one can form the abelian category
$$ \scrU_R = \bigoplus_{n \geqslant 0} \mathrm{GU}_n(q)\umod$$
of unipotent representations of the various finite unitary groups. Furthermore,
under mild assumption on $\ell$, we can modify the Harish-Chandra induction and
restriction functors to obtain a adjoint pair $(E,F)$ of functors on $\scrU_R$.
The functor $F$ corresponds to a Harish-Chandra induction from $\mathrm{GU}_n(q)$
to $\mathrm{GU}_{n+2}(q)$ whereas $E$ corresponds to the restriction. Note that
only specific Levi subgroups are considered, and we must work with a variation of
the usual Harish-Chandra theory (the \emph{weak Harish-Chandra theory}) introduced in
\cite{GHJ}. 

\smallskip

In this framework, problem (a) amounts to finding the modules $V$ such that
$EV=0$ and problem (b) is about the structure of $\mathrm{End}_G(F^m V)$
for such cuspidal modules $V$. As mentioned before, most of the structure
of this endomorphism algebra is already contained in $\mathrm{End}(F^m)$. 
In \S\ref{sec:rep-datum}, we construct natural transformations $X$ of $F$ and
$T$ of $F^2$ where $X$ should be thought of as a Jucys-Murphy element and
$T$ satisfies a quadratic relation with eigenvalues $q^2$ and $-1$. This
endows $\mathrm{End}(F^m)$ with a morphism from an affine Hecke algebra
$\bfH_{m}^{q^2}$ of type $A_{m-1}^{(1)}$ with parameter $q^2$. Back to our
original problem, the evaluation at a cuspidal module $V$ provides
a natural map $\bfH_{m}^{q^2} \longrightarrow \mathrm{End}_G(F^mV)$.
Then, we prove that this map induces a natural isomorphism between
$\mathrm{End}_G(F^m V)$ and the quotient of $\bfH_{m}^{q^2}$ by the ideal 
generated by the relation of $X$ on $FV$ if $V$ is an unipotent representation in characteristic zero (see Theorem \ref{thm:cat}).
In that case, $V$ is the unique cuspidal representation of
$\mathrm{GU}_{n}(q)$ with $n =t(t+1)/2$ for some $t \geqslant 0$ and the eigenvalues of
$X$ on $FV$ are $Q_t= \{(-q)^t,(-q)^{-1-t}\}$. The result was already
proved in \cite{HL80} but with some ambiguity on the eigenvalues of $X$. We 
remove that ambiguity by using the structure of unipotent blocks with
cyclic defects.

\smallskip

Having proved that the eigenvalues of $X$ are powers of $-q$, we can
form a Lie algebra $\frakg$ corresponding to the quiver
with vertices $(-q)^\bbZ$ and arrows given by multiplication by $q^2$.
When working in characteristic zero, $\frakg$ is isomorphic to
two copies of $\fraks\frakl_\bbZ$, whereas in positive characteristic
$\ell$ it will depend on the parity of $e$, the order of $(-q)$ modulo
$\ell$. When $e$ is even (\emph{linear
prime} case), it is a subalgebra of $(\widehat{\fraks\frakl}_{e/2})^{\oplus2}$
whereas it isomorphic to $\widehat{\fraks\frakl}_{e}$ when $e$ is
odd (\emph{unitary prime} case). Our main result is that $E$ and $F$
induce a categorical action of $\frakg$ on $\scrU_R$ (see 
Theorems \ref{thm:char0} and \ref{thm:charl}, and 
\S \ref{subsec:categoricalaction} for the definition of categorical actions). 

\begin{thmA} The representation datum $(E,F,X,T)$ given by
Harish-Chandra induction and restriction endows $\scrU_R$
with a structure of categorical $\frakg$-module.
\end{thmA}

Let $E = \bigoplus E_i$ and $F = \bigoplus F_i$ be the decomposition of
the functors into generalized $i$-eigenspaces for $X$. Then $\{[E_i],
[F_i]\}_{i \in (-q)^\bbZ}$ act as the Chevalley generators of $\frakg$
on the Grothendieck group $[\scrU_R]$ of $\scrU_R$ and many
problems on $\scrU_R$ have a Lie-theoretic counterpart. For example,
\begin{itemize}[leftmargin=8mm]
  \item weakly cuspidal modules correspond to highest weight vectors;
  \item the decomposition of $\scrU_R$ into Harish-Chandra series corresponds
  to the decomposition of the $\frakg$-module $[\scrU_R]$ into a direct sum
  of irreducible highest weight modules,
  \item the parameters of the ramified Hecke algebras $\mathrm{End}_G(F^mV)$
  are given by the weight of $[V]$,
  \item the blocks of $\scrU_R$, or equivalently the unipotent $\ell$-blocks, correspond
  to the weight spaces for the action of $\frakg$ (inside a Harish-Chandra series
  if $e$ is even). 
\end{itemize}
Such observations were already used in other situations
(for symmetric groups, cyclotomic rational double affine Hecke algebras or cyclotomic $q$-Schur algebras, etc).

\smallskip

For this dictionnary to be efficient one needs to determine the $\frakg$-module
structure on $[\scrU_R]$. This is done in \S \ref{sec:g-e} by looking at the action
of $[E_i]$ and $[F_i]$ on the basis of $[\scrU_R]$ formed by unipotent characters
(if char$(R) = 0$) or their $\ell$-reduction (if char$(R)=\ell$). 
On this basis the action can be made explicit, and we prove that there is a natural
$\frakg$-module isomorphism
$$  [\scrU_R] \, \mathop{\longrightarrow}\limits^\sim \, \bigoplus_{t \geqslant 0} \bfF(Q_t)$$
between the Grothendieck group of $\scrU_R$ and a direct sum of
level 2 Fock spaces $\bfF(Q_t)$, each of which corresponds to an ordinary
Harish-Chandra series (see Corollary \ref{cor:derivedaction-uk}). Through this
isomorphism, the basis of unipotent characters (or their $\ell$-reduction) is
sent to the standard monomial basis. 

\smallskip

Our original motivation for constructing a categorical action of $\frakg$
on $\scrU_R$ comes from a conjecture of Gerber-Hiss-Jacon \cite{GHJ}, which
predicts an explicit relation between the Harish-Chandra branching graph
and the crystal graph of the Fock spaces $\bigoplus_{t \geqslant 0} \bfF(Q_t)$. 
See also  \cite{GH}.
Using our categorical methods and the unitriangularity of the 
decomposition matrix we obtain a complete proof of the conjecture
(see Theorem \ref{thm:wHC}).

\begin{thmB} Assume $e >1$ is odd. Then the Harish-Chandra branching
graph coincides with the union of the crystal graphs of the Fock
spaces $\bfF(Q_t)$.
\end{thmB}

Note that the construction of these crystal graphs depend on the 
choice of a charge, which is made explicit in \S \ref{sec:isocrystals} and which indeed differs slightly
from the charge used in \cite{GHJ}. Note also that the proof is based on the following two basic ingredients :
\begin{itemize}[leftmargin=8mm]
  \item a unicity statement for crystals of perfect bases which seems to be new (Proposition \ref{prop:perfect-bases}),
  \item a particular choice of partial order on the basis elements of the Fock space which comes for the representation theory of 
rational double affine Hecke algebras and uses the main theorem in \cite{RSVV}.
\end{itemize}

\smallskip

A similar result can be deduced when $e$ is even. However, in 
that case, the situation is already well-understood by work 
of Gruber-Hiss \cite{GrH} on classical groups. The case where 
$e$ is odd (unitary primes) is considered as more challenging
and Theorem B is the first major result in that direction since
the case of $\mathrm{GL}_n(q)$ was solved by Dipper-Du
\cite{DiDu}. This solves completely the problem of classification of
irreducible unipotent modules for unitary groups mentioned at the
beginning of the introduction. More precisely, there are two notions of Harish-Chandra series for
unipotent modules in non defining positive characteristic. Our work describes the weak Harish-Chandra series.
Another categorical construction can be used in order to get the usual (non weak) Harish-Chandra series, by adapting some technics from
\cite{SV}. We mention it very briefly in a conjectural form in \S \ref{sec:heisenberg} in order that the paper remains of a reasonable length.
We will come back to this elsewhere.

\smallskip

By the work of Chuang-Rouquier, categorical actions also provide
derived equivalences between weight spaces. In our situation, these
weight spaces are exactly the unipotent $\ell$-blocks and we obtain
many derived equivalences between blocks with the same local structure.
Together with Livesey's construction of good blocks in the linear
prime case, we deduce a proof of Brou\'e's abelian defect group
conjecture (see Theorem \ref{thm:broue}). 

\begin{thmC} Assume $e > 2$ is even. Then any unipotent $\ell$-block
of a unitary group with abelian defect group is derived equivalent to its Brauer correspondent.
\end{thmC}

Many results and construction can be applied to other classical
groups. This will be the purpose of a subsequent paper. 

\smallskip

The paper is organized as follows. In Section \ref{part:categoricalactions}
we set our notations and recall the definition of categorical
actions. We also record several applications for later use, including
the existence of perfect bases and the construction of derived 
equivalences. In Section \ref{part:fock} we introduce the Fock spaces,
which are certain the level $l$ representations of Kac-Moody algebras.
They have a basis given by charged $l$-partitions, and a crystal
graph which can be defined combinatorially. In Section \ref{part:modular}
we recall standard results on unipotent representations of finite
reductive groups in non-defining characteristic. Section \ref{part:unitarygroups}
is  the core of our paper. We define a representation datum on the
category of representations, which is then shown to induce a categorical
action on unipotent modules. We give two main applications of our
construction, solving the recent conjecture of Gerber-Hiss-Jacon,
and Brou\'e's abelian defect groups conjecture for unipotent
$\ell$-blocks of unitary groups at linear primes $\ell$.
In the last section, we sketch a strategy towards the determination of
the usual (non weak) Harish-Chandra series using the action of a Heisenberg
algebra.

\medskip

\noindent{\bf Acknowledgements.} We are grateful to P. Shan for helpful discussions concerning the material in \S \ref{sec:heisenberg}.

\section{Categorical representations}\label{part:categoricalactions}

Throughout this section, $R$ will denote a noetherian commutative domain (with unit).

\subsection{Rings and categories}
An \emph{$R$-category} $\scrC$ is an additive category enriched over the tensor
category of $R$-modules. All the functors $F$ on $\scrC$ will be assumed to be $R$-linear.
Given such functor, we denote by $1_F$ or sometimes $F$ the identity element in the endomorphism
ring $\End(F)$. The identity functor on $\scrC$ will be denoted by by $1_\scrC$.
A composition of functors $E$ and $F$ is written as $EF$, while a composition of
morphisms of functors (or natural transformations) $\psi$ and $\phi$ is written as
$\psi\circ\phi$. We say that $\scrC$ is \emph{Hom-finite} if the Hom spaces are finitely
generated over $R$. Since $\scrC$ is additive, it is an exact category with split exact
sequences. If the category $\scrC$ is abelian or exact, we denote by $[\scrC]$ the
complexified Grothendieck group and by $\Irr(\scrC)$ the set of isomorphism classes of
simple objects of $\scrC$. The class of an object $M$ of $\scrC$ in the Grothendieck group
is denoted by $[M]$. An exact endofunctor $F$ of $\scrC$ induces a linear map
on $[\scrC]$ which we will denote by $[F]$.

\smallskip

Assume that $\scrC$ is Hom-finite. Given an object $M\in\scrC$ we set $\scrH(M)=
\End_\scrC(M)^\text{op}$. It is an $R$-algebra which is finitely generated as an $R$-module.
Consider the adjoint pair $(\frakE_M,\frakF_M)$ of functors given by
$$  \begin{aligned}
  &\frakF_M=\Hom_\scrC(M,-):\scrC\longrightarrow\scrH(M)\mod,\\
  &\frakE_M=M\otimes_{\scrH(M)}-:\scrH(M)\mod\longrightarrow\scrC.
  \end{aligned}$$
Let $\add(M)\subset\scrC$ be the smallest $R$-subcategory containing $M$ which is
closed under direct summands. Then the functors $\frakE_M,$ $\frakF_M$ satisfy the
following properties : 

\begin{itemize}[leftmargin=8mm]
  \item $\frakE_M,$ $\frakF_M$ are equivalences of $R$-categories between
  $\add(M)$ and $\scrH(M)\proj$,
  \item $\frakF_M\frakE_M=1_{\scrH(M)\mod},$
  \item if $R$ is a field, $M$ is projective and $\scrC$ is abelian and has
  enough projectives, then $\frakF_M$ yields a bijection
  $$\{L\in\Irr(\scrC)\,\mid\,M \twoheadrightarrow L \}\ 
  \mathop{\longleftrightarrow}\limits^{1:1} \ \Irr(\scrH(M)).$$ 
\end{itemize}

Assume now that $\scrC=H\mod$, where $H$ is an $R$-algebra with 1 which is finitely
generated and free over $R$. We abbreviate $\Irr(H)=\Irr(\scrC)$. Given an homomorphism
$R\to S$, we can form the $S$-category $S\scrC=SH\mod$ where $SH=S\otimes_R H$. 
Given another $R$-category $\scrC'$ as above and an exact ($R$-linear) functor
$F:\scrC\to\scrC'$, then $F$ is represented by a projective object $P\in\scrC$. 
We set $SF=\Hom_{S\scrC}(SP,-):S\scrC\to S\scrC'$. Let $K$ be the field of fraction
of $R$, $A\subset R$ be a subring which is integrally closed in $K$ and $\theta:R\to \k$
be a ring homomorphism into a field $\k$ such that $\k$ is the field of fractions of
$\theta(A)$. If $\k H$ is split, then there is a \emph{decomposition map} $d_\theta\,:
\,[KH\mod]\longrightarrow[\k H\mod]$, see e.g. \cite[sec.~3.1]{GJ} for more details.

\subsection{Kac-Moody algebras of type $A$ and their representations}\label{sec:quivers}

The Lie algebras which will act on the categories we will study will always be finite
sums of Kac-Moody algebras of type $A_\infty$ or $A_{e-1}^{(1)}$. They will arise
from quivers of the same type. 

\subsubsection{Lie algebra associated with a quiver}\label{subsec:quivers}
Let $v \in R^\times$ and $\scrI \subset R^\times$. We assume that $v \neq 1$ and that
$\scrI$ is stable by multiplication by $v$ and $v^{-1}$ with finitely many orbits.
To the pair $(\scrI,v)$ we associate a quiver $\scrI(v)$ (also denoted by $\scrI$)
as follows:
\begin{itemize}[leftmargin=8mm]
  \item the vertices of $\scrI(v)$ are the elements of $\scrI$;
  \item the arrows of $\scrI(v)$ are $i \to iv$ for $i \in \scrI$. 
\end{itemize}
Since $\scrI$ is assumed to be stable by multiplication by $v$ and $v^{-1}$, such a
quiver is the disjoint union of quivers of type $A_\infty$ if $v$ is not a root of unity,
or of cyclic quivers of type $A_{e-1}^{(1)}$ if $v$ is a primitive $e$-th root of $1$.

\smallskip

The quiver $\scrI(v)$ defines a symmetric generalized Cartan matrix
$A = (a_{ij})_{i,j\in \scrI}$ with $a_{ii}= 2$, $a_{ij} =-1$ when
$i \rightarrow j$ or $j \rightarrow i$ and $a_{ij}=0$ otherwise. To this Cartan matrix one can
associate the (derived) Kac-Moody algebra $\frakg_\scrI'$ over $\bbC$, which
has Chevalley generators $e_i,f_i$ for $i\in \scrI$, subject to the usual
relations. 

\smallskip

More generally, let $(\X_\scrI,\X_\scrI^\vee,\langle\bullet,\bullet\rangle_\scrI,
\{\alpha_i\}_{i \in \scrI}, \{\alpha_i^\vee\}_{i \in \scrI})$ be a \emph{Cartan datum} associated with $A$, \emph{i.e.}, 
we assume that
\begin{itemize}[leftmargin=8mm]
 \item $\X_\scrI$ and $\X^\vee_\scrI$ are free abelian groups,
 \item the simple coroots $\{\alpha_i^\vee\}$ are linearly independant in $\X_\scrI^\vee$, 
 \item for each $i\in \scrI$ there exists a fundamental weight $\Lambda_i\in\X_\scrI$
 satisfying $\langle \alpha_j^\vee,\Lambda_i \rangle_\scrI = \delta_{ij}$ for all $j\in\scrI$,
 \item $\langle \bullet , \bullet \rangle_\scrI : \X_\scrI^\vee\times \X_\scrI 
 \longrightarrow \bbZ$ is a perfect pairing such that $\langle \alpha_j^\vee,\alpha_i\rangle_\scrI =~a_{ij}$. 
\end{itemize}
Let $\Q_\scrI^\vee = \bigoplus \bbZ \alpha_i^\vee$ be the coroot lattice and $\P_\scrI = \bigoplus \bbZ \Lambda_i$ be the weight lattice.
Then, the Kac-Moody algebra $\frakg_\scrI$ corresponding to this datum is the Lie algebra generated by
the Chevalley generators $e_i, f_i$ for $i \in \scrI$ and the Cartan algebra
$\frakh = \bbC\otimes \X_\scrI^{\vee}$. An element $h \in \frakh$ acts
by $[h,e_i] = \langle h,\alpha_i\rangle e_i$.
The Lie algebra $\frakg_\scrI'$ is the derived 
subalgebra $[\frakg_\scrI,\frakg_\scrI]$.  

\begin{example}\label{ex:onegenerator}
When $\scrI = v^{\bbZ}$ two cases arise.
 \begin{itemize}[leftmargin=8mm]
  \item[(a)] If $\scrI$ is infinite, then $\frakg_\scrI'$ is isomorphic to 
  $\mathfrak{sl}_\bbZ$, the Lie algebra of traceless matrices with finitely many
  non-zero entries.
  \item[(b)] If $v$ has finite order $e$, then $\scrI$ is isomorphic
  to a cyclic quiver of type $A_{e-1}^{(1)}$. We can form $\X^\vee = \Q^\vee 
  \oplus \bbZ \partial$ and $\X = \P \oplus \bbZ \delta$ with $\langle \partial,
  \Lambda_i\rangle = 0$, $\langle \partial,\alpha_i\rangle = \delta_{i1}$ and
  $\delta = \sum_{i \in \scrI} \alpha_i$. The pairing is non-degenerate, and
  $\frakg_\scrI$ is isomorphic to the Kac-Moody algebra
  $$\widehat{\fraks\frakl}_e = \fraks\frakl_e(\bbC) \otimes \bbC[t,t^{-1}] \oplus
  \bbC c \oplus \bbC \partial.$$
  An explicit isomorphism sends $e_{v^i}$ (resp. $f_{v^i}$) to the matrix
  $E_{i,i+1} \otimes 1$ (resp. $E_{i+1,i} \otimes 1$) if $i \neq e$ and $e_1$ 
  (resp. $f_1$) to $E_{e,1}\otimes t$ (resp. $E_{1,e}\otimes t^{-1}$).
  Via this isomophism the central element $c$ corresponds to $\sum_{i \in \scrI}
  \alpha_i^\vee$, and the derived algebra $\frakg_\scrI'$
  to $\widetilde{\fraks\frakl}_e=\fraks\frakl_e(\bbC) \otimes \bbC[t,t^{-1}] \oplus \bbC c$.
 \end{itemize}
\end{example} 

\smallskip

When $\scrI$ is infinite, it will be sometimes useful to consider a completion of
$\frakg_\scrI$ denoted by $\overline{\frakg}_\scrI$, which has $\prod \bbC \alpha_i^\vee \simeq \bbC^\scrI$ as a Cartan subalgebra. This allows
to consider some infinite sums of the generators, such as $c = \sum \alpha_i^\vee$
which is a central element in $\overline{\frakg}_\scrI$. This will not affect
the representation theory of $\frakg_\scrI$ as we will be working with integrable
representations only (see the following section).

\smallskip

Let $S$ be another commutative domain with unit, and $\theta : R \longrightarrow S$ be
a ring homomorphism. Then there is a Lie algebra homorphism $\overline{\frakg}_{\theta(\scrI)}'
\longrightarrow \overline{\frakg}_{\scrI}'$ defined on the Chevalley generators by 
$$e_i   \longmapsto \displaystyle\sum_{\theta(j) = i} e_j \qquad \text{and} \qquad
	 f_i \longmapsto \displaystyle\sum_{\theta(j) = i} f_j.$$

\begin{example} Take $R = \bbZ_\ell$ and $v \in \bbZ_\ell^\times$ which is not
a root of unity. Assume however that the image of $v$ in $S = \bbZ_\ell/\ell \bbZ_\ell$
is an $e$-th root of unity. Then the canonical map $\theta : \bbZ_\ell 
\longrightarrow \bbZ_\ell/\ell \bbZ_\ell$ induces a Lie algebra homomorphism
$\widetilde{\mathfrak{sl}}_e
\longrightarrow \overline{\mathfrak{sl}}_\bbZ$ which sends $e_i$ to $\sum_{j\equiv i} e_j$. 
\end{example}

To avoid cumbersome notation, we may write $\frakg=\frakg_\scrI$, $\P=\P_{\scrI}$,
$\Q^\vee=\Q_{\scrI}^\vee$, etc. when there is no risk of confusion.

\subsubsection{Integrable representations}\label{subsec:oint}
Let $V$ be a $\frakg$-module. Given $\omega \in \X$, the $\omega$-\emph{weight space}
of $V$ is 
  $$V_\omega= \{v\in V\,\mid\,\alpha^\vee\cdot v=\langle\alpha^\vee,\omega\rangle\,v,\,
  \forall \alpha^\vee\in\Q^\vee\}.$$
An \emph{integrable $\frakg$-module} $V$ is a $\frakg$-module on which the action
of the Chevalley generators is locally nilpotent, and which has a weight decomposition 
$V = \bigoplus_{\omega \in \X} V_\omega.$ 
The set $\mathrm{wt}(V) = \{\omega \in \X\, \mid \, V_\omega \neq 0\}$ is the set 
of weights in $V$.

\smallskip
We denote by $\calO^\text{int}$ the category of \emph{integrable highest weight} modules,
\emph{i.e.} $\frakg$-modules $V$ satisfying
\begin{itemize}[leftmargin=8mm]
  \item $V = \bigoplus_{\omega \in \X} V_\omega$ and $\dim V_\omega < \infty$ for all
  $\omega \in \X$,
  \item the action of $e_i$ and $f_i$ is locally nilpotent for all $i \in \scrI$,
  \item there exists a finite set $F \subset \X$ such that 
  $\mathrm{wt}(V) \subset F + \sum_{i \in \scrI} \bbZ_{\leqslant 0} \alpha_i$.
\end{itemize}
Let $\X^+ = \{\omega \in \X\, \mid\, \langle \alpha_i^\vee,\omega\rangle \in \bbN \text{ for all } i \in \scrI\}$ be the set of \emph{integral dominant weights}. 
Given $\Lambda \in  \X^+$, we denote
by $\bfL(\Lambda)$ the simple integrable highest weight module with highest weight $\Lambda$.
Then $\calO^\text{int}$ is semisimple, and any object in $\calO^\text{int}$ is a direct
sum of such simple modules.

\subsubsection{Quantized enveloping algebras}

Let $u$ be a formal variable and $A=\bbC[u,u^{-1}]$.
Let $U_u(\frakg)$ the quantized enveloping algebra over $\bbC(u)$.
Let $U_A(\frakg)\subset U_u(\frakg)$ be
Lusztig's divided power version of $U_u(\frakg)$.
For each integral weight $\Lambda$ the module $\bfL(\Lambda)$ admits a deformed version $\bfL_u(\Lambda)$ over $U_u(\frakg)$
and an integral form $\bfL_A(\Lambda)$ which is the $U_A(\frakg)$-submodule of
$\bfL_u(\Lambda)$ generated by the highest vector $|\Lambda\rangle$.
Let $\calO^\text{int}_u$ be the category consisting of the $\frakg$-modules which are (possibly infinite) direct sums of $\bfL_u(\Lambda)$'s.
If $V_u\in\calO^\text{int}_u,$ then its integral form $V_A$ is the corresponding sum of the modules $\bfL_A(\Lambda)$.
It depends of the choice a family of highest weight vectors of the constituents of $V_u$.

\subsection{Categorical representations on abelian categories}
\label{sec:categoricalactions}

In this section we recall from \cite{CR,R08} the notion of a categorical action of $\frakg$.
It consists of the data of functors $E_i$, $F_i$ lifting the Chevalley
generators $e_i$, $f_i$ of $\mathfrak{g}$, together with an action of an affine
Hecke algebra on $(\bigoplus F_i)^{m}$.

\subsubsection{Affine Hecke algebras and representation data}\label{subsec:hecke}
Let $\scrC$ be an abelian $R$-category and $v \in R^\times$.

\begin{definition} 
A \emph{representation datum} on $\scrC$ is a tuple $(E,F,X,T)$ where
$E$, $F$ are bi-adjoint functors $\scrC\to\scrC$ and $X\in\End(F)^\times$, $T\in\End(F^2)$
are endomorphisms of functors satisfying the following conditions: 
\begin{itemize}[leftmargin=8mm]
  \item[(a)] $1_FT\circ T1_F\circ 1_FT=T1_F\circ 1_FT\circ T1_F$,

  \item[(b)] $(T+1_{F^2})\circ(T-v1_{F^2})=0$,

  \item[(c)] $T\circ(1_FX)\circ T=vX1_F$.
\end{itemize}
\end{definition}

This definition can also be formulated in terms of actions of affine Hecke algebras.
For $m\geqs 1,$ the \emph{affine Hecke algebra} $\bfH_{R,m}^v$ is the $R$-algebra
generated by $T_1,\ldots,T_{m-1}$, $X^{\pm 1}_1,\ldots,X^{\pm 1}_m$
subject to the relations
\begin{itemize}[leftmargin=8mm]
  \item Type $A_{m-1}$ Hecke relations for $T_1,\ldots,T_{m-1}$:
  $$ \begin{aligned}
    &(T_i+1)(T_i-v)=0,\\
    & T_iT_{i+1}T_i=T_{i+1}T_iT_{i+1}\  \ \text{and} \ \ 
    T_iT_j=T_jT_i \ \ \textrm{if}\ \ |i-j|>1,\\
    \end{aligned}$$
  
  \item  Laurent polynomial ring relations for $X_1^{\pm1},\ldots,X_m^{\pm1}$:
  $$X_iX_j=X_jX_i\  \ \text{and} \ \  
    X_iX_i^{-1}=X_i^{-1}X_i=1, $$
   
  \item Mixed relations:
  $$T_iX_{i}T_i=vX_{i+1} \  \ \text{and} \ \ 
    X_iT_j=T_jX_i\ \ \textrm{if}\ \ i-j\neq 0,1.$$
\end{itemize} 
Therefore $\bfH_{R,m}^v$ contains both a finite dimensional Hecke algebra of 
finite type $A_{m-1}$ and a Laurent polynomial ring in $m$ variables. 
We will also set $\bfH_{R,0}^v=R$. 

\smallskip

Given $(E,F)$ a pair of biadjoint functors, and $X\in\End(F)$, $T\in\End(F^2)$,
the tuple $(E,F,X,T)$ is a representation datum if and only if 
for each $m\in\bbN$, the map
$$ \begin{array}{rcl} \phi_{F^m}\, :\, \bfH_{R,m}^v & \longrightarrow & \End(F^m) \\[4pt]
  X_k & \longmapsto &  1_{F^{m-k}} X 1_{F^{k-1}}\\
  T_l & \longmapsto & 1_{F^{m-l-1}}T1_{F^{l-1}}\\
  \end{array}$$
is a well-defined $R$-algebra homomorphism.

\begin{remark}
Transposing an endomorphism of $F^m$ relatively to the adjunction $(F,E)$ yields
a canonical $R$-algebra isomorphism $\End(E^m)=\End(F^m)^\op$, see e.g., 
\cite[sec.~4.1.2]{CR}. Therefore, if $(E,F,X,T)$ is a representation datum,
the morphisms $X,$ $T$ yield morphisms $X\in\End(E)$, $T\in\End(E^2)$ which induce
an $R$-algebra homomorphism $\phi_{E^m}:\bfH_{R,m}^v\longrightarrow\End(E^m)^\op.$
\end{remark}

\subsubsection{Categorical representations}\label{subsec:categoricalaction}
We assume now that $R$ is a field and that
$\scrC$ is Hom-finite. We fix a pair $(\scrI,v)$ as in \S \ref{sec:quivers} and
we denote by $\frakg = \frakg_\scrI$ the Lie algebra associated to that pair.

\begin{definition}[\cite{R08}]\label{df:cat1}
A \emph{$\frakg$-representation} on $\scrC$ consists of a representation datum
$(E,F,X,T)$ on $\scrC$ and of a decomposition $\scrC=\bigoplus_{\omega\in\X}\scrC_\omega$.
For each $i\in \scrI,$ let $F_i,$ $E_i$ be the generalized $i$-eigenspaces
of $X$ acting on $F,$ $E$ respectively. We assume in addition that
\begin{itemize}[leftmargin=8mm]

\item[(a)] $F=\bigoplus_{i\in \scrI} F_i$ and $E=\bigoplus_{i\in \scrI} E_i$,

\item[(b)] the action of $[E_i]$ and $[F_i]$ for $i\in \scrI$ endow $[\scrC]$ with 
a structure of integrable $\frakg$-module such that $[\scrC]_\omega=[\scrC_\omega]$,

\item[(c)] $E_i(\scrC_\omega)\subset\scrC_{\omega+\alpha_i}$ and
$F_i(\scrC_\omega)\subset\scrC_{\omega-\alpha_i}$.
\end{itemize}
\end{definition}

We say that the tuple $(E, F, X, T)$ and the decomposition $\scrC=\bigoplus_{\omega\in \Y}
\scrC_{\omega}$ is a \emph{$\frakg$-categorification} of the integrable $\frakg$-module
$[\scrC]$.

\subsection{Outcomes} 

In Section \ref{part:unitarygroups} we will endow the category of unipotent representations of finite
unitary groups with a structure of $\frakg$-representation. We give here three main
applications of the existence of a categorical action, which we will use in Sections 
\ref{sec:broue} and \ref{sec:isocrystals} to determine:
 \begin{itemize}[leftmargin=8mm]
  \item[(1)] the Hecke algebras associated to cuspidal representations,
  \item[(2)] the branching graph for the parabolic induction and restriction,
  \item[(3)] derived equivalences between blocks.
 \end{itemize}
Note that for most of the results in this section we will assume that $R$ is a field
and that $\scrI$ is finite. In particular $v \in R^\times$ will be a root of unity.

\subsubsection{Minimal categorical representations} \label{subsec:minimalcat}
Let $m \geqslant 0$, $v \in R^\times$ and $\bfH_{R,m}^v$ be the affine Hecke algebra as
defined in \S \ref{subsec:hecke}. We fix a tuple $Q=(Q_1,\ldots,Q_l)$ in $(R^\times)^l$. The
\emph{cyclotomic Hecke algebra} $\bfH^{Q,\,v}_{R,m}$ is the quotient of $\bfH_{R,m}^v$ 
by the two-sided ideal generated by $\prod_{i=1}^l(X_1-Q_{i})$.

\begin{example} (a) If $l=2$, then $\bfH^{Q,\,v}_{R,m}$ is generated by
$X_1$ and $\{T_i\}_{i =1,\ldots,m-1}$. Set $T_0=-Q_{2}^{-1}X_1$ and $u=-Q_{1}Q_{2}^{-1}.$ Then, we have $(T_i+1)(T_i-v)=0$ and $(T_0+1)(T_0-u)=0$. Therefore $\bfH^{Q,\,v}_{R,m}$ is
isomorphic to an Iwahori-Hecke algebra of type $B_{m}$ with parameters $(u,v)$. 

\smallskip
\noindent 
(b) If $q=1$ and $Q_{p}=\zeta^{p-1},$  then $\bfH^{Q,\,v}_{R,m}=RG(l,1,m)$ is the group
algebra of the complex reflection group $(\bbZ/l\bbZ)^m \rtimes \mathfrak{S}_m$.
\end{example}

Assume now that $R$ is a field. Any finite dimensional $\bfH^{Q,\,v}_{R,m}$-module $M$
is the direct sum of the weight subspaces
  $$M_\nu=\{v\in M\,\mid\,(X_r-i_r)^dv=0,\
  r\in[1,m],\, d\gg 0\},\quad\nu=(i_1,\dots,i_m)\in R^m.$$
Decomposing the regular module, we get a system of orthogonal idempotents $\{e_\nu\,;\,\nu
\in R^m\}$ in $\bfH^{Q,\,v}_{R,m}$ such that $e_\nu M=M_\nu$ for each $M$. 
The eigenvalues of $X_r$ are always of the form $Q_i v^j$ for some $i \in \{1,\ldots,j\}$
and $j \in \bbZ$. As a consequence, if we set $\scrI = \bigcup Q_i v^\bbZ$, then
$e_\nu = 0$ unless $\nu \in \scrI$. The pair $(\scrI, v)$ satisfies the assumptions of
\S \ref{sec:quivers} and we can consider a corresponding Kac-Moody algebra $\frakg_\scrI$ and
its root lattice $\Q_\scrI$. Given $\alpha\in \Q_\scrI^+$ of height $m$, let
$e_\alpha=\sum_{\nu} e_\nu$ where the sum runs over the set of all tuples such that
$\sum_{r=1}^m \alpha_{i_r}=\alpha$. The nonzero $e_\alpha$'s are the primitive central
idempotents in $\bfH^{Q,\,v}_{R,m}$. 

\smallskip

To the dominant weight $\Lambda_Q=\sum_{i=1}^l\Lambda_{Q_{i}}$ of $\frakg_\scrI$ and 
any $\alpha \in \Q_\scrI^+$ we associate the following abelian categories: 
  $$\scrL(\Lambda_Q) = \bigoplus_{m\in\bbN}\bfH^{Q,\, v}_{R,m}\mod \quad \text{and} \quad
  \scrL(\Lambda_Q)_{\Lambda_Q-\alpha}=e_\alpha\bfH^{Q,\,v}_{R,m}\mod.$$
For any $m<n$, the $R$-algebra embedding of the affine Hecke algebras $\bfH_{R,m}^v
\hookrightarrow \bfH_{R,n}^v$ given by $T_i\longmapsto T_i$ and $X_j\longmapsto X_j$ 
induces an embedding $\bfH^{Q,\, v}_{R,m}\hookrightarrow \bfH^{Q,\,v}_{R,n}$.
The $R$-algebra $\bfH^{Q,\,v}_{R,n}$ is free as a left and as a right
$\bfH^{Q,\,v}_{R,m}$-module. This yields a pair of exact adjoint functors 
$(\Ind^{n}_m,\,\Res^{n}_m)$ between $\bfH^{Q,\,v}_{R,n}\mod$ and $\bfH^{Q,\,v}_{R,m}\mod$.
They induce endofunctors $E$ and $F$ of $\scrL(\Lambda_Q)$ by $E=\bigoplus_{m\in\bbN}
\Res_m^{m+1}$ and $F=\bigoplus_{m\in\bbN}\Ind_m^{m+1}$. The right multiplication on 
$\bfH^{Q}_{\scrI,m+1}$ by $X_{m+1}$ yields an endomorphism of the functor $\Ind_m^{m+1}$.
The right multiplication by $T_{m+1}$ yields an endomorphism of $\Ind_m^{m+2}.$
We define $X\in\End(F)$ and $T\in\End(F^2)$ by $X=\bigoplus_{m}X_{m+1}$ and 
$T=\bigoplus_{m}T_{m+1}$.

\smallskip

This construction yields a categorification of the simple
highest module $\bfL(\Lambda_Q)$ of $\frakg_\scrI$. 
Indeed, a theorem of Kang and Kashiwara implies that this holds in the more
general setting of cyclotomic quiver Hecke algebras of arbitrary type. 

\begin{theorem}[\cite{KK12}, \cite{K12}]\label{thm:minimalcat}  \hfill  
  \begin{itemize}[leftmargin=8mm]
  \item[$\mathrm{(a)}$] The endofunctors $E$ and $F$ of $\scrL(\Lambda_Q)$ are biadjoint. 
  \item[$\mathrm{(b)}$] The tuple $(E,F,X,T)$ and the decomposition $\scrL(\Lambda_Q)=
  \bigoplus_{\omega\in \X} \scrL(\Lambda_Q)_{\,\omega}$ is a $\frakg_\scrI$-categorification
  of $\bfL(\Lambda_Q)$.
  \qed
  \end{itemize}
\end{theorem}

This categorical representation is called the \emph{minimal categorical
$\frakg_\scrI$-representation} of highest weight $\Lambda_Q$. 

\smallskip

The $\frakg_\scrI$-modules we are interested in are direct sums of various irreducible
highest weight modules $\bfL(\Lambda_Q)$.
This decomposition admits the following categorical counterpart.
Let $(\scrI,v)$ as in \S \ref{sec:quivers},
and $\frakg = \frakg_\scrI$ be a corresponding Kac-Moody algebra. 
Let $(E,F,X,T)$ be a $\frakg$-representation on an abelian $R$-category $\scrC$.
Recall that for any $m \geqslant 0$ we have an $R$-algebra homomorphism 
$\phi_{F^m} \, : \, \bfH_{I,m}^v\ \longrightarrow \End(F^m)^\op$. Given an
object $M$ in $\scrC$, it specializes to an $R$-algebra homomorphism 
  $$ \bfH_{I,m}^v\ \longrightarrow \End(F^m M)^\op =: \scrH(F^mM).$$

\begin{proposition}[\cite{R08}]\label{prop:unicite}
Assume that the simple
roots are linearly independent in $\X$. 
Let $(E,F,X,T)$ be a  representation of $\frakg$ in a abelian $R$-category $\scrC$,
and $M\in\scrC_{\omega}$. Assume that $EM=0$ and $\End_\scrC(M)=R$. Then
  \begin{itemize}[leftmargin=8mm]
    \item[$\mathrm{(a)}$] $\omega \in \X^+$ is an integral dominant weight,    
    \item[$\mathrm{(b)}$] if we write $\Lambda_Q = \sum_{i \in \scrI} \langle 
    \alpha_i^\vee,\omega\rangle \Lambda_i = 
    \sum_{p=1}^l \Lambda_{Q_p}$ for some $Q=(Q_1,\ldots,Q_l) \in \scrI^l$ and $l \geqslant 1$,
    then the map $\phi_{F^m}$ factors to an $R$-algebra isomorphism 
    $$\bfH_{R,m}^{Q,v} \ \mathop{\longrightarrow}\limits^\sim \ \scrH(F^mM).$$
    \qed
  \end{itemize}
\end{proposition}

In the framework of \S \ref{part:unitarygroups}, the functors $E$ and $F$ are the parabolic induction and
restriction functors. Therefore from Proposition \ref{prop:unicite} we deduce that the
endomorphism algebra of the induction of a cuspidal module is a cyclotomic Hecke algebra
whose parameters are given by the weight of the cuspidal module. 

\subsubsection{Perfect bases, crystals and branching graphs}
\label{subsec:perfect}
We start by a review of Kashiwara's theory of perfect bases and crystals.
A good reference is \cite{K95}, or \cite{KOP} for a short review.
We will be working with the Kac-Moody algebra $\frakg$ coming from a pair $(\scrI,v)$ as
in \S \ref{sec:quivers}. 

\begin{definition}
An {\em abstract crystal} is a set $B$ together with maps
$\mathrm{wt} : B \to \P$, $\varepsilon_i,$ $\varphi_i : B \to {\bbZ}
\sqcup \{-\infty \}$ and $\widetilde{e_i},\widetilde{f_i}$ : $B \to
B\sqcup\{0\}$ for all $i \in \scrI$ satisfiying the following properties:
 \begin{itemize}[leftmargin=8mm]
  \item[(a)] $\varphi_i(b)=\varepsilon_i(b)+\left<\alpha^\vee_i,\mathrm{wt}(b)\right>,$
    
  \item[(b)] $\mathrm{wt}(\widetilde{e_i}b)=\mathrm{wt}(b)+\alpha_i$ and $\mathrm{wt}
  (\widetilde{f_i}b)=\mathrm{wt}(b)-\alpha_i,$
  
  \item[(c)] $b=\widetilde{e_i}b'$ if and only if $ \widetilde{f_i}b=b'$, where
  $b, b'  \in B$, $i \in \scrI$,

  \item[(d)] if $\varphi_i(b)=-\infty$, then $\widetilde{e_i}b=\widetilde{f_i}b=0$,
  
  \item[(e)] if $b\in B$ and $\widetilde{e_i}b \in B$,  then
  $\varepsilon_i(\widetilde{e_i}b)= \varepsilon_i(b)-1$ and $\varphi_i(\widetilde{e_i}b)=
  \varphi_i(b)+1$,

  \item[(f)] if $b\in B$ and  $\widetilde{f_i}b \in B$, then 
  $\varepsilon_i(\widetilde{f_i}b)= \varepsilon_i(b)+1$ and $\varphi_i(\widetilde{f_i}b)=
  \varphi_i(b)-1$.
\end{itemize}
\end{definition}

Note that by (a), the map $\varphi_i$ is entirely determined by $\varepsilon_i$ and
$\mathrm{wt}$. We may therefore omit $\varphi_i$ in the data of an abstract crystal
and denote it by $(B,\widetilde e_i,\widetilde f_i, \varepsilon_i, \mathrm{wt})$.   

\smallskip

An isomorphism between crystals $B_1$, $B_2$ is a bijection ${\psi}:B_1\sqcup\{0\} \longrightarrow B_2\sqcup \{0\}$ such that $\psi(0)=0$ which commutes with $\mathrm{wt}$,
$\varepsilon_i$, $\varphi_i$ $\widetilde{f}_i,$, $\widetilde{e}_i$.

\smallskip

To an abstract crystal we can associate a \emph{crystal graph} as follows: the vertices
of the graph are indexed by the elements of $B$ and the arrows are $b 
\mathop{\longrightarrow}\limits^i b'$ for each $b,b' \in B$ and $i \in \scrI$ with
$b = \widetilde{e_i} b'$ (or equivalently $ b' = \widetilde{f_i}b$). The operators
$\widetilde{e_i}$ and $\widetilde{f_i}$ are thus entirely determined by the graph.

\smallskip

Let $V_u$ be an integrable $U_u(\frakg)$-module in $\calO^\text{int}_u$.
Let $V_A$ be an integral form of $V_u$.
A \emph{lower crystal lattice} in $V_u$ is a free $\bbC[u]$-submodule $\calL$ of $V_A$ such that
$V_A=A\calL$, $\calL=\bigoplus_{\lambda\in \X}\calL_\lambda$ with $\calL_\lambda=\calL\cap(V_A)_\lambda$ 
and $\calL$ is preserved by the \emph{lower Kashiwara crystal operators}  $\widetilde e_i^\lw$, $\widetilde f_i^\lw$ on $V_u$.
A \emph{lower crystal basis} of $V_u$ is a pair $(\calL,B)$ where $\calL$ is a lower crystal lattice of $V_u$ and $B\subset\calL/u\calL$ is a basis
such that we have $B=\bigsqcup_{\lambda\in\X} B_\lambda$ where $B_\lambda=B\cap(\calL_\lambda/u\,\calL_\lambda)$,
$\widetilde e_i^\lw(B),\,\widetilde f_i^\lw(B)\subset B\sqcup\{0\}$ and $b'=\widetilde f_i^\lw b$ if and only if $b=\widetilde e_i^\lw b'$ for each
$b,b'\in B$.
A \emph{lower global basis} of $V_u$ (or, in Lusztig terminology, a \emph{canonical basis}) is an $A$-basis $\bfB$ of $V_A$ such that the lattice 
$\calL=\bigoplus_{b\in\bfB}\bbC[u]\,b$ and the basis $B=\{b\,\text{mod}\,u\,\calL\,|\,b\in\bfB\}$ of $\calL/u\,\calL$ form a lower crystal basis.

\smallskip

One defines in a similar way an \emph{upper crystal lattice}, an \emph{upper crystal basis} and an 
\emph{upper global basis} (or a \emph{dual canonical basis}) using the
\emph{upper Kashiwara crystal operators} $\widetilde e_i^\up$, $\widetilde f_i^\up$ on $V_u$, see, e.g., \cite[def.~4.1,4.2]{KOP}.

\smallskip

If $(\calL,B)$, $(\calL^\vee,B^\vee)$ are lower, upper crystal bases, then $(B,\widetilde e_i^\lw,\widetilde f_i^\lw)$, $(B^\vee,\widetilde e_i^\up,\widetilde f_i^\up)$
are abstract crystals. Therefore, the datum of a lower global basis or an upper global basis determines an abstract crystal.

\smallskip

Let $E_i,$ $F_i,$ $u^h$ with $i\in I$, $h\in\X^\vee,$ be the standard generators of $U_u(\frakg)$.
There exists a unique non-degenerate symmetric bilinear form $(\bullet,\bullet)$ on the module $\bfL_u(\Lambda)$ with highest 
weight vector $|\Lambda\rangle$ satisfying
\begin{itemize}[leftmargin=8mm]
\item $(|\Lambda\rangle,|\Lambda\rangle)=1$,
\item $(E_ix,y)=(x,F_iy)$, $(F_ix,y)=(x,E_iy)$, $(u^{h}x,y)=(x,u^{h}y)$, 
\item $(\bfL_u(\Lambda)_\lambda,\bfL_u(\Lambda)_\mu)=0$ if $\lambda\neq\mu$.
\end{itemize}
\smallskip

Any $U_u(\frakg)$-module in $\calO^\text{int}_u$ admits a lower and an upper crystal and global basis.
If $(\calL, B)$ is a lower crystal basis of $\bfL_u(\Lambda)$ then the pair $(\calL^\vee,B^\vee)$ such that
$\calL^\vee=\{x\in\bfL_u(\Lambda)\,|\,(x,\calL)\subset \bbC[u]\}$ and 
$B^\vee$ is the basis of $\calL^\vee/u\,\calL^\vee$ which is dual to $B$ with respect to the non-degenerate bilinear form
$\calL^\vee/u\,\calL^\vee\times\calL/u\,\calL\to\bbC$ induced by $(\bullet,\bullet),$ is an upper crystal basis, see, e.g., \cite[prop.~4.4]{KOP}.
Finally, taking a basis element in $B$ to the dual basis element in $B^\vee$ is a crystal isomorphism 
$(B,\widetilde e_i^\lw,\widetilde f_i^\lw)\to(B^\vee,\widetilde e_i^\up,\widetilde f_i^\up)$ by \cite[lem.~4.3]{KOP}.
Therefore, if $\bfB$ is a lower global basis of $\bfL_u(\Lambda)$ then the dual basis $\bfB^\vee$ with respect to the non-degenerate bilinear form
$(\bullet,\bullet)$ is an upper global basis and the corresponding abstract crystals 
$(B,\widetilde e_i^\lw,\widetilde f_i^\lw)$ and $(B^\vee,\widetilde e_i^\up,\widetilde f_i^\up)$ are canonically isomorphic.

\smallskip

The crystals that we will consider in this paper all come from particular bases of $\frakg_\scrI$-modules called \emph{perfect bases}. 
Let us define them.
Let $V \in \calO^\mathrm{int}$ be an integrable highest weight $\frakg_\scrI$-module. 
Under this assumption we define, for $i \in \scrI$ and $x \in V$
  $$\ell_i(x)= \max\{k\in \bbN\,\mid\, e_i^{k}x\not=0\} = 
  \min\{k\in \bbN\,\mid\, e_i^{k+1}x=0\}$$
with the convention that $\ell_i(0) = -\infty$.
For each integer $k$, we also consider the vector spaces
  $$ V_i^{\leqslant\, k}=\{x\in V\,\mid\,\ell_i(x)\leqslant k\} \quad \text{and} \quad
   V^{\leqslant\, k}=\bigcap_{i\in I}V_i^{\leqslant\, k}.
  $$
Note that $V_i^{\leqslant\, k} = \mathrm{ker}\, e_i^{k+1}$ when $k\geqslant 0$.
Finally we set $V_i^{k}=V_i^{\leqslant\, k}/\,V_i^{<\, k}$.

\begin{definition} \label{def:perfect}
A basis $B$ of $V$ is {\em perfect} if
 \begin{itemize}[leftmargin=8mm]
  \item[(a)] $B=\bigsqcup_{\mu \in \X_\scrI} B_\mu$ where $B_\mu= B \cap V_\mu,$
  \item[(b)]  for any $i\in \scrI$, there exists a map $\bfe_i : B\to B\sqcup\{0\}$
  such that for any $b \in B$, we have
  \begin{itemize}[leftmargin=8mm]
    \item[(i)] if $\ell_i(b)=0$, then $\bfe_ib=0$,
    \item[(ii)] if $\ell_i(b)>0$, then $\bfe_i b\in B$ and
    $e_ib\in\bbC^\times\,\bfe_ib+V_i^{<\,\ell_i(b)-1}$,
  \end{itemize}
  \item[(c)] if $\bfe_i b= \bfe_i b' \neq 0$ for $b, b' \in B$, then $b=b'.$
\end{itemize}
\end{definition}

\smallskip

Any $\frakg$-module in $\calO^{\text{int}}$ admits a perfect basis.
More precisely, we have the following, see, e.g., \cite[sec.~4]{KOP} for a proof.

\begin{proposition}\label{prop:upper-perfect}
 If $V$ is an integrable $\frakg$-module in $\calO^\text{int}$ with a quantum deformation $V_u$, then
the specialization at $u=1$ of an upper global basis of $V_u$ is a perfect basis of $V$. 
\qed
\end{proposition}

\smallskip

To any categorical representation we associate a perfect basis as in
\cite[prop. 6.2]{S}. More precisely, let $R$ be a field (of any characteristic) and
consider a $\frakg$-representation on an abelian artinian $R$-category $\scrC$.
Then, for each $i\in\scrI$ we define the maps
\begin{align*}
&\widetilde E_i\,:\,\Irr(\scrC)\to\Irr(\scrC)\sqcup\{0\},\quad [L]\mapsto [\soc(E_i(L))],\\
&\widetilde F_i\,:\,\Irr(\scrC)\to\Irr(\scrC)\sqcup\{0\},\quad [L]\mapsto [\top(F_i(L))].
\end{align*}

\begin{proposition}\label{prop:PBfromcategorification}
The tuple $\big(\Irr(\scrC),\widetilde E_i,\widetilde F_i\big)$
defines a perfect basis of $[\scrC]$.
\qed
\end{proposition}

We now recall how to construct an abstract crystal from a perfect basis. We 
set $\widetilde e_i=\bfe_i$. For all $b\in B$ we set $\widetilde{f_i}b=b'$ if $\bfe_ib'=b$
for some $b'\in B$, and $0$ otherwise. Then it follows easily from the definition that
$(B,\widetilde{e_i}, \widetilde{f_i},\ell_i, \mathrm{wt})$ is an abstract crystal. 
In the case where the perfect basis comes from a categorical $\frakg$-representation,
the corresponding crystal graph is the \emph{branching graph} of the exact functors $E_i$,
$F_i$. More precisely, it is the colored graph with vertices labelled by $\Irr(\scrC)$
and arrows $[L] \mathop{\longrightarrow}\limits^i [L']$ whenever $L'$ appears in the
head of $F_i L$,
or equivalently when $L$ appears in the socle of $E_i L'$. 

\smallskip

We finish this section with two results which will be important to identify the crystal 
graph obtained by the categorification with the crystal graph of some Fock space (see
\S \ref{subsec:crystalfock} for the definition of the crystal of a charged Fock space).
For each $i\in \scrI$ and $k\in\bbN$, we set
$B^{\leqslant\,k}=V^{\leqslant k}\cap B$ and $B^{\leqslant\,k}_i=V^{\leqslant k}_i\cap B$.
For a given $i\in \scrI$ and for $b\in B,$ let $[b]_i$ be the image of $b$ in
$V_i^{\ell_i(b)}$. We have the following well-known facts. 

\begin{lemma}\label{lem:perfect}
Let $B$ be a perfect basis of $V \in \calO^{\mathrm{int}}$. Let $i \in\scrI$ and $b,b' \in B$.
 \begin{itemize}[leftmargin=8mm]
  \item[$\mathrm{(a)}$]  $b=b'$ if and only if  $\ell_i(b)=\ell_i(b')$ and $[b]_i=[b']_i$.
  \item[$\mathrm{(b)}$] $B^{\leqslant\,k}$ and $B^{\leqslant\,k}_i$ are bases of $V^{\leqslant\,k}$
  and $V^{\leqslant\,k}_i.$
 \end{itemize}
\end{lemma}

\begin{proof}
Let $i \in \scrI$. For each $b\in B$, we set $\bfe_i^{\,+} b=
\bfe_i^{\,\ell_i(b)} b$ and $e_i^{+} b=e_i^{\ell_i(b)} b$. Note that
$e_i V_i^{\leqslant\,k}\subset V_i^{<\,k}$. Applying successively the axiom (b)(ii)
of perfect bases, we get
\begin{align}\label{axiom}
  e_i^{k}b\in\bbC^\times\,\bfe_i^{\,k}b+V_i^{<\,\ell_i(b)-k},\quad 
  \forall k=1,2,\dots,\ell_i(b).
\end{align}
In particular we have $e_i^{+} b\in\bbC^\times\,\bfe_i^{\,+} b$
and $\bfe_i^{\,+} b\in B_i^{\,\leqslant 0}$. Furthermore, the axiom (c) of perfect
bases implies that $\bfe_i^{\,+} b\neq \bfe_i^{\,+} b'$ whenever $b\neq b'$.

\smallskip
Next, let us prove that $B^{\leqslant\,k}_i$ is a basis of $V^{\leqslant\,k}_i$.
It is enough to check that $B^{\leqslant\,k}_i$ spans $V^{\leqslant\,k}_i$, which
we prove by induction on $k$. If $k<0$ this is obvious. Assume that $k\geqslant 0$.
Given $x\in V^{\leqslant\,k}_i$, $x\neq 0$, let $x_b\in\bbC$ be such that
$x=\sum_{b\in B}x_b\,b$. Let $\ell=\max\{\ell_i(b)\,\mid\,x_b\neq 0\}$.
It is enough to check that $\ell\leqslant k$.
Assume that $\ell>k$.
By \eqref{axiom}, there are elements $c_b\in\bbC^\times$ such that
$$0=e_i^\ell x=\sum_{\ell_i(b)=\ell}x_b\,c_b\,\bfe_i^{\,+}b.$$ 
However, the elements  $\bfe_i^{\,+} b$ such that $\ell_i(b)=\ell$ belong to
$B^{\,\leqslant 0}$ and are distinct, hence linearly independent, yielding a contradiction.
Furthermore, since $B^{\leqslant\,k}_i$ is a basis of $V^{\leqslant\,k}_i$ for all
$i\in \scrI$, we deduce that $B^{\leqslant\,k} = \bigcap B^{\leqslant\,k}_i 
 = \bigcap (B \cap V^{\leqslant\,k}_i)$ is a basis of $V^{\leqslant\,k}_i =
 \bigcap V^{\leqslant\,k}_i$ which proves (b).

\smallskip
Now, let us prove that if $b,b'\in B$ are such that $\ell_i(b)=\ell_i(b')=k$ and 
$[b]_i=[b']_i$ in $V_i^k$, then we have $b=b'$. By \eqref{axiom}, we have
$e_i^{+}[b]_i=e_i^{+}b\in\bbC^\times\,\bfe_i^{\,+}b$ and
$e_i^{+}[b']_i=e_i^{+}b'\in\bbC^\times\,\bfe_i^{\,+}b'.$
Thus, if $[b]_i=[b']_i$ we deduce that $\bfe_i^{\,+}b=\bfe_i^{\,+}b'$,
from which we get $b=b'$.
\end{proof}

We deduce the following proposition.

\begin{proposition} \label{prop:perfect-bases}
Let $B$ and $B'$ be perfect bases of $V$. Assume that there is a bijection 
$\varphi\,:\,B\to B'$ and a partial order $\leqslant$ on $B$ such that
$\varphi( b)\in b + \sum_{c > b}\bbC\,c$ for each $b\in B$. Then the map
$\varphi$ is a crystal isomorphism $B\mathop{\to}\limits^\sim B'$.

\end{proposition}

\begin{proof}
For convenience we shall write $\ell_i'$ for the restriction of the map $\ell_i$ to
the basis $B'$ of $V$. 

\smallskip
First observe that $\ell'_i\circ\varphi=\ell_i$ for all $i\in \scrI$.
Indeed, by Lemma \ref{lem:perfect}(b), the sets $B^{\leqslant k}_i=\{b\,\mid\,b\in B,
\,\ell_i(b)\leqslant k\}$ and $(B')^{\leqslant k}_i=\{\varphi(b)\,\mid\,b\in B,\,
\ell'_i(\varphi(b)) \leqslant k\}$ are both bases of $V^{\leqslant\,k}_i$. 
Since $\varphi(b)\in b + \sum_{c>b}\bbC\,c$, we deduce that for all $k \in \bbZ$ 
  $$\ell'_i(\varphi(b))\leqslant k\iff \varphi(b) \in V^{\leqslant\,k}_i\Longrightarrow
  b\in V^{\leqslant\,k}_i \iff \ell_i(b)\leqslant k.$$
Therefore we have $\ell_i(b) \leqslant \ell'_i(\varphi(b))$ and $(B')^{\leqslant k}_i 
\subset \varphi(B^{\leqslant k}_i)$ for all $k \in \bbZ$ and $b \in B$. Using 
$\varphi^{-1}$ with the order on $B'$ induced by $\leq$ and $\varphi$ we get equalities.
We also deduce that $\varphi(b) \in b + \sum_{c}\bbC\,c$ where $c$ runs over elements of
$B$ satisfying 
\begin{align}\label{C0}
  \ell_i(c)\leqslant\ell_i(b)=
  \ell'_i(\varphi(b)),\ \ c>b \ \ \text{and} \ \ \wt(b)=\wt(\varphi(b))=\wt(c).
\end{align}
In particular, the map $\varphi$ yields a weight preserving bijection
$B^{\leqslant 0}\to(B')^{\leqslant 0}$. Therefore it extends to an automorphism
of the $\frakg$-module $V$, and in turn, to a crystal isomorphism
$\psi:B \mathop{\to}\limits^\sim B'$ such that $\psi (b) = \varphi (b)$ for all $b\in B^{\leqslant\,0}$ (see \cite[thm.~5.37]{BK}). 

\smallskip

We claim that $\varphi=\psi$. We will prove it by induction with respect
to $\ell_i$. Fix an element $b\in B$ with $k:=\ell_i(b)>0$ for some $i\in \scrI$,
and assume that $\varphi(c)=\psi (c),\ \forall c\in B_i^{<\,k}.$
We have $\ell'_i(\varphi(b))=\ell_i(b)=k>0$, hence $\bfe_i(b),$ $\bfe_i'(\varphi(b))$ are
both non-zero. Since $\psi$ is a crystal isomorphism, the induction hypothesis applied
to $c=\bfe_i(b)$ yields $\varphi(\bfe_i(b)) = \psi(\bfe_i(b))=\bfe_i'(\psi(b)).$
The axiom (c) of perfect bases would imply that $\varphi(b)=\psi (b)$ if we can
show that
  \begin{align}\label{C1}\varphi(\bfe_i(b))=\bfe_i'(\varphi(b)).
  \end{align}

Let us prove this equality. We have $\ell'_i(\varphi(\bfe_i(b)))=\ell_i(\bfe_i(b))=k-1$ and
$\ell'_i(\bfe_i'(\varphi(b)))=\ell'_i(\varphi(b))-1=k-1$. Therefore by Lemma \ref{lem:perfect}(a)
it is enough to check that
  \begin{align}\label{C2}
  [\varphi(\bfe_i(b))]_i=[\bfe_i'(\varphi(b))]_i\ \ \text{in}\ \,V_i^{k-1}
  \end{align}
which we shall prove by induction with respect to the order $\geqslant$.
\smallskip

Recall that the map $\varphi$ is unitriangular in the basis $B$.
Therefore the same holds for $\varphi^{-1}$ in the basis $B'$ and we have, by projection
to $V_i^k$ and using \eqref{C0} 
  \begin{align}\label{C3}
  [b]_i\in [\varphi(b)]_i + \sum_{\begin{subarray}{c} \ell_i(c) = k \\ c> b\end{subarray}}
  \bbC\,[\varphi(c)]_i.
  \end{align}  
Since $e_i V_i^{\leqslant\,k}\subset V_i^{<\,k}$, the map $e_i$ factors through a linear map
$e_i\,:\,V_i^{k}\longrightarrow V_i^{k-1}.$
The axiom (b)(ii) of perfect bases implies that
$e_i([b]_i)\in\bbC^\times\,[\bfe_i(b)]_i\ \text{in}\ V_i^{k-1}.$
Now, applying $e_i$ to \eqref{C3} yields following relation in $V_i^{k-1}$
  \begin{align}\label{C4}
  [\bfe_i (b)]_i\in\bbC^\times\, [\bfe_i'(\varphi(b))]_i + 
  \sum_{\begin{subarray}{c} \ell_i(c) = k \\ c> b\end{subarray}}\bbC\,
  [\bfe_i' (\varphi(c))]_i.
  \end{align}
Assume now by induction that \eqref{C2} holds for any $c > b$. Then we can rewrite \eqref{C4}
as
  \begin{align}\label{C6}
  [\bfe_i (b)]_i\in \bbC^\times\,[\bfe_i'(\varphi(b))]_i + 
  \sum_{\begin{subarray}{c} \ell_i(c) = k \\ c> b\end{subarray}}\bbC\,[\varphi(\bfe_i (c))]_i.  \end{align}
On the other hand, applying \eqref{C3} to $\bfe_i (b)$ instead of $b$, we also get the following relation in $V_i^{k-1}$
  \begin{align}\label{C5}
  [\bfe_i (b)]_i\in [\varphi(\bfe_i (b))]_i + \sum_{\begin{subarray}{c} \ell_i(c) = k-1 \\ 
  c> \bfe_i (b)\end{subarray}} \bbC\,[\varphi(c)]_i.
  \end{align}  
Now, observe that $[\varphi(\bfe_i(b))]_i\notin\bbC [\varphi(\bfe_i (c))]_i$
whenever $\ell_i(c)=k$ and $c\neq b$, by Lemma \ref{lem:perfect}(a) and the axiom (c)
of perfect bases. Therefore, comparing \eqref{C6} and \eqref{C5}, we get the identity \eqref{C2}.

\smallskip
Finally, in the case where $b$ is maximal with respect to $\geq$, the relation \eqref{C4}
becomes $[\bfe_i (b)]_i\in \bbC^\times\,[\bfe_i'(\varphi(b))]_i$ and therefore \eqref{C6}
still holds. Hence \eqref{C2} holds in that case as well. 
\end{proof}

\subsubsection{Derived equivalences} 
Given $V$ an integrable $\frakg$-module, and $i \in \scrI$ one can consider
the action of the simple reflection $s_i = \exp(-f_i) \exp(e_i) \exp(-f_i)$
on $V$. 
For each weight $\omega\in \X$, this action maps a weight space $V_\omega$ to $V_{s_i(\omega)}$ with
$s_i(\omega) = \omega - \langle \alpha_i^\vee,\omega \rangle \alpha_i$. 
If $\scrC$ is a categorification of $V$, then it restricts to an $\fraks\frakl_2(\bbC)$-categorification in the sense 
of Chuang-Rouquier. In particular, the simple objects are weight vectors for the categorical $\fraks\frakl_2(\bbC)$-action.
Thus, the theory of Chuang-Rouquier can be applied and
\cite[thm. 6.6]{CR} implies that $s_i$ 
can be lifted to a derived equivalence $\Theta_i$ of $\scrC$.

\begin{theorem}\label{thm:reflection}
Assume that $R$ is a field. 
Let $(E,F,X,T)$
be a representation of $\frakg$ in a abelian $R$-category $\scrC$, and
$i \in \scrI$. Then there exists a derived self-equivalence $\Theta_i$ of
$\scrC$ which restricts to derived equivalences
 $$ \Theta_i \, : \, D^b(\scrC_\omega) \mathop{\longrightarrow}\limits^\sim
 D^b(\scrC_{s_i(\omega)})$$
for all weight $\omega \in \X$. Furtermore, $[\Theta_i] = s_i$ as a linear map
of $[\scrC]$.
\qed
\end{theorem}

In the context of \S \ref{part:unitarygroups}, each weight space $\scrC_\omega$ will be a unipotent
block of a finite unitary group. As a consequence of this theorem we will obtain
many derived equivalences between unipotent blocks, in the spirit of Brou\'e's 
abelian defect group conjecture (see Section \ref{subsec:broue}).

\section{Representations on Fock spaces}\label{part:fock}

Let $R$ be a noetherian commutative domain with unit. As in \S\ref{sec:quivers}, we fix an element $v \in R^\times$
and a subset $\scrI$ of $R^\times$ which is stable by multiplication by $v$ and $v^{-1}$.
We explained in \S\ref{subsec:quivers} how one can associate a Lie algebra $\frakg=\frakg_\scrI$ to this data. 
In this section we recall the construction of (charged) Fock spaces which are particular
integrable representations of $\frakg$. These will be the representations that we shall
categorify using unipotent representations of finite unitary groups (see 
\S \ref{sec:g-infty} and \S \ref{sec:g-e}). 

\subsection{Combinatorics of $l$-partitions}\label{sec:combinatorics}

\subsubsection{Partitions and $l$-partitions}\label{subsec:partitions}
A \emph{partition} of $n$ is a non-increasing sequence of non-negative integers
$\lambda = (\lambda_1 \geqslant \lambda_2 \geqslant \cdots)$ whose terms add up to $n$. We denote by 
$\scrP_n$ be the set of partitions of $n$ and by $\scrP=\bigsqcup_n\scrP_n$ be the set
of all partitions. Given a partition $\lambda$, we write $|\lambda|$ for the \emph{weight}
of $\lambda$, $l(\lambda)$ for the number of non-zero parts in $\lambda$ and
${}^t\lambda$ for the transposed partition.
We associate to $\lambda=(\lambda_1,\lambda_2,\dots)$ the \emph{Young diagram} $Y(\lambda)$
defined by $Y(\lambda)=\{(x,y)\in\bbZ_{>0}\times\bbZ_{>0}\,\mid\,y\leqslant \lambda_x\}.$
It may be visualised by an array of boxes in left justified rows with $\lambda_x$ boxes in
the $x$-th row. If $\lambda$, $\mu$ are partitions of $n$ then we write $\lambda\geqslant\mu$
if for all $n\geqslant i\geqslant 1$ we have $\sum_{j=1}^i\lambda_j\geqslant\sum_{j=1}^i\mu_j$.
This relation defines a partial order on $\scrP$ called the \emph{dominance order}.

\smallskip

An \emph{$l$-partition} of $n$ is an $l$-tuple of partitions whose weights add up to $n$.
We denote by $\scrP^l_n$ be the set of $l$-partitions of $n$ and by $\scrP^l=\bigsqcup_n\scrP^l_n$
the set of all $l$-partitions. The Young diagram of the $l$-partition $\lambda=(\lambda^{1},\ldots,\lambda^{l})$ is the set $Y(\lambda)=\bigsqcup_{p=1}^lY(\lambda^p)\times\{p\}.$
Its weight is the integer $|\lambda|=\sum_p|\lambda^{p}|.$ 

\subsubsection{Residues and content}\label{subsec:contents}
We fix $Q = (Q_1,\ldots,Q_l) \in \scrI^l$.
Let $\lambda$ be an $l$-partition and $A=(x,y,p)$ be a node in $Y(\lambda)$. 
The \emph{$(Q,v)$-shifted residue} of the node $A$ is the element of $\scrI$ given by
$\res(A,Q)_\scrI=v^{y-x}Q_{p}.$  Let $n_i(\lambda,Q)_\scrI$ be the number of nodes of
$(Q,v)$-shifted residue $i$ in $Y(\lambda)$. If $\lambda$, $\mu$ are $l$-partitions such
that $|\mu|=|\lambda|+1$ we write $\res(\mu-\lambda,Q)_\scrI=i$ if $Y(\mu)$ is obtained by
adding a node of $(Q,v)$-shifted residue $i$ to $Y(\lambda)$.
We denote by $\add_i(\lambda,Q)_\scrI$ (resp. $\rem_i(\lambda,Q)_\scrI$) the set of 
\emph{addable nodes} (resp. \emph{removable nodes}) of $(Q,v)$-shifted residue $i$.
With $n=|\lambda|$ we have
\begin{align*}
\add_i(\lambda,Q)_\scrI=\{A=Y(\mu)\setminus Y(\lambda)\,\mid\,\mu\in\scrP^l_{n+1}\,\text{s.t.}\, 
\res(\mu-\lambda,Q)_\scrI=i\},\\
\rem_i(\lambda,Q)_\scrI=\{A=Y(\lambda)\setminus Y(\mu)\,\mid\,\mu\in\scrP^l_{n-1}\,\text{s.t.}
\,\res(\lambda-\mu,Q)_\scrI=i\}.
\end{align*}

A \emph{charge} of the tuple $Q = (Q_1,\ldots,Q_l)$ is an $l$-tuple of integers
$s = (s_1,\ldots,s_l)$ such that $Q_p = v^{s_p}$ for all $p = 1,\ldots,l$.
Conversely, given $\scrI \subset R^\times$ and $v \in R^\times$ as in \S\ref{sec:quivers}, any
$\ell$-tuple of integers $s \in \bbZ^l$ define a tuple $Q = (v^{s_1},\ldots,v^{s_l})$ with
charge $s$. 
The \emph{$s$-shifted content} of the box $A=(x,y,p)$ is the integer $\text{ct}^s(A)=s_p+y-x$.
It is related to the residue of $A$ by the formula $\res(A,Q)_\scrI=v^{\text{ct}^s(A)}$.
We will also write $p(A)=p$.
We will call \emph{charged $l$-partition} a pair $(\mu,s)$ in $\scrP^l\times\bbZ^l$.

\subsubsection{$l$-cores and $l$-quotients}\label{subsec:l-core}
We start with the case $l=1$. The set of \emph{$\beta$-numbers} of a charged partition
$(\lambda,d) \in \scrP\times\bbZ$ is the set given by $\beta_d(\lambda)=\{\lambda_u+d+1-u\,\mid
\,u\geqslant 1\}.$ 
The charged partition $(\lambda,d)$ is uniquely determined by the set $\beta_d(\lambda)$.

\smallskip

If $Q = v^d$, then the set $\mathrm{rem}_i(\lambda,Q)$ of removable nodes of $(Q,v)$-shifted
residue $i$ is the set of integers $j \notin \beta_d(\lambda)$ such that
$i = v^{j}$ and $j+1 \in \beta_d(\lambda)$. Removing such a node has the effect of replacing
$j+1$ by $j$. We have an analogue description for addable nodes.
More generally, for any positive integer $e $, a \emph{$e$-hook} of $(\lambda,d)$ is a pair $(x,x+e)$ such
that $x+e\in \beta_d(\lambda)$ and $x\not\in \beta_d(\lambda)$. Removing the $e$-hook
$(x,x+e)$ corresponds to replacing $x+e$ with $x$ in $\beta_d(\lambda)$. We say
that the charged partition $(\lambda,d)$ is an \emph{$e$-core} if it does not have any
$e$-hook. This does not depend on $d$. 

\smallskip

Next, we construct a bijection $\tau_l:\scrP\times\bbZ\to\scrP^l\times\bbZ^l$.
It takes the pair $(\lambda,d)$ to $(\mu,s),$ where $\mu=(\mu^1,\dots,\mu^l)$ is an
$l$-partition and $s=(s_1,\dots,s_l)$ is a $l$-tuple in
  $$\bbZ^l(d)=\{s\in\bbZ^l\,\mid\,s_1+\cdots+s_l=d\}.$$ 
The bijection is uniquely determined by the following relation
$$  \beta_d(\lambda)=\bigsqcup_{p=1}^l\big(p-l+l\beta_{s_p}(\mu^p)\big).$$
See \cite{Y05} for details.

\smallskip

The bijection $\tau_l$ takes the pair $(\lambda,0)$ to $(\lambda^{[l]},\lambda_{[l]}),$ 
where $\lambda^{[l]}$ is the \emph{$l$-quotient} of $\lambda$ and $\lambda_{[l]}$ lies
in $\bbZ^l(0)$. Since $\lambda$ is an $l$-core if and only if $\lambda^{[l]}=\emptyset$,
this bijection identifies the set of $l$-cores and $\bbZ^l(0)$. We define the \emph{$l$-weight} $w_l(\lambda)$ of the partition $\lambda$ to be the weight of its $l$-quotient. Equivalently,
it is the number of $l$-hooks that can be successively removed from $\lambda$ to get its $l$-core.
So, if we view $\lambda_{[l]}$ as the $l$-core of $\lambda,$ we get 
  $$w_l(\lambda)=|\lambda^{[l]}|=\big(|\lambda|-|\lambda_{[l]}|\big)/l.$$

We will mostly consider the bijection $\tau_l$ for $l=2$.
In particular, a 2-core is either $\Delta_0=\emptyset$ or a triangular partition $\Delta_t=(t,t-1,\dots,1)$ with $t\in\bbN$.
We abbreviate $\sigma_t=(\Delta_t)_{[2]}$, and we write $\sigma_t=(\sigma_1,\sigma_2)$. We have
\begin{align}\label{sigma}
\sigma_t=\begin{cases}\big(-t/2,\,t/2\big)&\ \text{if}\ t\ \text{is\ even},\\
\big((1+t)/2,-(1+t)/2\big)&\ \text{if}\ t\ \text{is\ odd}.
\end{cases}
\end{align}
For each bipartition $\mu$, let $\varpi_t(\mu)$ denote the unique partition 
with 2-quotient $\mu$ and 2-core $\Delta_t.$ 
Thus, the bijection $\tau_2$ maps $(\varpi_t(\mu),0)$ to the pair $(\mu,\sigma_t)$.

\subsection{Fock spaces}
\label{subsec:Fock}
For a reference for the results presented in this section, see for example \cite{U}, \cite{Y05}.
Let $Q = (Q_1,\ldots,Q_l) \in \scrI^l$. It defines an integral dominant weight $\Lambda_Q=\sum_{p=1}^l\Lambda_{Q_{p}} \in \P^+$. The \emph{Fock space} $\bfF(Q)_{\scrI}$ is the $\bbC$-vector
space with basis $\{|\lambda,Q\rangle_{\scrI}\,|\,\lambda\in\scrP^l\}$ called  the 
\emph{standard monomial basis}, and action of $e_i, f_i$ for all $i \in \scrI$ given by
\begin{equation}\label{eq:EF}
  f_i(|\lambda,Q\rangle_{\scrI})=\sum_{\mu}|\mu,Q\rangle_{\scrI},\qquad
  e_i(|\mu,Q\rangle_{\scrI})=\sum_{\lambda}|\lambda,Q\rangle_{\scrI},
\end{equation}
where the sums run over all partitions such that $\res(\mu-\lambda,Q)_\scrI=i$.
This endows $\bfF(Q)_{\scrI}$ with a structure of $\frakg'$-module.
The Fock space $\bfF(Q)_{\scrI}$ can also be equipped with a symmetric non-degenerate bilinear
form $\langle\bullet,\bullet\rangle_{\scrI}$ for which the standard monomial basis is orthonormal.
To avoid cumbersome notation, we shall omit the subscript $\scrI$ when not necessary.

\smallskip

It is easy to see that each element of the standard monomial basis is a weight vector whose
weight can be explicitely computed by the formula
$$\alpha_i^\vee(|\lambda,Q\rangle)=\big(|\add_i(\lambda,Q)|- 
 |\rem_i(\lambda,Q)|\big)\,|\lambda,Q\rangle$$
for $i \in \scrI$. In particular, the vector $|\emptyset,Q\rangle_{\scrI}$ has weight
$\Lambda_Q$.

\begin{proposition}\label{prop:subrepinfock}
 The $\frakg'$-submodule of $\bfF(Q)$ generated by $|\emptyset,Q\rangle$ is
 isomorphic to $\bfL(\Lambda_Q)$. Furthermore, if $\scrI= A_\infty$, then
 $\bfF(Q)=\bfL(\Lambda_Q)$. 
 \qed
\end{proposition}

Using the minimal categorification $\scrL(\Lambda_Q)$ of $\bfL(\Lambda_Q),$ the isomorphism $\bfF(Q)=\bfL(\Lambda_Q)$
can be made more explicit. 
To explain this, let us first recall briefly the definition of the \emph{Specht modules}.
Assume that $R$ has characteristic 0 and contains a primitive $l$-th root $\zeta$ of 1, so
$R$ is a splitting field of the complex reflection group $G(l,1,m)$.
Let $\Irr(R\frakS_m)=\{\phi_\lambda\,|\,\lambda\in\scrP_m\}$
be the standard labelling of the characters of the symmetric group.
Then $$\Irr(RG(l,1,m))=\{\scrX_{\lambda}\,|\,\lambda\in\scrP^l_m\}$$
is the labelling of the simple modules such that $\scrX_\lambda$ is induced from the 
$G(l,1,|\lambda_1|)\times\ldots\times G(l,1,|\lambda_l|)$-module
$$\phi_{\lambda^{(1)}}\chi^0\otimes
\phi_{\lambda^{(2)}}\chi^1\otimes\cdots\otimes\phi_{\lambda^{(l)}}\chi^{l-1}.$$
Here, we denote by $\chi^p$ the one dimensional module of the $|\lambda_p|$-th cartesian power of 
the cyclic group $G(l,1,1)$ given by the $p$-th 
power of the determinant, see, e.g., \cite[sec.~5.1.3]{GJ}.
Recall that (for every field $R$) the $R$-algebra $\bfH^{Q,\,v}_{R,m}$ is split and
that it is semi-simple if and only if we have, see, e.g., \cite[sec.~3.2]{Ma},
\begin{equation*}\label{(A)}
\prod_{i=1}^m(1+v+\cdots+v^{i-1})\,\prod_{a<b}\prod_{-m<r<m}(v^r\,Q_{a}-Q_{b})\neq 0.
\end{equation*}
Thus, by Tits' deformation theorem, under the evaluation $v\mapsto 1$ and $Q_p\mapsto\zeta^{p-1}$, the labelling of $\Irr(R G(l,1,m))$ yields a 
canonical labelling  
$$\Irr(\bfH^{Q,\,v}_{R,m})=\{S(\lambda)_R^{Q,v}\,|\,\lambda\in\scrP^l_m\}.$$
Now, if $R$ is a commutative domain with fraction field $K$ of characteristic 0 as above, we define 
the $\bfH^{Q,\,v}_{R,m}$-module $S(\lambda)^{Q,\,v}_R$ as in \cite[sec.~2.4.3]{RSVV} or \cite[sec.~5.3]{GJ},
using $S(\lambda)^{Q,\,v}_K$ and the dominance order on $\scrP^l_m$,
and if $\theta\,:\,R\to\k$ 
is a ring homomorphism such that $\k$ is the fraction field of $\theta(R)$
we set $S(\lambda)^{Q,\,v}_\k=\k S(\lambda)^{Q,\,v}_R$.
Then, we have the following, see, e.g., \cite{RSVV}.

\begin{proposition}\label{prop:explicitiso}
Let $R$ be a field of characteristic 0 which contains a primitive $l$-th root of 1.
The composition $[\scrL(\Lambda_Q)] \simto \bfL(\Lambda_Q) \to \bfF(Q)$ obtained from
Theorem \ref{thm:minimalcat} and Proposition \ref{prop:subrepinfock} sends the class
of $S(\lambda)^{Q,\,v}_R$ to the standard monomial $|\lambda,Q\rangle$.
\qed
\end{proposition}

For each $p = 1,\ldots, l$, let $\scrI_p$ be the subquiver of $\scrI$ corresponding
to the subset $v^\bbZ Q_p$ of $\scrI$. We define a relation on $\{1,\ldots,l\}$ by 
$i \sim j \iff \scrI_i = \scrI_j$. Let  $\Omega = \{1,\ldots,l\}/\sim$ be the set of 
equivalence classes for this action. Given $p \in \Omega$, we denote by $Q_p$ the tuple
of $(Q_{i_1},\ldots,Q_{i_r})$ where $(i_1,\ldots,i_r)$ is the set of ordered elements
in $p$. The decomposition $\scrI = \bigsqcup_{p \in \Omega} \scrI_p$ yields a canonical
decomposition of Lie algebras
 $\frakg'_\scrI  = \bigoplus_{p \in \Omega} \frakg'_{\scrI_p}$.
The corresponding decomposition of Fock spaces is given in the following proposition.

\begin{proposition}\label{prop:tensorfock}
The map $|\lambda,Q\rangle_\scrI \longmapsto \otimes_{p \in \Omega} |\lambda^p,
Q_p\rangle_{\scrI_p}$ yields an isomorphism of $\frakg'_I$-modules
 $$ \bfF(Q)_\scrI \ \mathop{\longrightarrow}\limits^\sim\ \bigotimes_{p\in \Omega}\bfF({Q_p})_{\scrI_p}.$$
 \qed
\end{proposition}

\subsection{Charged Fock spaces}\label{sec:chargedfock}

A \emph{charged Fock space} is a pair $\bfF(s)=(\bfF(Q),s)$ such that $s\in\bbZ^l$ is a charge
of $Q$, that is  $Q = (v^{s_1},\ldots,v^{s_l})$. Throughout this section, we will always
assume that $\scrI$ is either of type $A_\infty$ or a cyclic quiver. For more general quivers we
can invoke Proposition \ref{prop:tensorfock} to reduce to that case.

\subsubsection{The $\frakg$-action on the Fock space}\label{subsec:xgrading}
The action of $\frakg'$ on $\bfF(Q)$ can be extended to an action of $\frakg$ 
when $Q$ admits a charge $s$. We describe this action in the case where $v$ has finite
order $e$, and $l=1$. In that case $\scrI = v^{\bbZ}$ is isomorphic to the cyclic quiver  $A^{(1)}_{e-1}$
and the charge $s$ is just an integer $d \in \bbZ$ such that $Q = v^d$. If we fix the affine
simple root to be $\alpha_1$, then $\X = \P \oplus \bbZ \delta$ and
$\X^\vee = \Q^\vee \oplus \bbZ \partial$ with $\delta = \sum_{i \in \scrI} \alpha_i$
and $\partial = \Lambda_1^\vee$ (see Example \ref{ex:onegenerator} for more details). 

\smallskip

Given $l \in \bbN$, $l\neq 0$,  and $s = (s_1,\ldots, s_l) \in \bbZ^\ell$, we define 
\begin{align}\label{Delta}
  \Delta(s,e)={1\over 2}\sum_{j=1}^l\big(\bar s_j(1-\bar s_j/e)+s_j(s_j/e-1)\big),
\end{align}
where $\bar s_j$ is the residue of $s_j$ modulo $e$ in $[0,e-1]$. 
Then, we define the action of the derivation $\partial$ on $\bfF(d)=(\bfF(Q),d)$
by
 $$ \partial(|\lambda,Q\rangle) = -\big(n_1(\lambda,Q)+\Delta(d,e)\big)|\lambda,Q\rangle.$$
For this action the weight of a standard basis element is 
\begin{equation}\label{eq:weightlevel1}
 \mathrm{wt}(|\lambda,Q\rangle) = \Lambda_Q-\sum_{i\in \scrI} n_i(\lambda,Q)\,\alpha_i-
 \Delta(d,e)\,\delta.
\end{equation}

Recall from \S\ref{subsec:l-core} that to a charged partition $(\lambda,d)$ we can associate
via $\tau_e$ a pair consisting of an $e$-partition (the $e$-quotient) and a $e$-tuple
of integers adding up to $d$. When $\lambda$ is an $e$-core, we have the following
formula.

\begin{lemma}\label{lem:Delta} If 
$\tau_e(\lambda,d)=(\emptyset,s)$ with $s\in\bbZ^e(d)$, then
$n_1(\lambda,Q)=\Delta(s,1)-\Delta(d,e)$.
\qed
\end{lemma}

\noindent 
Consequently, on an $e$-core $\lambda$ the action of $\partial$ is given by multiplication
by $-\Delta(s,1)$. 

\smallskip

We now describe the action of the affine Weyl group of $\frakg$ on $\bfF(d)$. 
For $i \in \scrI \smallsetminus\{1\}$, we denote by $\alpha_i^\cl = 2 \Lambda_i - 
\Lambda_{iv}-\Lambda_{iv^{-1}}$ and $\Lambda_i^\cl=\Lambda_i-\Lambda_1$ the $i$-th simple root and fundamental weight of
$\fraks \frakl_e$. These (classical) simple roots span the lattice of classical roots $\Q^\cl$.
It is a sublattice of $\P$ of rank $e-1$. The affine Weyl group of $\frakg$ is $W = 
\mathfrak{S}_\scrI \ltimes \Q^\cl$. It acts linearly on $\X$. We will denote by
$t_\gamma \in \mathrm{End}(\X)$ the action of an element $\gamma \in \Q^\cl$ 
as defined in \cite[chap. 6]{Kac}, i.e., for each $\alpha\in\X$ we set
$$t_\gamma(\alpha)=\alpha+(\alpha:\delta)\,\gamma-(\alpha:\gamma)\,\delta-\frac{1}{2}(\alpha:\delta)(\gamma:\gamma)\,\delta$$
where $(\bullet:\bullet)$ is the standard symmetric non-degenerate bilinear form on $\X\times \X$.
For each tuple $s\in \bbZ^{\scrI}$ we consider the element
$\pi_s=\sum_{i\in\scrI}(s_i-s_{iv})\,\Lambda_{i}$. If $s \in \bbZ^{\scrI}(d)$,
then $\pi_s - \Lambda_Q^\cl \in \Q^\cl$ and we can consider the 
corresponding operator $t_{\pi_s- \Lambda_Q^\cl} \in \mathrm{End}(\X)$, from which
we can compute the weight of $|\lambda,Q\rangle$, and the action of $W$ as follows.

\begin{proposition}\label{prop:Delta} 
Let $\lambda$ and $\nu$ be two partitions. Let $(\lambda^{[e]},s) = \tau_e(\lambda,d)$ where
$\lambda^{[e]}$ is the $e$-quotient of $\lambda$ and $s \in \bbZ^\scrI(d)$.
\begin{itemize}[leftmargin=8mm]
  \item[$\mathrm{(a)}$] The weight of $|\lambda,Q\rangle$ equals
   $$ \mathrm{wt}(|\lambda,Q\rangle)  = t_{\pi_s-\Lambda_Q^\cl}(\Lambda_Q)-w_e(\lambda)
   \,\delta$$
 where $w_e(\lambda) = |\lambda^{[e]}|$ is the $e$-weight of $\lambda$.
  \item[$\mathrm{(b)}$] The weights of $|\lambda,Q\rangle$ and $|\nu,Q\rangle$
  are $W$-conjugate if and only if $w_e(\lambda)=w_e(\nu)$.
\end{itemize}
\end{proposition}

\begin{proof}
Recall that $\lambda_{[e]}$ denotes the $e$-core of $\lambda$. We have
$n_i(\lambda,Q)=n_i(\lambda_{[e]},Q)+w_e(\lambda)$ for each $i$.
Hence, from \eqref{eq:weightlevel1} we deduce that the weight of $|\lambda,Q\rangle$ is 
$-w_e(\lambda)\,\delta$ plus the weight of $|\lambda_{[e]},Q\rangle$.

\smallskip
Now assume that $\lambda$ is an $e$-core. Since $\tau_e(\lambda,0)=
(\emptyset,\lambda_{[e]})$, we have $\tau_e(\lambda,d)=(\emptyset,s)$ for some
tuple $s\in\bbZ^\scrI(d)$.
Hence, the weight of the element $|\lambda,Q\rangle$ in $\bfF(d)$ equals
\begin{align}\label{formA}\Lambda_1+\pi_s-\Delta(s,1)\,\delta=t_{\pi_s-\Lambda_Q^\cl}(\Lambda_Q)\end{align}
by Uglov's formulas, see e.g.,  \cite[prop.~3.7]{Y05}. The discussion above implies
part (a). Part (b) is a direct consequence of (a) since $W$ acts trivially on $\delta$.
\end{proof}

In the particular case where the charge is zero (forcing $Q$ to be $1$), then
$s = \lambda_{[e]}$ is the $e$-core of $\lambda$, and the weight of $|\lambda,1\rangle$
is given by
\begin{align}\label{formB} \mathrm{wt}(|\lambda,1\rangle)  = t_{\pi_{\lambda_{[e]}}}(\Lambda_1)-w_e(\lambda)
   \,\delta.\end{align}
Therefore weight spaces are parametrized by pairs $(\nu,w)$ where $\nu$ is an $e$-core
and $w$ is a non-negative integer. The basis element $|\lambda,1\rangle$ is in the weight
space corresponding to $(\lambda_{[e]},w_e(\lambda))$. We will see later that these weight
spaces correspond to the unipotent $\ell$-blocks of finite unitary groups $\mathrm{GU}_n(q)$ 
when $e$ is the order of $-q$ modulo $\ell$.

\subsubsection{The crystal of the Fock space}\label{subsec:crystalfock}
We explain here how to associate an abstract crystal to a charged Fock space $\bfF(s)$. 
By Proposition \ref{prop:tensorfock}, we can assume that $\scrI$ is either cyclic or of
type $A_\infty$. Then, the abstract crystal of $\bfF(s)$ is the abstract crystal associated with Uglov's canonical basis of $\bfF(s)$. We assume that the reader is familiar with \cite{U}. Another good reference is \cite{Y05}. 

\smallskip

When $\scrI$ has type $A_\infty$, Uglov's bases coincide with the standard monomial
basis and the discussion is trivial in that case. We will therefore assume that
$v$ has finite order $e$ and $\scrI = v^\bbZ$, so that $\scrI$ has type $A_{e-1}^{(1)}$.
Let $u$ be a formal parameter and $A = \bbC[u,u^{-1}]$. 
The $\frakg$-module $\bfF(Q)$ admits a quantum deformation $\bfF_{\!u}(s)$ with an
$A$-lattice $\bfF_{\!A}(s)$,  which is a free $A$-module with basis
$\{|\mu,s\rangle\,|\,\mu\in\scrP^l\}$.
It is equipped with an integrable representation of $U_A(\frakg)$ which is
given by the formulas (29), (35), (36) in \cite{U}.
Note that the action of the Chevalley generators $e_i$ and $f_i$ 
depends on the choice of the charge $s$. The representation of $\frakg$ on $\bfF(s)$ given in
\S \ref{subsec:xgrading} is recovered by specializing the parameter $u$ to $1$.

\smallskip

Uglov has constructed a remarkable $A$-basis of $\bfF_{\!A}(s)$ in \cite[p.~283]{U}
$$\bfB^+_u(s)=\{\bfb_u^+(\mu,s)\,|\,\mu\in\scrP^l\}.$$
It depends on $s$ and it is a lower global basis for the representation of $U_u(\frakg)$ on $\bfF_{\!u}(s)$.
In \cite{U}, this basis is denoted by the symbol 
$$\calG^+(s_l)=\{\calG^+(\lambda_l,s_l)\,|\,\lambda_l\in\scrP^l\}.$$
In order to match Uglov's parameters with ours, we set $q,l,n,s_l=u,l,e,s$ in the definition of
$\calG^+(s_l)$ to get our basis $\bfB_u^+(s)$.

\smallskip

Next, we consider the pairing $(\bullet,\bullet)$ on $\bfF_{\!u}(s)$ defined in \cite[sec.~4.3]{Y05}.
It is a modified version of the pairing $\langle\bullet,\bullet\rangle$ introduced in \S \ref{subsec:Fock}.
Then, let $\bfB_u^\vee(s)=\{\bfb_u^\vee(\mu,s)\,|\,\mu\in\scrP^l\}$ be the $\bbC(u)$-basis of
$\bfF_{\!u}(s)$ dual to $\bfB_u^+(s)$ relatively to the bilinear form $(\bullet,\bullet)$.
By Kashiwara's theory of global bases, we deduce that
$\bfB_u^\vee(s)$ is an upper global basis of $\bfF_{\!u}(s),$ compare \cite[lem.~4.13]{Y05}.
Let $\bfB^\vee(s)=\{\bfb^\vee(\mu,s)\,|\,\mu\in\scrP^l\}$ be the specialization at $u=1$ of
$\bfB_u^\vee(s)$, with the obvious labeling of its elements. 
It is a perfect basis of the usual Fock space $\bfF(s) =  \bfF_u(s)|_{u=1}$ by Proposition \ref{prop:upper-perfect}.

\smallskip

Next, we equip the set of $l$-partitions
with the abstract crystal structure $B(s)=\big(\scrP^l,\widetilde{e_i},\widetilde{f_i}\big)$ defined in \cite{JMMO}.
See \cite{FLOTW} for a reformulation closer to our notations.
Let $B(s)=\{b(\mu,s)\,|\,\mu\in\scrP^l\}$ be the obvious labeling.
The operators $\widetilde {e_i},\widetilde {f_i}$ are described in a combinatorial way:
we have $\widetilde{ f_i}(b(\mu,s))=b(\gamma,s)$ 
if and only if $\gamma$ is obtained from $\mu$ by adding a 
\emph{good $i$-node}. The definition of a good $i$-node depends on the charge $s$.
See \cite[sec.~3.4.2]{Y05} for more details.

\smallskip

We will need the following well-known result.

\begin{proposition}\label{prop:Uglov}
The map $\bfb^\vee(\mu,s) \in \bfB^\vee(s)\longmapsto b(\mu,s) \in B(s)$ is a crystal isomorphism.
\end{proposition}

\begin{proof}
The proposition follows from \cite{JMMO}.
More precisely, it is proved there that the formulae (27), (33), (34) in \cite{U} for the action of the quantum group on $\bfF_{\!A}(s)$
imply that the pair formed by $\calL(s)=\bigoplus_{\mu\in\scrP^l}\bbC[u]|\mu,s\rangle$ 
and the basis $\{|\mu,s\rangle\,\text{mod}\, u\,\calL(s)\,|\,\mu\in\scrP^l\}$ of $\calL(s)/u\,\calL(s)$
is a lower crystal basis of $\bfF_{\!u}(s)$ and that the assignment 
$|\mu,s\rangle\mapsto b(\mu,s)$ is a crystal isomorphism onto $B(s)$. 
In other words, the abstract crystal associated with the lower global basis $\bfB^+_A(s)$ is canonically isomorphic to $B(s)$, \emph{i.e.}, 
the map $\bfb_u^+(\mu,s) \in \bfB^+_A(s)
\longmapsto  b(\mu,s) \in B(s)$ realizes this isomorphism.
To conclude, we use the canonical bijections $\bfB^+_u(s)\to\bfB_u^\vee(s) \to \bfB^\vee(s)$ whose composition is an isomorphim of crystals as well.
\end{proof}

\section{Unipotent representations}
\label{part:modular}
In this section we record standard results on unipotent representations of finite
reductive groups in non-defining characteristic. A good reference is \cite{CE}.  

\subsection{Basics}\label{sec:basics}
By an \emph{$\ell$-modular system} we will mean a triple $(K,\scrO,\k)$ where 
$K$ is a field of characteristic zero, $\scrO$ is a complete discrete valuation
ring with fraction field $K$, and $\k$ is the residue field of $\scrO$ with
$\text{char}(\k)=\ell$. When working with representations of a finite group $\Gamma$,
we will always assume that $(K,\scrO,\k)$ is a \emph{splitting $\ell$-modular system
for $\Gamma$}, which means that $K$ and $\k$ are splitting fields for all subgroups
of $\Gamma$. When $\Gamma$ comes from an algebraic group in characteristic $p$,
we will in addition assume that $\ell \neq p$. This case is usually referred to as
the \emph{non-defining characteristic} case.

\smallskip

Let $R$ be any commutative domain (with 1) and $\Gamma$ be a finite group.
We will assume that $p$ is invertible in $R$ and $\scrO$.
Let $R\Gamma$ denote the group ring of $\Gamma$ over $R.$ 
For any subset $S\subseteq\Gamma$ such that $|S|$ is invertible in $R$, let $e_S$ 
be the idempotent $e_S=|S|^{-1}\sum_{g\in S}g$ in $R\Gamma$. If $R$ is not a field,
an $RG$-module which is free as an $R$-module will be called an \emph{$RG$-lattice}.

\smallskip

The $R$-module of class functions $\Gamma\to R$ is denoted by 
$R\,\Irr(K\Gamma)$ or $R\,\Irr(\Gamma)$. If $R$ is a field, this vector
space is endowed with the canonical scalar product, $\langle-,-\rangle_\Gamma$,
for which the set of irreducible characters
$\Irr(K\Gamma)$ of $K\Gamma$ is an orthonormal basis.

\subsection{Unipotent $KG$-modules}
Let $\bfG$ be a connected reductive group over $\overline\bbF_q$ 
with a Frobenius endomorphism $F : \bfG\to\bfG$.
Fix a parabolic subgroup $\bfP$ of $\bfG$ and an $F$-stable Levi complement
$\bfL$ of $\bfP$. We \emph{do not} assume $\bfP$ to be $F$-stable.
Write $L=\bfL^F$ and $G=\bfG^F$.

\smallskip

Let $R^{\bfG}_{\bfL\subset\bfP}$ and ${}^*\!R^{\bfG}_{\bfL\subset\bfP}$ 
denote respectively the Lusztig induction and restriction maps
from $\bbZ\,\Irr(KL)$ to $\bbZ\,\Irr(KG)$. 
We will assume that the Mackey formula holds for $R^{\bfG}_{\bfL\subset\bfP}$ and
${}^*\!R^{\bfG}_{\bfL\subset\bfP}$, which we know for the groups we will focus on 
later (see \cite{BM11} for more details). 
Under this condition, the Lusztig induction and restriction do not depend on the 
choice of the parabolic subgroup $\bfP$, see \cite[chap.~6]{DM}. We abbreviate
$R^{\bfG}_{\bfL\subset\bfP}=R^{\bfG}_{\bfL}$ and ${}^*\!R^{\bfG}_{\bfL\subset\bfP}={}^*\!R^{\bfG}_{\bfL}$.

\smallskip

Let $\bfT$ be an $F$-stable maximal torus of $\bfG$ and let $\bfN$ be the normalizer of $\bfT$ in $\bfG$.
Fix an $F$-stable Borel subgroup $\bfB$ of $\bfG$ containing $\bfT$.
Write $B=\bfB^F$, $T=\bfT^F$ and $N=\bfN^F$.
The groups $\bfB$, $\bfN$ form a reductive $BN$-pair of $\bfG$ with Weyl group 
$\bfW=\bfW_\bfG$ given by $\bfW=\bfN/\bfT$.
Since $\bfB$, $\bfN$ are stable by $F$ and $\bfG$ is connected, 
the finite groups $B$, $N$ form a split BN-pair of $G$ whose Weyl group
$W=W(T)$ is given by $W=\bfW^F=N/T$ (see \cite[sec.~4]{Geck} for more details).

\smallskip

The $G$-conjugacy classes of
$F$-stable maximal tori of $\bfG$ are parametrized by the $F$-conjugacy classes in $\bfW$.
For each $w\in\bfW$ let  $\bfT_{w}$ be
an $F$-stable maximal tori in the $G$-conjugacy class parametrized by $w$. 
Under conjugation by some element of $\bfG$, the pair $(\bfT_{w}, F)$ is identified with the pair $(\bfT, wF)$. 
In particular, we have $T_w\simeq\bfT^{wF}$ and $W(T_{w})\simeq\bfW^{wF}$.
The virtual characters $R^\bfG_{\bfT_{w}}(1)$ obtained by induction of the 
trivial representation of the tori $T_w$ are called the \emph{Deligne-Lusztig
characters}. They satisfy the following orthogonality relations:
$$\langle R^\bfG_{\bfT_{w}}(1)\,,\,R^\bfG_{\bfT_{w'}}(1)\rangle_G=
\begin{cases}|\bfW^{wF}|&\text{if\ $w$\ and\ $w'$\ are\ $F$-conjugate\ in\ $\bfW$},\\
0&\text{otherwise.}
\end{cases}$$

\begin{definition}
An irreducible $KG$-module is \emph{unipotent} if its character, say $\chi$, occurs 
as a constituent of a Deligne-Lusztig character $R^\bfG_{\bfT_{w}}(1)$ for some element
$w\in W$, \emph{i.e.}, if we have $\langle\chi,R^\bfG_{\bfT_{w}}(1)\rangle_G\neq 0$.
\end{definition}

We denote by $KG\umod$ the full subcategory of $KG\mod$ consisting of the modules
which are sums of irreducible unipotent modules. The objects of this category are
the \emph{unipotent $KG$-modules}.

\subsection{Unipotent $\k G$-modules and $\ell$-blocks}\label{sec:unipblocks}
As a result of the lifting of idempotents, the blocks of $\scrO G$ and
$\k G$ correspond by reduction. Both are usually called the \emph{$\ell$-blocks} of
$G$. For $R=\scrO$ or $\k$, any block $B$ of $R G$ is of the form $B=R G\cdot b,$
where $b$ is a central primitive idempotent of $R G$.
The unit $b$ of $B$ is called the \emph{block idempotent of  $B$}.
We will also call \emph{block of $R G\mod$ associated with $B$}
the Serre subcategory generated by the simple modules on which $b$ acts non-trivially.
The $\ell$-blocks of $G$ induce a partition of $\Irr(KG)$ such that the piece
associated with $B$ is the set of all irreducible characters $\chi$ of $KG$ with
$\chi(b)=\chi(1)$. If $\chi\in\Irr(KG)$, we will write $B(\chi)\subseteq \Irr(KG)$
for the piece containing $\chi$. When there is no risk of confusion, we will also
call $B(\chi)$ an $\ell$-block of $G$.

\begin{definition}\label{def:lblock}
An $\ell$-block of $\scrO G$ is \emph{unipotent} if it contains 
at least one unipotent $KG$-module.
A simple $\k G$-module is \emph{unipotent} if it lies in a unipotent block of $\k G$.
\end{definition}

We denote by $\k G\umod$ be the Serre subcategory of $\k G\mod$ generated by
the simple unipotent $\k G$-modules. It correspond to the sum of unipotent
blocks of $\k G\mod$. The \emph{unipotent $\k G$-modules} are by definition the objects
of this category.

\smallskip

Recall that $(K,\scrO,\k)$ is a splitting $\ell$-modular system. To this system
one can associate a decomposition map $d_{\scrO G}:[KG\mod]\to[\k G\mod]$.
From now on we will assume that the centre $\bfZ(\bfG)$ of $\bfG$ is connected.
Then, by \cite{H90}, a simple $\k G$-module is unipotent if and only if it is a
constituent of the $\ell$-reduction of a unipotent $K G$-module, see
also \cite{BM}. In other words, the classes of unipotent modules are 
exactly the image of unipotent characters through the decomposition map.
We will denote by $d_{\scrU}:[KG\umod]\to[\k G\umod]$ the restriction
of this map to unipotent characters.

\begin{proposition}[\cite{G93}, \cite{GH91}] \label{prop:isodecmap}
Assume $\ell$ is good for $\bfG$. Then the map
$d_{\scrU}$ is a linear isomorphism $[KG\umod]\simto[\k G\umod]$.
\qed
\end{proposition}

Given a positive integer $f$, let $\Phi_f$ be the $f$th cyclotomic polynomial.
A torus $T\subset G$ is a \emph{$\Phi_f$-torus} if its order is a power of $\Phi_f(q)$.
An $F$-stable Levi subgroup $\bfL\subseteq\bfG$ is \emph{$f$-split} 
if $L=\bfL^F$ is the centralizer in $G$ of a $\Phi_f$-torus.
A \emph{unipotent $f$-pair} is a pair $(\bfL,\chi)$ where $\bfL$ is an $f$-split
Levi subgroup and $\chi$ is an irreducible unipotent $KL$-module.
The pair $(\bfL,\chi)$ is \emph{$f$-cuspidal} if for every proper
$f$-split Levi subgroup $\bfM\subseteq\bfL$ we have ${}^*\!R^\bfL_\bfM(\chi)=0$.

\smallskip

Now assume that $f$ is the smallest positive integer such that $\ell$ divides
$q^f-1$. In other words, $f$ is the order of the class of $q$ in $\k$.
Under the assumption that $\ell$ is good, unipotent $\ell$-blocks correspond
to $f$-cuspidal $f$-pairs (see for example \cite[thm. 22.9]{CE}).

\begin{proposition}\label{prop:fcuspidalpairs}
Assume $\ell$ is good for $\bfG$, and $\ell \neq 3$ if $\bfG$ has a constituent
of type ${}^3 D_4$. 
Then there is a bijection between the $G$-conjugacy classes of unipotent
$f$-cuspidal $f$-pairs and the set of unipotent
$\ell$-blocks of $G$ which takes the class of $(\bfL,\chi)$ to the $\ell$-block
$B_{\bfL,\chi}$ such that the irreducible unipotent characters in $B_{\bfL,\chi}$
are exactly the irreducible constituents of $R^\bfG_\bfL(\chi).$
\qed
\end{proposition}

\subsection{Harish-Chandra series}
Assume now that the parabolic subgroup $\bfP\subseteq\bfG$ is $F$-stable.
In that case the group $L$ is $G$-conjugate to a standard Levi subgroup of $G$.
Let $R_L^{G}$ and ${}^*\!R^{G}_L$ be the corresponding Harish-Chandra induction and restriction functors from $RL\mod$ to $RG\mod$.  
Let $P=\bfP^F$ and $U=\bfU^F$, where $\bfU\subset \bfP$ is the unipotent radical of $\bfP$. 
Notice that the Harish-Chandra induction is the special case of Lusztig induction for 1-split Levi subgroups.

\smallskip

The order of $U$ is a power of $q$, hence it is invertible in $R$.
Thus, for all $M\in RL\mod$, $N\in RG\mod$ we have 
  $$R_L^{G}(M)=RG_n\cdot e_U\otimes_{RL}M \quad \text{and} \quad
  {}^*\!R_L^{G}(N)=e_U\cdot RG_n\otimes_{RG_n}N.$$
We will say that the functors $R_L^{G}$ and ${}^*\!R^{G}_L$ are \emph{represented} by the
$(RG,RL)$-bimodule $RG\cdot e_U$ and the $(RL,RG)$-bimodule $e_U\cdot RG$ respectively.

\smallskip

Here are some well-known basic properties of the functors $R_L^{G}$, ${}^*\!R^{G}_L$,
see for example \cite[prop. 1.5]{CE}. 
\begin{itemize}[leftmargin=8mm]
  \item[(a)] $R_L^{G}$, ${}^*\!R^{G}_L$ do not depend on $P$,
  \item[(b)] $R_L^{G}$, ${}^*\!R^{G}_L$ are exact and left and right adjoint
  to one another,
  \item[(c)] if $L\subseteq M\subseteq G$ there are isomorphisms of
  functors $R_L^{G}=R_M^{G}R_L^{M}$ and ${}^*\!R_L^{G}={}^*\!R_M^{G}{}^*\!R_M^{G}$.
\end{itemize}

Let $R=K$ or $\k$. An irreducible $RG$-module $E$ is \emph{cuspidal} if 
${}^*\!R^{G}_L(E)=0$ for all standard Levi subgroup $L\subsetneq G$.
A \emph{cuspidal pair of $RG$} is a pair $(L,E)$ where $L$ is as above and
$E\in\Irr(RL)$ is cuspidal. 
Since the group $L$ is uniquely recovered from $G$, $E$, from now on we may omit it from the notation.
Then, the set $\Irr(R G,E)\subseteq\Irr(R G)$ consisting of the constituents of the top of $R^G_L(E)$ 
is equal to the set of the constituents of the socle of $R^G_L(E)$ and
is called the \emph{Harish-Chandra series} of $(L,E)$.
The $R$-algebra $\scrH(RG,E)=\End_{RG}(R^{G}_L(E))^\op$ is
its \emph{ramified Hecke algebra}.
We have the following facts: 
\begin{itemize}[leftmargin=8mm]
\item[(d)] the Harish-Chandra series form a partition of $\Irr(R G)$,
\item[(e)] the functor $\frakF_{R^{G}_L(E)}$ yields a bijection $\Irr(R G,E)
\mathop{\longleftrightarrow}\limits^{1:1} \Irr(\scrH(RG,E)).$
\end{itemize}
See \cite{H93}, \cite{GHM96} and \cite[thm.~4.2.6, 4.2.9]{GJ} for details.

\smallskip

The following result is well-known. It follows from the 
fact that the $\ell$-rational series are stable by Harish-Chandra induction
(see for example \cite[thm. 10.3]{BoRo03} for a more general statement). 

\begin{proposition}\label{prop:stable} If $R=K$ or $\k$, then 
the Harish-Chandra induction and restriction functors preserve the category
of unipotent $RG$-modules. 
\end{proposition}

\section{Finite unitary groups}\label{part:unitarygroups}

This section is devoted to the construction of categorical actions on the
category of unipotent representations of finite unitary groups $\mathrm{GU}_n(q)$.
It contains the main results of this paper. 

\smallskip

Let $R$ be a commutative domain with unit. 
Under mild assumptions on $R$, we construct in \S\ref{sec:rep-datum} a representation
datum on the abelian category 
 $$ RG\mod = \bigoplus_{n \geqslant 0} \mathrm{GU}_n(q)\mod$$
given by Harish-Chandra induction and restriction. It consists of the adjoint pair 
$(E,F)$ of the functors themselves, together with natural transformations $X$ and $T$
of $F$ and $F^2$. The construction of the latter are similar to the case of
$\mathrm{GL}_n(q)$ given in \cite{CR}, and $X$ should be thought of as a Jucys-Murphy element,
whereas $T$ satisfies a Hecke relation with parameter $q^2$.
When $R$ is an extension of $\bbQ_\ell$ or $\bbF_\ell$, the categorical datum
restricts to the category $\scrU_R$ of unipotent representations of $RG\mod$. 
On this smaller category, the eigenvalues of $X$ are powers of $-q$, which hints
that the Lie algebra $\frakg$ that should act on $[\scrU_R]$ corresponds
to the quiver with vertices $(-q)^\bbZ$ and arrows given by multiplication by
$q^2$. 

\smallskip

We prove that the representation datum lifts indeed to a categorical action of $\frakg$
on $\scrU_R$. We start in \S\ref{sec:g-e} with the case where $\bbQ_\ell\subset R$. Then
$\frakg \simeq (\fraks\frakl_\bbZ)^{\oplus 2}$ and $[\scrU_R]$ is isomorphic to
a direct sum of level 2 Fock spaces, each of which corresponds to an ordinary
Harish-Chandra series. When $\bbF_\ell\subset R$, the category $\scrU_R$ is no
longer semisimple but we show in \S\ref{sec:g-e} a compatibility between weights for
the action of $\frakg$ and unipotent $\ell$-blocks which yields our second
categorification result. The situation depends on the parity of $e$, the order of
$-q$ modulo $\ell$. When $e$ is even (linear prime case), $\frakg$ is a subalgebra of
$(\widehat{\fraks\frakl_{e/2}})^{\oplus 2}$ and each ordinary Harish-Chandra series
categorifies a level 2 Fock space for $\frakg$. When $e$ is odd (unitary prime
case), $\frakg \simeq  \widehat{\fraks\frakl_\e}$ and weight spaces correspond
to unipotent $\ell$-blocks, which are now transverse to the ordinary Harish-Chandra series.

\smallskip

For studying the weight space decomposition of $[\scrU_R]$ as well as the action of 
the affine Weyl group we use the action of a bigger Lie algebra $\frakg_\circ$, which 
comes from Harish-Chandra induction and restriction for general linear groups. 
Going from linear to unitary groups by Ennola duality introduces signs for this action
(see Lemma \ref{lem:isofocks}). It is to be expected that the action of $\frakg_\circ$
on $[\scrU_R]$ lifts to an action by triangulated functors on $D^b(\scrU_R)$ coming
from Lusztig induction, although we will not use it.

\smallskip

We give two main applications of our categorical construction. In \S \ref{sec:broue}
we use Chuang-Rouquier's framework to produce derived equivalences between blocks, 
from which we deduce Brou\'e's abelian defect conjecture when $e$ is even.
Finally, we show in \S \ref{sec:isocrystals} that the crystal graph of the level 2 Fock
spaces that we categorified coincide with the Harish-Chandra branching graph.
This solves a recent conjecture of Gerber-Hiss-Jacon \cite{GHJ} and gives a combinatorial
way to compute the modular Harish-Chandra series and the parameters of the various
ramified Hecke algebras.

\subsection{Definition}\label{sec:def-unitary}
Fix a positive integer $n$.
We equip the reductive algebraic group $\bfG\bfL_n=\GL_n(\overline\bbF_q)$ with
the standard Frobenius map $F_q:\bfG\bfL_n\to\bfG\bfL_n$, $(a_{ij})\mapsto (a_{ij}^q)$
given by raising every coefficient to the $q$th power. The finite general linear group
$\mathrm{GL_n}(q)$ is given by the fixed points of $\bfG\bfL_n$ under $F_q$. 
In this section we will work with a twisted version of this group obtained by
twisting the Frobenius map.  Let $J_n$ be the
$n\times n$ matrix with entry $1$ in $(i,n-i+1)$ and zero elsewhere. We will often
write $J = J_n$ when there is no risk of confusion on the size of the matrices.
We define a new Frobenius map $F$ on $\bfG\bfL_n$, called the \emph{twisted} Frobenius
map, by setting $F=F_q\circ\alpha$ where $\alpha(g)=J\cdot{}^tg^{-1}\cdot J$ for each
$g\in \bfG\bfL_n$. The \emph{finite unitary group $G_n=\GU_n(q)$} is then given by
$$G_n=(\bfG\bfL_n)^F=\{g\in\bfG\bfL_n\,;\,F(g)=g\}.$$
Since $F^2=(F_q)^2$ we have $G_n\subset GL_n$, where we abbreviate $GL_n=\GL_n(q^2):=(\bfG\bfL_n)^{F^2}$. By convention we also define $G_0 = \{1\}$
to be the trivial group.

\smallskip
We equip $\bfG\bfL_n$ with the standard split BN-pair such that $\bfB$ is the
subgroup of upper triangular matrices and $\bfN$ is the subgroup of all monomial
matrices. Since $\bfB$, $\bfN$ are stable by $F$ and $\bfG\bfL_n$ is connected, the
groups $B=\bfB^F$, $N=\bfN^F$ form a split BN-pair of the finite group $G_n$.
Let $\bfT$ be the diagonal torus in $\bfG\bfL_n$ and $T=\bfT^F$.
Let $\bfW=\bfW_n$ be the Weyl group of $\bfG\bfL_n$ and $W=W_n$ be the Weyl group of
$G_n$. We have $\bfW\simeq\frakS_n$, and $F$ induces on $\bfW$ the automorphism
given by conjugation with the longest element $w_0$. We will embed $\bfW$ in $G_n$
using permutation matrices, so that $w_0$ corresponds to $J$. We have $W=\bfW^F=
C_\bfW(w_0)$. It is a Weyl group of type $B_m$ if $n=2m$ or $2m+1$.

\smallskip

Let $\varepsilon_1,\dots,\varepsilon_n$ be the characters of $\bfT$ such that
$t=\diag(\varepsilon_1(t), \varepsilon_2(t), \dots, \varepsilon_n(t))$.
The roots (resp. simple roots) of $\bfG\bfL_n$ are given by $\{\varepsilon_i-\varepsilon_j \, | \, i \neq j\}$ (resp. $\{\varepsilon_r-\varepsilon_{r+1}\}$).
Let $s_r=(r,r+1)$ be the simple reflection in $\bfW$ associated with the simple root
$\alpha_r = \varepsilon_r-\varepsilon_{r+1}$. The action of $F$ on the roots
induces an automorphism $\sigma$ of the Dynkin diagram of $\bfG\bfL_n$ such that
 $$ F \circ \alpha_r^\vee= q \sigma (\alpha_r)^\vee\quad \text{with} \quad \sigma(\alpha_r)=\alpha_{n-r}.$$
For each root $\alpha$, let $\bfU_\alpha$ and $\alpha^\vee\in\Hom(\bbG_m,\bfG\bfL_n)$
be the corresponding root subgroup and cocharacter. We also choose an
isomorphism $u_\alpha \, : \, \mathbb{G}_a \simto \bfU_\alpha$ such that
$F(u_{\sigma(\alpha)}(t))=u_{\alpha}(-t^q)$. Note that a one-parameter subgroup
of $T$ has either $(q^2-1)$ or $(q+1)$ elements, and a root subgroup of $G_n$ has either $q$ or $q^2$ elements.

\smallskip

The standard Levi subgroups $L_{r,m}$ of $G_n$ are parametrized by pairs
$(r,m)$ where $r$ is a non-negative integer and $m=(m_1,m_2,\dots,m_s)$ is a tuple of
positive integers such that $n=r+2\sum_{u=1}^sm_u$. The group $L_{r,m}$ consists
of all matrices of $\bfG\bfL_n$ which belong to $G_n$ and are of block-diagonal type 
$$\prod_{u=1}^s\bfG\bfL_{m_u}\times\bfG\bfL_r\times \prod_{u=s}^1\bfG\bfL_{m_u}.$$ 
Consequently we have a group isomorphism $L_{r,m}\simeq G_r\times\prod_uGL_{m_u}.$
If $m$ is a positive integer we abbreviate $L_{r,1^m}=L_{r,(1^m)}$ and
$L_{r,m}=L_{r,(m)}$.

\subsection{The representation datum on $RG$-mod}
\label{sec:rep-datum}
Let $R$ be a commutative domain with unit. We assume that $q(q^2-1)$ in invertible
in $R$. Using parabolic induction and restriction, we show in this section how
to construct a representation datum on 
 $$ RG\mod=\bigoplus_{n\in\bbN} RG_n\mod.$$

Fix a positive integer $n$. Parabolic (or Harish-Chandra) induction provides functors
between $L\mod$ and $G_n\mod$ for any standard Levi subgroup $L = L_{r,m}
\subset G_n$. Since we want functors between $G_{r}\mod$ and $G_n\mod$
we will consider a slight variation of the usual parabolic induction.

\smallskip 

Let $0 \leqslant r < n$. We denote by $V_r$ the unipotent radical of the standard
parabolic  subgroup $P_{r,1} \subset G_{r+2}$ with Levi complement $L_{r,1}$.
Let $U_{r}\subset G_{r+2}$ be the subgroup given by
\begin{align}\label{U}
  U_{r}=V_{r}\rtimes\bbF^\times_{q^2}= 
  \begin{pmatrix}*&*&\cdots&\cdots&*\\
	&1&&&\vdots\\
	&&\ddots&&\vdots\\
	&&&1&*\\
	&&&&*
  \end{pmatrix}
\end{align}
so that $P_{r,1} =V_{r} \rtimes  L_{r,1}  \simeq U_r \rtimes G_r$.
If $n-r$ even, we set $U_{n,r}=U_{n-2}\rtimes\dots\rtimes U_{r}$ and we define
$e_{n,r}=e_{U_{n,r}}$ to be the corresponding idempotent. In particular,
we have $e_{r+2,r}=e_{U_{r}}$. We embed $G_r$ into the Levi subgroup
$L_{r,1^m}=G_r\times GL_1^{m}$ in the obvious way. This yields an embedding
$G_r\subset G_n$ and functors
$$\begin{aligned}
  F_{n,r}=& \, RG_n\cdot e_{{n,r}}\otimes_{RG_r}-\,:\,RG_r\mod\longrightarrow RG_n\mod,\\
  E_{r,n}=& \, e_{{n,r}}\cdot RG_n\otimes_{RG_n}-\,:\,RG_n\mod\longrightarrow RG_r\mod.
\end{aligned}$$
Note that $F_{n,r}$ can be seen as the composition of the inflation $G_r\mod
\longrightarrow L_{r,1^m}\mod$ with the parabolic induction from $L_{r,1^m}$ to $G_n$. 

\smallskip

An endomorphism of the functor $F_{n,r}$ can be represented by an $(RG_n,RG_r)$-bimodule
endomorphism of $RG_n\cdot e_{{n,r}}$, or equivalently by an element of 
$e_{{n,r}}\cdot RG_n\cdot e_{{n,r}}$ centralizing $RG_m$. Thus, the elements
$$ X_{r+2,r}= (-q)^r e_{{r+2,r}} (1,r+2)\, e_{{r+2,r}},\qquad
  T_{r+4,r}= q^2 e_{{r+4,r}}(1,2)(r+3,r+4)\,e_{{r+4,r}}$$
define respectively natural transformations of the functors $F_{r+2,r}$ and $F_{r+4,r}$.
We set 
$$F=\bigoplus_{r\geqslant 0}F_{r+2,r},\qquad 
  X=\bigoplus_{r\geqslant 0}X_{r+2,r},\qquad
  T=\bigoplus_{r\geqslant 0}T_{r+4,r}.$$

\begin{proposition}\label{prop:rep-datum}
The endomorphisms $X\in\End(F)$ and $T\in\End(F^2)$ satisfy the following relations: 
\begin{itemize}[leftmargin=8mm]
  \item[$\mathrm{(a)}$] $1_FT\circ T1_F\circ 1_FT=T1_F\circ 1_FT\circ T1_F$,
  \item[$\mathrm{(b)}$] $(T+1_{F^2})\circ(T-q^21_{F^2})=0$,
  \item[$\mathrm{(c)}$] $T\circ(1_FX)\circ T=q^2X1_F$.
\end{itemize}
\end{proposition}

\begin{proof}
The relation (a) comes from the usual braid relations. For (b), we compute
\begin{align*}
  (T_{r+4,r})^2&=q^4 e_{{r+4,r}}(1,2)(r+3,r+4)\,e_{{r+4,r}}(1,2)(r+3,r+4)\,e_{{r+4,r}}\\
  &=q^4 e_{{r+4,r}}\,e_V\,e_{{r+4,r}},
\end{align*}
where $V\subset G_{r+4}$ is the subgroup consisting of the matrices with diagonal
entries equal to 1 and off-diagonal entries equal to zero, except for the entries
$(2,1)$ and $(r+4,r+3)$, \emph{i.e.}
  $$V=\begin{pmatrix}1&&&&\\ *&1&&&\\&&\ddots&&\\&&&1&\\&&&*&1\end{pmatrix}.$$
The group $V$ is the root subgroup $u_{-\alpha}(\bbF_{q^2})$ of $G_{r+4}$ associated
with the negative root say $-\alpha$. The corresponding simple reflexion $s_\alpha$ 
is given by the permutation $(1,2)(r+3,r+4).$ Let $B_\alpha$ be the (finite) Borel
subgroup of upper triangular matrices in the copy of $GL_2$ in $G_{r+4}$ associated
with $\alpha$. We have $B_\alpha u_{-\alpha}(t)\,B_\alpha=B_\alpha s_\alpha\,B_\alpha$
if and only if $t\neq 0$, so that 
$$q^2e_{B_\alpha}e_{U_{-\alpha}}e_{B_\alpha}= e_{B_\alpha}+
  (q^2-1)\,e_{B_\alpha}s_\alpha e_{B_\alpha}$$ 
in $RGL_2$. Now, the image of $B_\alpha$ through the embedding $GL_2\subset G_{r+4}$
lies in $U_{r+4,r}$. As a consequence, 
$$U_{r+4,r}e_V\,U_{r+4,r}=q^{-2}e_{{r+4,r}}+(1-q^{-2})e_{{r+4,r}}(1,2)(r+3,r+4)
\,e_{{r+4,r}},$$
which yields the expected formula
$$(T_{r+4,r})^2=(q^2-1)\,T_{r-4}+q^2\,e_{{r+4,r}}.$$

Finally, to prove (c) we must compute $T_{r+4,r}X_{r+2,r}T_{r+4,r}$. Using the group $V$ introduced above, it equals
$$(-q)^{r+4} e_{{r+4,r}}e_V(1,2)(r+3,r+4)(2,r+3)(1,2)(r+3,r+4)\,e_Ve_{{r+4,r}},$$
which simplifies to $(-q)^{r+4} e_{{r+4,r}}e_V(1,r+4)\,e_Ve_{{r+4,r}}.$
Let $V'=^{(1,r+4)}\!V\subset U_{r+4,r}$. The only off-diagonal entries which are
non-zero in $V'$ are the entries $(2,r+4)$, $(1,r+3)$ which lie at the top right corner
$$V'=\begin{pmatrix}1&&&*&\\ &1&&&*\\&&\ddots&&\\&&&1&\\&&&&1\end{pmatrix}.$$
By Chevalley's commutator formula, we have $[V',V]\subset U_{r+2,r}$.
This proves that
$$e_V(1,r+4)\,e_Ve_{{r+4,r}}=(1,r+4)\,e_Ve_{{r+4,r}}.$$
Finally, by moving $e_V$ to the left, we obtain 
$$T_{r+4,r}X_{r+2,r}T_{r+4,r}=(-q)^{r+4}e_{{r+4,r}}(1,r+4)e_{{r+4,r}}=q^2X_{r+4,r+2},$$
from which (c) follows.
\end{proof}

\begin{corollary}\label{cor:pre-cat}
  The tuple $(E,F,X,T)$ defines a representation datum on $RG\mod =
  \bigoplus_{n\in\bbN} RG_n\mod$. \qed
\end{corollary}

The next step will be to show that this categorical datum yields
a $\frakg$-representation on the full subcategory of unipotent modules
of $RG\mod$ (see Theorems \ref{thm:char0} and \ref{thm:charl}). 

\begin{remark}
The reader might argue that we did not check that $X$ was invertible. 
Nevertheless, we shall only be working with unipotent modules over
a field, in which case the eigenvalues of $X$ are powers of $-q$, hence
non-zero (see Theorem \ref{thm:cat}). This will ensure that the
restriction of $X$ to this category is indeed invertible. 
\end{remark}

\subsection{The categories of unipotent modules $\scrU_K$ and $\scrU_\k$}
\label{sec:uK-and-uk}
From now on, we fix a prime number $\ell$ such that $\ell \nmid q\,(q^2-1)$,
and an $\ell$-modular system $(K,\scrO,\k)$ with $\bbQ_\ell\subset K$. We will
assume that the modular system is large enough, which means that $KG_n$ and
$kG_n$ are split for all $n \geqslant 0$. One can take for example
$(K,\scrO,\k) = (\overline{\bbQ}_\ell, \overline{\bbZ}_\ell,\overline{\bbF}_\ell)$.

\smallskip

Throughout the following sections, we will denote by $d$, $e$ and $f$ the order
of $q^2$, $-q$ and $q$ in $\k$. In particular $e \neq 1,2$. If $e$ is odd,
then $d=e$ and $f=2e$; if $e$ is even, then $d = e/2$ and either $f= e/2$ if
$e \equiv 2$ mod $4$ or $f=e$ if $e \equiv 0$ mod $4$. 

\subsubsection{The category $\scrU_K$}\label{subsec:uK}
Fix a positive integer $n$. By \cite{LS77}, the irreducible unipotent $KG_n$-modules
are labelled by partitions of $n$. Their character can be directly constructed from
the Deligne-Lusztig characters. Namely, for each $w\in\frakS_n$, fix an $F$-stable
maximal torus $\bfT_{w}\subset\bfG\bfL_n$ in the $G_n$-conjugacy class parametrized
by $w$, with the convention that $\bfT_1 = \bfT$. Then the class function
\begin{align}\label{unip.char}  \chi_\lambda=|\frakS_n|^{-1}\sum_{w\in\frakS_n}\phi_\lambda(ww_0)
  R_{\bfT_{w}}^{\bfG\bfL_n}  (1)\in\bbZ\Irr(KG_n)
\end{align}
is, up to a sign, an irreducible unipotent character. 
For each $\lambda$ we choose a corresponding irreducible $KG_n$-module $E_\lambda$.
By abuse of notation we will still denote by $E_\lambda$ its isomorphism class. 

\smallskip

Recall that $G_0=\{1\}$. We define the category of unipotent $KG$-modules by  
  $$\scrU_K=\bigoplus_{n\in\bbN}K G_n\umod.$$
This category is abelian semisimple. 
From the previous discussion we have $\Irr(\scrU_K)=\{E_\lambda\}_{\lambda\in\scrP}$,
where by convention $\Irr(KG_0)=\{E_\emptyset\}$.

\subsubsection{The category $\scrU_\k$}\label{subsec:uk}
Using the $\ell$-modular system we have decomposition maps $d_{\scrO G_n}$
which by Proposition \ref{prop:isodecmap} restrict to linear isomorphisms 
 $$ d_{\scrO G_n}:[KG_n\umod]\simto[\k G_n\umod]. $$
For unitary groups, this map is actually unitriangular with respect 
to the basis of irreducible modules and the dominance order.

\begin{proposition}[\cite{G91}]\label{prop:unitriangularunip} 
  There is a unique labelling $\{D_\lambda\}_{\lambda\in\scrP_n}$ of the unipotent
  simple $\k G_n$-modules such that
    $$d_{\scrO G_n}([E_\lambda]) \in \ [D_\lambda] + \sum_{\nu > \lambda} \bbZ[D_\nu]$$
  where $>$ is the dominance order on partitions of $n$. 
  \qed
\end{proposition}

We define the category of unipotent $\k G$-modules by  
  $$\scrU_\k=\bigoplus_{n\in\bbN}\k G_n\umod.$$
This is an abelian category which is not semi-simple. The isomorphism classes of
simple objects are $\Irr(\scrU_\k)=\{D_\lambda\}_{\lambda\in\scrP}$ and the
decomposition map yields a $\bbZ$-linear isomorphism $d_\scrU : [\scrU_K]\simto[\scrU_\k]$.
The following result is a consequence of the unitriangularity of this map. 

\begin{proposition}\label{prop:liftingstandards}
Given a partition $\lambda \vdash n$ there is a unique $\scrO G_n$-lattice
  $\widetilde E_\lambda$ such that 
  \begin{itemize}[leftmargin=8mm]
   \item[$\mathrm{(a)}$] $K\widetilde E_\lambda=E_\lambda$ as a $K G_n$-module,
   \item[$\mathrm{(b)}$] $V_\lambda:=\k\widetilde E_\lambda$ is an indecomposable $\k G_n$-module with
   top isomorphic to $D_\lambda$,
   \item[$\mathrm{(c)}$] $d_\scrU([E_\lambda])=[V_\lambda].$
  \end{itemize}
\end{proposition}

\begin{proof}
By \cite[ex.~6.16]{CR81}, any system of orthogonal idempotents of $\k G_n$ lifts to a system of orthogonal idempotents of $\scrO G_n$. Let $P_\lambda$ be the projective cover of $D_\lambda$,
and let $\widetilde P_\lambda$ be a lift of $P_\lambda$ to a projective indecomposable
$\scrO G_n$-lattice. By \cite{G91}, we have
  $$K \widetilde P_\lambda=E_\lambda\oplus \bigoplus_{\mu < \lambda} E_\mu^{\,\oplus\, 
  d_{\lambda,\mu}} \quad  \text{with} \ d_{\lambda,\mu}\in\bbN.$$
Let $e_\chi\in K G$ be the idempotent associated with an irreducible $K G$-character $\chi$.
Set $e_{\neq \lambda} = \sum_{\chi \not\simeq E_\lambda} e_{\chi}$. The $\scrO G_n$-submodule 
$N:=\widetilde P_\lambda \cap e_{\neq \lambda} K\widetilde P_\lambda$ of $\widetilde P_\lambda$
is pure. We define $\widetilde E_\lambda = \widetilde P_\lambda / N.$
It is an $\scrO G_n$-lattice. The module $\k\widetilde E_\lambda$ has a simple top equal to
$D_\lambda$, because $\top(\k\widetilde P_\lambda)=D_\lambda$. Furthermore, we have
$K \widetilde E_\lambda = K\widetilde P_\lambda / e_{\neq \lambda} K\widetilde P_\lambda=
E_\lambda.$

\smallskip
Now, let us concentrate on the unicity of $\widetilde E_\lambda$.
Let $\widetilde E$ be an $\scrO G_n$-lattice satisfying the properties (a), (b) above.
Let $\phi\,:\,\k\widetilde E\twoheadrightarrow D_\lambda$ be the obvious surjective map.
Since $\widetilde P_\lambda$ is projective, there is a morphism $\psi\,:\,\widetilde
P_\lambda\longrightarrow \widetilde E$ such that the following triangle commutes
$$\xymatrix{\k\widetilde P_\lambda\ar[r]\ar[dr]_-{\k\psi}&D_\lambda\\
&\k\widetilde E.\ar[u]_-{\phi}}$$
Since $\k\widetilde E$ has a simple top and the $\k$-algebra $\k G_n$ is finite dimensional,
the map $\phi$ is an essential epimorphism. Thus $\k\psi$ in onto, and so is $\psi$
by Nakayama's Lemma. We are left with proving that $\Ker(\psi)=N$.
Since the module $\widetilde E$ is free over $\scrO$, it is enough to check that $\Ker(K\psi)=e_{\neq \lambda} K\widetilde P_\lambda$.

\smallskip

To prove this, observe first that since $\widetilde P_\lambda$ is projective and $\widetilde E$ is free over $\scrO$,
the $\scrO$-module $\Hom_{\scrO G_n}(\widetilde P_\lambda,\widetilde E)$ is free and we have
$$\k\Hom_{\scrO G_n}(\widetilde P_\lambda,\widetilde E) =\Hom_{\k G_n}(\k\widetilde P_\lambda,\k\widetilde E) =\Hom_{\k G_n}(\k\widetilde P_\lambda,D_\lambda)=\k.$$
We deduce that the $\scrO$-module $\Hom_{\scrO G_n}(\widetilde P_\lambda,\widetilde E)$ is free of rank 1 and generated by the map $\psi$.
Consequently
$$K\psi =K\Hom_{\scrO G_n}(\widetilde P_\lambda,\widetilde E)
=\Hom_{KG_n}(K\widetilde P_\lambda,K\widetilde E)
=\Hom_{KG_n}(K\widetilde P_\lambda,E_\lambda)$$
and the claim is proved.
\end{proof}

\begin{remark}
In the case where $V_\lambda$ is simple, property (b) is superfluous and
$V_\lambda$ is automatically isomorphic to $D_\lambda$. This is for example
true when $E_\lambda$ is cuspidal (see \cite[thm.~6.10]{GHM1}).
\end{remark}

\subsubsection{Blocks of $\scrU_\k$}\label{subsec:blocksofuk}
A block of $\scrU_\k$ is an indecomposable summand of $\scrU_\k$. Therefore 
blocks of $\scrU_\k$ correspond to the unipotent blocks of $\k G_n$
where $n$ runs over $\bbN$. These were first obtained in \cite{FS82}, before the 
general classification was given in \cite{BMM} (see Proposition \ref{prop:fcuspidalpairs}). 

\smallskip

Recall that $e$ is the order of $-q$ modulo $\ell$. Given a partition $\lambda$,
we defined in \S\ref{subsec:l-core} its $e$-core $\lambda_{[e]}$, its $e$-quotient $\lambda^{[e]}$ and
its $e$-weight $w_e(\lambda) = |\lambda^{[e]}|$.

\begin{proposition}\label{prop:lblocks}
The map $E_\lambda \longmapsto (\lambda_{[e]},w_e(\lambda))$
yields a bijection between unipotent 
$\ell$-blocks and pairs $(s,w)$ where $s \in \bbZ^e(0)$ and $w \in \bbN$.
\qed
\end{proposition}

Recall that $\nu \mapsto \nu_{[e]}$ induces a bijection between $e$-cores
and $\bbZ^e(0)$. Given $\nu$ an $e$-core we will denote by $B_{\nu,w}$ or $B_{\nu_{[e]},w}$
the unipotent $\ell$-block containing all the unipotent characters $E_\lambda$ such that 
$\nu$ is the $e$-core of $\lambda$ and $w_e(\lambda) =w$. It is a block of $\k G_n$
with $n = |\nu|+ew$.

\smallskip

It also follows from the classification of blocks that when $e  < \ell$ and $e \neq 1$,
the defect group of $B_{\nu,w}$ is an elementary abelian $\ell$-group of rank $w$ (see
\cite{BMM}). In particular, when $w=0$ the module $D_\lambda$ is simple and projective, 
isomorphic to $V_\lambda$, and when $w =1$ the defect group of $B_{\nu,w}$ is a cyclic
group. The structure of such blocks was determined in \cite{FS}.

\subsubsection{The weak Harish-Chandra series}\label{subsec:defweakhc}
For this section we assume that $R$ is one of the fields $K$ or $\k$. 

\smallskip

Let $r,m \geqslant 0$ and $n = r+2m$. The inflation from $G_r$ to
$L_{r,1^m}$ yields a equivalence between $RG_r\umod$ and $RL_{r,1^m}\umod$, since
the Deligne-Lusztig varieties depend only on the semisimple type of the reductive
group. This equivalence intertwines the functors $E_{n,r}$, $F_{n,r}$ with the
parabolic restriction and induction ${}^*R_{L_{r,1^m}}^{G_n}$, $R_{L_{r,1^m}}^{G_n}$.
Therefore working with $\scrU_R$ and the functors $E$ and $F$ is the same as
working in the usual framework of unipotent representations and Harish-Chandra
induction/restriction from Levi subgroups. Note however that one does not
consider all the standard Levi subgroups, but only the ones that are conjugate
to $L_{r,1^m}$. Therefore one needs to consider a slight variation of the
usual Harish-Chandra theory.

\begin{definition}[\cite{GHJ}] Fix a non-negative integer $n$. 
 \begin{itemize}[leftmargin=8mm]
  \item[(a)] An $R G_n$-module $D$ is \emph{weakly cuspidal} if ${}^*\!R^{G_n}_L(D)=0$
   for any Levi subgroup $L\subsetneq G_n$ which is $G_n$-conjugated to a subgroup of
   the form $L_{r,1^m}$.
  \item[(b)] A \emph{weak cuspidal pair} of $R G_n$ is a pair which is $G_n$-conjugated
  to $(L_{r,1^m}\,,\, D)$ for some $n=r+2m$ and some weakly cuspidal irreducible $R G_r$-module
  $D$ which is viewed as a $\k L_{r,1^m}$-module by inflation.
  \item[(c)] The \emph{weak Harish-Chandra series} $\wIrr(R G_n, D)$ of $R G_n$
  determined by the weak cuspidal pair $(L,D)$ is the set of the constituents of the top
  of $R_L^{G_n}(D)$ and it coincides with the set of the constituents of its socle.
 \end{itemize}
\end{definition}

If $D$ is a unipotent $RG_n$-module, then $D$ is weakly cuspidal if and only if $E(D) =0$.
Moreover, if $D$ is irreducible, the weak Harish-Chandra series coincides with the set
of irreducible constituents in the top of $F^m D$. Therefore it makes sense to define
the weak Harish-Chandra of $\scrU_R$ by
  $$\wIrr(R G,D)=\bigsqcup_{n\geqslant r}\wIrr(R G_n,D).$$
It is the set of irreducible consituents in the top of $F^kD$ for some $k \geqslant 0$
(or equivalently in the socle). 
As in the case of the usual theory, the weak Harish-Chandra series of $\scrU_R$ form
a partition of $\Irr(\scrU_R)$ (see \cite[prop. 2.3]{GHJ} for a proof).

\begin{proposition}\label{prop:partitionwcusp}
  Assume that $R= K$ or $\k$. Then
      $$\Irr(\scrU_R) =\bigsqcup \wIrr(R G,D)$$
  where the sum runs over the set of isomorphism classes of weakly cuspidal unipotent
  irreducible modules $D$. 
  \qed
\end{proposition}

When $R=K$, weakly cuspidal unipotent modules coincide with cuspidal modules and were
determined in \cite{L77} (see also Corollary \ref{cor:cuspidal}). One of the main result
in this paper is the classification of the weakly cuspidal modules in the case where $R =\k$.

\subsection{The $\frakg_\infty$-representation on $\scrU_K$}
\label{sec:g-infty}
In this section we show how the categorical datum defined in \S \ref{sec:rep-datum} yields
a categorical representation on $\scrU_R$ in the case where $R=K$. 
This is achieved by translating in our framework the theory of Howlett-Lehrer
on endomorphism algebras of induced representations. 

\subsubsection{Action of $E$ and $F$}
Since every unipotent character is a linear combination of Deligne-Lusztig characters
$-$ we call such class function a \emph{uniform function} $-$ it is well-known how to
compute the action of $E$ and $F$ on the category $\scrU_K$. 

\begin{lemma}\label{lem:indresunip}
  Let $\lambda$ be a partition. Then 
  $$ [E] (E_\lambda) = \sum_{\mu} E_\mu \qquad \text{and} \quad [F] (E_\lambda) = \sum_{\mu'} E_{\mu'}$$
  where $\mu$ (resp. $\mu'$) runs over the partitions that are obtained from $\lambda$
  by removing (resp. adding) a $2$-hook. 
\end{lemma}

\begin{proof} 
By adjunction, it is enough to prove the formula for the restriction functor $E$. 
Let $n = |\lambda|$ and let $\bfL$ be the standard Levi with $\bfL^F = L = L_{n-2,1}$.
For $x \in G_n$, if ${}^x \bfT_w$ is a maximal torus of $\bfL$ then $w$ is $F$-conjugate
to an element $v$ of the Weyl group $\bfW_\bfL$ of $\bfL$ (which is isomorphic to $\frakS_{n-2}$).
In addition, the set $L\backslash \{x \in G_n \, |\, {}^xT_v \subset \bfL\}/T_v$ is in bijection
with $\bfW^{vF}/\bfW_\bfL^{vF}$. Thus the Mackey formula \cite[thm. 11.13]{DM} yields
\begin{align} \label{rest} {}^*R_{L}^{G_n} R_{\bfT_w}^{\bfG\bfL_n}(1) = \frac{|\bfW^{vF}|}{|\bfW_\bfL^{vF}|}
  R_{\bfT_v}^\bfL(1).
  \end{align}
Now if $w$ is not $F$-conjugate to an element of $\bfW_\bfL$, then
${}^*R_{L}^{G_n} R_{\bfT_w}^{\bfG\bfL_n}(1) = 0$. Using the fact that $w$ is $F$-conjugate to
$v$ if and only if $ww_0$ is conjugate to $vw_0$, 
we deduce from \eqref{unip.char}, \eqref{rest} the following formula for the restriction of the character $\chi_\lambda$:
  $${}^*R_{L}^{G_n}(\chi_\lambda) = \frac{1}{|\bfW_\bfL|} \sum_{v \in \bfW_\bfL}
  \phi_\lambda(vw_0) R_{\bfT_v}^\bfL(1).$$
To conclude, we write $\phi_\lambda(vw_0) = \phi_\lambda(vw_\bfL (1,n))$, where $w_\bfL$
is the longest element of $\bfW_\bfL$, and we apply the Murnaghan-Nakayama rule for
the $2$-cycle $(1,n)$. This gives a decomposition of ${}^*R_{L}^{G_n}(\chi_\lambda)$ in terms
of the $\chi_\mu$'s, where $\mu$ is obtained from $\lambda$ by removing a $2$-hook. 
Note that there is no need to worry about the signs in this formula: ${}^*R_{L}^{G_n}(\chi_\lambda)$
is a virtual character but $[E](E_\lambda)$ is a true character. 
\end{proof}

Since the only partitions from which one cannot remove any $2$-hook are 
the triangular partitions, we have the following immediate corollary. 

\begin{corollary}[\cite{L77}] \label{cor:cuspidal} 
 A unipotent $KG_n$ module $E_\lambda$ is weakly cuspidal if and only if there exists
 $t \geqslant 0$ such that $n=t(t+1)/2$ and $\lambda = \Delta_t := (t,t-1,t-2,\ldots,1)$.
\end{corollary}

We will set $E_t := E_{\Delta_t}$. Consequently, the weakly cuspidal pairs are,
up to conjugation, the pairs $(L_{r,1^m},E_t)$ with $r =t(t+1)/2$. The partition
into series is thus
  $$\Irr(\scrU_K) =\bigsqcup_{t \in \bbN} \wIrr(R G,E_t).$$
Note that it is proven in \cite{L77} that this is also the partition into
usual Harish-Chandra series. However we shall not use this fact.

\subsubsection{The Howlett-Lehrer isomorphism}\label{subsec:hl-iso}
Let $r,m$ be non-negative integers and $n = r+2m$. As mentioned in \S \ref{subsec:defweakhc},
the inflation
from $G_r$ to $L_{r,1^m}$ yields a equivalence between $KG_r\umod$ and
$KL_{r,1^m}\umod$ which intertwines the functor $F_{n,r}$ with the parabolic
induction $R_{L_{r,1^m}}^{G_n}$. In particular, we have a canonical isomorphism
  $$\scrH(KG_n,E_t):= \End_{KG_n}(F^m(E_t))^\op \simto 
  \End_{K G_n}(R^{G_n}_{L_{r,1^m}}(E_t))^\op.$$
Now recall from \S \ref{subsec:hecke} that to the categorical datum $(E,F,X,T)$ is attached
a map $\phi_{F^m}:\bfH_{K,m}^{q^2}\to\End(F^m)$. The evaluation of this map at the
module $E_t$ yields a $K$-algebra homomorphism
  $$\phi_{K,m}:\bfH^{q^2}_{K,m}\to\scrH(KG_n,E_t),\quad 
  X_k\mapsto X_k(E_t),\quad T_l\mapsto T_l(E_t).$$
By \cite{HL80}, $\scrH(KG_n,E_t)$ is isomorphic to a Hecke algebra $\bfH^{Q_t,\,q^2}_{K,m}$
of type $B_m$ with 
\begin{align}\label{Q}
  Q_t=\begin{cases}\big((-q)^{-1-t}\,,\,(-q)^t\big)&\ \text{if}\ t\ \text{is\ even},\\
  \big((-q)^t\,,\,(-q)^{-1-t}\big)&\ \text{if}\ t\ \text{is\ odd}.
  \end{cases}
\end{align}
We show that the previous map provides such an isomorphism.

\begin{theorem}\label{thm:cat}
Let $t,m \geqslant 0$ and $n = t(t+1)/2+2m$. Then the map $\phi_{K,m}$ factors through
a $K$-algebra isomorphism 
  $$\bfH^{Q_t,\,q^2}_{K,m}\mathop{\longrightarrow}\limits^\sim\scrH(KG_n,E_t).$$
\end{theorem}

\begin{proof} 
Write $Q_t=(Q_1,Q_2)$ and $X=X(E_t)$. We must check that the operator $X$ on $F(E_t)$
satisfies the relation 
  $$(X-(-q)^{-1-t})(X-(-q)^{t})=0.$$
Then the invertibility of the morphism $\bfH^{Q_t,\,q^2}_{K,m}\longrightarrow\scrH(KG_n,E_t)$
follows from the general theory of Howlett and Lehrer. In fact, it is shown in \cite{HL80}
that $X$ satifies the relation
\begin{equation} \label{eq:xsign} 
(X-\epsilon_t\,(-q)^{-1-t})(X-\epsilon_t\,(-q)^{t})=0. 
\end{equation}
for some $\epsilon_t = \pm 1$. Therefore we must show that $\epsilon_t=1$,
which we will do by induction on $t$. First observe that 
the eigenvalues of $X_{2,0}$ on $R_{G_0}^{G_2} (K)$ and of $X_{3,1}$ on
$R_{G_1}^{G_3} (K)$ are $1$, $(-q)^{-1}$ and $-q$, $(-q)^{-2}$ respectively. 
So the eigenvalues of $X(E_0)$, $X(E_1)$ on $F(E_0)$, $F(E_1)$ 
are powers of $-q$. Now fix $t > 1$ and assume that for all $t>s\geqslant 0$
the eigenvalues of $X(E_s)$ on $F (E_s)$ are powers of $-q$. We will
show that the eigenvalues of $X(E_t)$ are also powers of $-q$ using
the modular representation theory of unitary groups. 

\smallskip

Recall that $K$ is chosen with respect to an $\ell$-modular system $(K,\scrO,\k)$.
Since the parametrization of unipotent characters does not depend on $\ell$
and $(E,F,X,T)$ are defined over $\bbZ[1/q(q^2-1)]$, we can first choose a specific prime
number $\ell$ and prove that the eigenvalues of $X(E_t)$ are powers of $-q$ modulo $\ell$. 
We choose $\ell$ to be odd and such that
the order of $-q$ in $\k^\times$ is $e:=2t-1$. 

\smallskip

Fong-Srinivasan determined
the repartition of the unipotent characters into $\ell$-blocks, and the structure
of the blocks with cyclic defect (see \S\ref{subsec:blocksofuk} and \cite{FS}). 
It follows from their result that the module $E_t$ lies in a block
with cyclic defect, since only one $e$-hook can be removed from $\Delta_t$, the
corresponding partition.   

\smallskip

Let $r= t(t-1)/2$.
We refer to \S \ref{sec:combinatorics} for the combinatorics of partitions and $\beta$-sets.
Let $\beta = \beta_0(\Delta_t)$ be the set of $\beta$-numbers of the charged partition
$(\Delta_t,0)$. It equals $\beta = \{t,t-2,\ldots,-t\}\cup\{-t-1,-t-2,\ldots\}$. 
Only one $e$-hook can be removed from $\Delta_t$, therefore $E_t$ lies
in a block with cyclic defect. Removing the $e$-hook corresponds to replacing
$t$ by $t-e = -t+1$ in $\beta$ to obtain the set $\overline{\beta} = \{t-2,t-4,\ldots,-t+2\} 
\cup \{-t+1,-t,\ldots\}$. It is exactly the set of $\beta$-numbers of the charged
partition $(\Delta_{t-2},0)$. 
From now on we will assume that all the partitions are charged with the charge $0$. 
Now, the unipotent characters in the block of $E_t$ correspond to the partitions
obtained by adding a $e$-hook to $\Delta_{t-2}$. In terms of $\beta$-sets, they 
are obtained by replacing $x$ by $x+e$ in $\overline{\beta}$, with
$x \in \{-3t+4,-3t+6,\ldots, t-4,t-2\}\cup\{-t+1\}$. We denote by $\beta_x$ 
the corresponding $\beta$-set. For $m=0,1,\ldots,2t-3$, let $\rho_{m}$ be the unipotent
character of $G_r$ whose partition corresponds to the $\beta$-set $\beta_{2m-3t+4}$.
The unipotent character $E_t$ correspond to $\beta = \beta_{-t+1}$. Let $\chi_{\text{exc}}$
be the sum of the non-unipotent characters in the block of $E_t$, called
\emph{exceptional characters}.


\smallskip

The Brauer tree of the block is described in \cite[sec.~6]{FS}. Recall that the
vertices in the Brauer tree are labeled by $\chi_\text{exc}$ and the unipotent
characters of the block, while edges are labeled by simple $\k G_r$-modules.
The simple $\k G_r$-module $S$ labels the edge $\chi{-\hskip-2mm - \hskip-2mm -} \chi'$
if and only if $S$ is a composition factor of the $\ell$-reduction of both $\chi$
and $\chi'$. The set of the composition factors of the $\ell$-reduction of $\chi$
is exactly the set of the labels of the edges adjacent to $\chi$. The parity of
$x$ determines on which branch of the Brauer tree the character lies. The Brauer
tree is

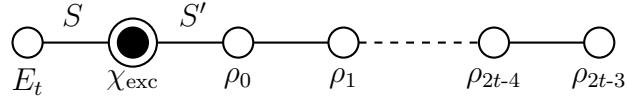
\begin{center}
\begin{tikzpicture}[scale=.4]
\draw[thick] (0 cm,0) circle (.5cm);
\node [below] at (0 cm,-.5cm) {$E_t$};
\node [above] at (1.5 cm,.2cm) {$S$};
\draw[thick] (0.5 cm,0) -- +(2.2 cm,0);
\draw[thick,fill=black] (3.5 cm,0) circle (.5cm);
\draw[thick] (3.5 cm,0) circle (.8cm);
\node [below] at (3.5 cm,-.55cm) {$\chi_{\text{exc}}$};
\node [above] at (5.5 cm,.2cm) {$S'$};
\draw[thick] (4.3 cm,0) -- +(2.2 cm,0);
\draw[thick] (7 cm,0) circle (.5cm);
\node [below] at (7 cm,-.5cm) {$\rho_0$};
\draw[thick] (7.5 cm,0) -- +(2.5 cm,0);
\draw[thick] (10.5 cm,0) circle (.5cm);
\node [below] at (10.5 cm,-.5cm) {$\rho_1$};
\draw[dashed,thick] (11 cm,0) -- +(4 cm,0);
\draw[thick] (15.5 cm,0) circle (.5cm);
\node [below] at (15.5 cm,-.5cm) {$\rho_{\text{2$t$-4}}$};
\draw[thick] (16 cm,0) -- +(2.5 cm,0);
\draw[thick] (19 cm,0) circle (.5cm);
\node [below] at (19 cm,-.5cm) {$\rho_{\text{2$t$-3}}$};
\end{tikzpicture}
\end{center}

Since $E_t$ is weakly cuspidal, \emph{i.e.} $E(E_t) = 0$, so is $S$. 
The character $E(\rho_m)$ can be explicitly computed using Lemma
\ref{lem:indresunip}, and the combinatorics on the $\beta_{2m-3t+4}$
(see \S \ref{sec:combinatorics} for the effect of removing $2$-hooks on a $\beta$-set). 
Two cases arise: if $m=0$ (resp. $m=2t-3$), then $E(\rho_m)$ is irreducible, and the
corresponding $\beta$-set is obtained by changing $-3t+6$ to $-3t+4$ (resp. $3t-3$ to $3t-5$);
otherwise $E(\rho_m)$ has two constituents whose $\beta$-sets are obtained
by changing $2m-3t+6$ to $2m-3t+4$ or $2m-t+3$ to $2m-t+1$. 
We deduce that $\sum (-1)^m E(\rho_m) = 0$ in $[K G_r\mod]$.
Since $[S'] = \sum (-1)^m d_\scrU([\rho_m])$ in $[\k G_r\mod]$, 
this implies that the $\k G_r$-module $S'$ is also weakly cuspidal. Therefore the two composition factors of the 
$\ell$-reduction of the exceptional characters are weakly
cuspidal, which forces the exceptional characters to be weakly cuspidal as well.


\smallskip

Let $\chi$ be one of the exceptionals characters. Since $\chi$ is weakly cuspidal,
we can use \cite{HL80} to show that the operator $X(\chi)$ on $F(\chi)$ satisfies a
quadratic relation and the product of its eigenvalues equals $-q^{-1}$. In particular,
if one of the eigenvalues of $X(\chi)$ is a power of $-q$ modulo $\ell$ then so is the other. Now $\rho_0$ is not cuspidal, because $E_t$ is the unique cuspidal unipotent character of $KG_r$. Therefore $\rho_0$ lies in the Harish-Chandra series
$\Irr(KG_r,E_s)$ for some $s <t$. The induction hypothesis implies that the
parameters of $\scrH(KG_r,E_s)$ are  $Q_s$ and $(-q)^2$. Hence the eigenvalues of
$X(\rho_0)$ are powers of $-q$ modulo $\ell$ by \eqref{Q}, see \S \ref{subsec:minimalcat}.
Since the $\ell$-reductions of $F(\rho_{0})$ and $F(\chi)$ share a common
composition factor, we deduce that the eigenvalues of $X(\chi)$ are powers
of $-q$ modulo $\ell$. On the other hand, the $\ell$-reduction of $E_t$
is isomorphic to $S$, which is a composition factor of the $\ell$-reduction of $\chi$. 
Therefore, the eigenvalues of $X(E_t)$ are also powers of $-q$ modulo $\ell$.

\smallskip

Finally, since $e=2t-1$ is odd and $-q$ is of order $e$ modulo $\ell$, $-1$ is not a power of $-q\equiv -q^{-t-1}$ modulo $\ell$. We deduce that the eigenvalues of $X(E_t)$ are powers of $-q$. 
\end{proof}

\subsubsection{Parametrization of the weak Harish-Chandra series of $\scrU_K$}
\label{subsec:hcseries-char0}
Let $W(B_m)$ be the Weyl group of type $B_m$, and $t_0,t_1,\ldots,t_{m-1}$ be the
generators corresponding to the following Dynkin diagram
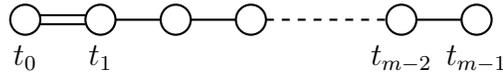
\begin{align*}
\begin{split}
\begin{tikzpicture}[scale=.4]
\draw[thick] (0 cm,0) circle (.5cm);
\node [below] at (0 cm,-.5cm) {$t_0$};
\draw[thick] (0.5 cm,-0.15cm) -- +(1.5 cm,0);
\draw[thick] (0.5 cm,0.15cm) -- +(1.5 cm,0);
\draw[thick] (2.5 cm,0) circle (.5cm);
\node [below] at (2.5 cm,-.5cm) {$t_1$};
\draw[thick] (3 cm,0) -- +(1.5 cm,0);
\draw[thick] (5 cm,0) circle (.5cm);
\draw[thick] (5.5 cm,0) -- +(1.5 cm,0);
\draw[thick] (7.5 cm,0) circle (.5cm);
\draw[dashed,thick] (8 cm,0) -- +(4 cm,0);
\draw[thick] (12.5 cm,0) circle (.5cm);
\node [below] at (12.5 cm,-.5cm) {$t_{m-2}$};
\draw[thick] (13 cm,0) -- +(1.5 cm,0);
\draw[thick] (15 cm,0) circle (.5cm);
\node [below] at (15 cm,-.5cm) {$t_{m-1}$};
\end{tikzpicture}
\end{split}
\end{align*}
Let us first recall the construction of the irreducible characters of $W(B_m)$
(see for example \cite[def. 5.4.4]{GP} or \cite[sec. I.7]{M}). We denote by 
$\sigma_m$ the linear character of $W(B_m)$ such that $\sigma_m(t_0) = -1$
and $\sigma_m(t_i) = 1$ for all $i>0$. Given $\lambda \vdash m$ a partition 
of $m$, we write $\widetilde{\phi}_\lambda$ for the inflation to $W(B_m)$
of the irreducible character of $\frakS_m$ corresponding to $\lambda$. 
Given $a,b$ such that $a + b = m$, one can consider the subgroups
$\frakS_a\times \frakS_b \subset W(B_a) \times W(B_b)$ of $W(B_m)$ where
$W(B_a)\times \frakS_b$ is a parabolic subgroup generated by $\{t_0,t_1,\ldots,t_{a-1}\}\cup \{t_{a+1},\ldots,t_{m-1}\}$, and $W(B_a) \times W(B_b)$ is obtained by adding the
reflection $t_{a}\cdots t_0 \cdots t_a$. The irreducible character of $W(B_m)$
associated to a bipartition $(\lambda,\mu)$ of $m$ is
 $$ \scrX_{\lambda,\mu} = \mathrm{Ind}_{W(B_{|\lambda|}) \times W(B_{|\mu|})}^{W(B_m)}
  \, (\widetilde{\phi}_\lambda \boxtimes \sigma_{|\mu|} \widetilde{\phi}_\mu).$$
For example, $\scrX_{(m),\emptyset}$ is the trivial character, whereas
$\scrX_{\emptyset,(1^m)}$ is the signature. We also have $\sigma_m = 
\scrX_{\emptyset,(m)}$ and more generally $\scrX_{\lambda,\mu} = \sigma_m
\scrX_{\mu,\lambda}$. 

\smallskip

By Tits deformation theorem, the evaluation $q \mapsto 1$ yields a bijection 
$$ \Irr( \bfH^{Q_t,\,q^2}_{K,m}) \, \mathop{\longleftrightarrow}\limits^{1:1} \, \Irr(W(B_m)) $$
from which we obtain a canonical labelling of the irreducible representations of
$\bfH^{Q_t,\,q^2}_{K,m}$. We write $\Irr( \bfH^{Q_t,\,q^2}_{K,m}) = 
\{S(\lambda,\mu)^{Q_t,q^2}_K\}_{(\lambda,\mu) \in \scrP^2_m}$, compare \S \ref{subsec:Fock}.
In \cite{HL80}, Howlett-Lehrer use a renormalization of \eqref{eq:xsign}. Setting
$T_0 = - \epsilon_t (-q)^{-1-t} X $ $= (-1)^t q^{-1-t} X$, we have now the quadratic relation
 $$(T_0 + 1)(T_0-q^{2t+1}) = 0.$$
Using this generator instead of $X$, we obtain the usual presentation for a Hecke algebra
of type $B_m$ with parameters $(q^{2t+1},q^2)$. The endomorphism of $KW(B_m)$ which is obtained
from the renormalization $\bfH^{Q_t,\,q^2}_{K,m} \simto  \bfH^{(q^{2t+1},1),\,q^2}_{K,m}$
at $q=1$ is the identity on $\frakS_m$ but sends $t_0$ to $(-1)^tt_0$. Therefore 
this renormalization sends $S(\lambda,\mu)$ to $S(\lambda,\mu)$ if $t$ is even, and
to $S(\mu,\lambda)$ is $t$ is odd. Combining this observation with \cite[Appendix]{FS},
we obtain the following parametrization of the unipotent characters in the series of
$E_t$. 

\begin{corollary}\label{cor:bijection} 
Let $t,m \geqslant 0$ and $n = t(t+1)/2+2m$. Then the map $\phi_{K,m}$ and the functor $\frakE_{F^m(E_t)}$
induce a bijection
  $$ \wIrr(K G_n,E_t) \, \mathop{\longleftrightarrow}\limits^{1:1} \, \Irr( \bfH^{Q_t,\,q^2}_{K,m})$$
sending $E_\lambda$ to $S(\lambda^{[2]})^{Q_t,q^2}_K$ for all partitions $\lambda \vdash n$ with
$2$-core $\Delta_t = (t,t-1,\ldots,1)$. 
\qed
\end{corollary}

\subsubsection{The $\frakg_\infty$-representation on $\scrU_K$}
\label{subsec:g-infty}
The functors $E,F$ preserve the subcategory $\scrU_K$ by Proposition \ref{prop:stable},
hence $(E,F,X,T)$ yields a representation datum on $\scrU_K$. In order to extend it to a
categorical representation on $\scrU_K$, one should consider the quiver
$\scrI(q^2)$ with vertices given by the various eigenvalues of $X$ (which
we showed to be all powers of $-q$ in the proof of Theorem \ref{thm:cat}) and arrows
$i \longrightarrow q^2i$.

\smallskip

In this section we will view the integer $q$ as an element
of $K^\times$ in the obvious way. For the construction of Kac-Moody algebras associated to quivers 
we refer to \S \ref{sec:quivers}. 

\begin{definition}\label{def:g-infty} 
Let $\scrI_\infty$ denote the subset $(-q)^\bbZ$ of $K^\times.$
We define $\frakg_{\infty}$  to be the (derived)
Kac-Moody algebra associated to the quiver  $\scrI_\infty(q^2)$.
\end{definition}

To avoid cumbersome notation, we will write for short $\scrI_{\infty} = \scrI_\infty(q^2)$,
and $(-)_\infty=(-)_{\scrI_\infty}$.
We denote by $\{\Lambda_i\}$, $\{\alpha_i\}$ and $\{\alpha_i^\vee\}$ the fundamental
weights, simple roots and simple coroots of $\frakg_\infty$. Here $\X_\infty$ coincides with
$\P_{\infty} = \bigoplus \bbZ \Lambda_i$. Consequently, there is a Lie algebra isomorphism $(\fraks\frakl_\bbZ)^{\oplus 2} \simto \frakg_{\infty}$ such that $(\alpha^\vee_d,0)\longmapsto\alpha^\vee_{-q^{2d-1}}$ and $(0,\alpha^\vee_d)\longmapsto\alpha^\vee_{q^{2d}}$.

\smallskip

For any $t,m,n\in\bbN$, let $(KG_n,E_t)\mod$ be the Serre subcategory of $\scrU_K$ generated by
the modules $F^m(E_t)$ with $n=r+2m$ and $r=t(t+1)/2$. We define
  $$\scrU_{K,t}=\bigoplus_{n\geqslant 0}\, (KG_n,E_t)\mod.$$
Then $\Irr((KG_n,E_t)\mod)=\wIrr(KG_n,E_t)$, which implies $\scrU_K=\bigoplus_{t\geqslant 0}\scrU_{K,t}$
by Proposition \ref{prop:partitionwcusp}. We can now prove our first categorification result, which 
says that $\scrU_{K,t}$ is a $\frakg_\infty$-representation which categorifies the Fock space $\bfF(Q_t)_\infty$.

\begin{theorem}\label{thm:char0}
 Let $t \geqslant 0$ and $Q_t$ be as in \eqref{Q}. 
 \begin{itemize}[leftmargin=8mm]
  \item[$\mathrm{(a)}$] 
  The Harish-Chandra induction and restriction functors yield a representation of
  $\frakg_{\infty}$ on $\scrU_{K,t}$ which is isomorphic to $\scrL(\Lambda_{Q_t})_\infty$. 
  \item[$\mathrm{(b)}$] The map $|\mu,Q_{t}\rangle_\infty\longmapsto [E_{\varpi_t(\mu)}]$ induces an isomorphism
  of $\frakg_\infty$-modules $$\bfF({Q_t})_{\infty} \simto [\scrU_{K,t}].$$
 \end{itemize}
\end{theorem}

\begin{proof}
Composing the functor $\frakE_{F^m(E_t)}$ with the algebra isomorphism in Theorem \ref{thm:cat}
and taking the sum over all $m\in\bbN$, we get an equivalence of semi-simple abelian $K$-categories
$\frakE_{t}\, :\, \scrL(\Lambda_{Q_t})_\infty \simto \scrU_{K,t}$. We claim that it is actually
an isomorphism of representation data. To see this, we first observe that the representation datum
$(E,F,X,T)$ restricts to a representation datum on $\scrU_{K,t}$: indeed, by definition $\scrU_{K,t}$
is stable by $F$ and by the Mackey formula and \cite[prop.~2.2]{GHJ} it is also stable by the adjoint
functor $E$. Therefore to prove our claim we must show that 
\begin{itemize}[leftmargin=8mm]
  \item[(i)] there are isomorphisms $\frakE_t E\simeq E\frakE_t$ and $\frakE_t F\simeq F\frakE_t$
   of functors $\scrL(\Lambda_{Q_t})_{\infty}\to\scrU_{K,t}$,
  \item[(ii)] the isomorphisms $\frakE_t F\simeq F\frakE_t$  and $\frakE_t F^2\simeq F^2\frakE_t$
   intertwine the endomorphisms $\frakE_t X$, $X\frakE_t$ and $\frakE_t T$, $T\frakE_t$.
\end{itemize}
Note that to prove (i), it is enough to show $\frakE_t F\simeq F\frakE_t$, since 
$\frakE_t E\simeq E\frakE_t$ will follow by adjunction.

\smallskip

For the rest of the proof we will write for short $\bfH_{m} := \bfH^{Q_t,q^2}_{K,m}$. 
Set $n=r+2m$ and $r=t(t+1)/2$. Then the functors
  $$F\frakE_{F^m(E_t)},\ \frakE_{F^{m+1}(E_t)}F\,:\,\bfH_{m}\mod\longrightarrow(KG_{n+2},E_t)\mod$$ 
are both obtained by tensoring with the $(KG_{n+2},\bfH_{m})$-bimodule
$F^{m+1}(E_t)$. More precisely, the left action of $KG_{n+2}$ is the same in both cases, while the
right action of $\bfH_{m}$ comes from the right action of $\bfH_m$ on $F^m(E_t)$ and the functoriality
of $F$ in the first case, and from the right action of $\bfH_{m+1}$ on $F^{m+1}(E_t)$ and the obvious
inclusion  $\bfH_m \subset\bfH_{m+1}$ in the second case. Assertion (a) follows.

\smallskip

Claim (ii) is obvious. Indeed, given an integer $d>0$, let $x\in\bfH^{q^2}_{K,d}$ be an element
in the affine Hecke algebra $\bfH^{q^2}_{K,d}$ and let $M\in \bfH_m\mod$.
Then, we have $F^d(M)\in \bfH_{d+m}\mod$.
The action of $x$ on $F^d\frakE_{F^m(E_t)}(X)$ and on $\frakE_{F^{d+m}(E_t)}F^d(M)$ are represented
respectively by the action of $\phi_{F^d}(x)\otimes 1$ on 
\begin{align*}
F^d(F^m(E_t))\otimes_{\bfH_m}M
&=F^{d+m}(E_t)\otimes_{\bfH_{d+m}}\bfH_{d+m}\otimes_{\bfH_{m}}M\\
&=F^{d+m}(E_t)\otimes_{\bfH_{d+m}}F^d(M)
\end{align*}
and by the action of $1\otimes x$ on $F^d(M)$ in
$F^{d+m}(E_t)\otimes_{\bfH_{d+m}}F^d(M)$. They obviously coincide. The claim follows, taking
$d=1$, $x=X$ or $d=2$, $x=T$.

\smallskip

Now we can finish the proof of the theorem. We equip $\scrU_{K,t}$ with the
$\frakg_\infty$-representation which is transferred from the $\frakg_\infty$-representation on
$\scrL(\Lambda_{Q_t})_{\infty}$ via the equivalence $\frakE_t$. This proves (a).
We deduce that $\frakE_t$ induces on the Grothendieck groups a $\frakg_{\infty}$-module isomorphism
$\bfL(\Lambda_{Q_t})_\infty = [\scrL(\Lambda_{Q_t})_{\infty}] \simto [\scrU_K]$. Then
(b) follows from Corollary \ref{cor:bijection} and the $\frakg_{\infty}$-module isomorphism
$\bfF({Q_t})_{\infty}=\bfL(\Lambda_{Q_t})_{\infty}$ given in Propositions \ref{prop:subrepinfock}, \ref{prop:explicitiso}.
\end{proof}

\begin{remark} (a) The functor $F$ on $\scrU_K$ is represented by the sum of bimodules
$\bigoplus_{n\in\bbN} K G_{n+2}\,e_{{n+2,n}}$ on which the endomorphism $X\in\End(F)$ acts by
right multiplication by the element
$$\label{X}\sum_n(-1)^nq^ne_{{n+2,n}}\,(1,n+2)\,e_{{n+2,n}}.$$
(b) The functor $F_i : (KG_n,E_t)\mod\longrightarrow (KG_{n+2},E_t)\mod$ is the
generalized eigenspace of $X\in\End(F)$ associated with the eigenvalue $i$. For each
bipartition $\mu$ we have
$$F_i (E_{\varpi_t(\mu)}) = \left\{\begin{array}{l} E_{\varpi_t(\nu)} \text{ if } 
\res(\nu-\mu,Q_t)_\infty = i, \\
0 \text{ otherwise.} \end{array}\right.$$
\end{remark}

\subsection{The $\frakg_e$-representation on $\scrU_\k$.}
\label{sec:g-e}
By Proposition \ref{prop:stable}, 
the representation datum $(E,F,X,T)$ on $\k G\mod$ induces a representation datum
on $\scrU_\k$. Since the abelian category $\scrU_\k$ is not semisimple,
to extend the representation datum to a categorical $\frakg$-representation
one needs to prove
that weight spaces of $\scrU_\k$ are sums of blocks. This will be done combinatorially by studying a
representation of a bigger Lie algebra $\frakg_\circ$, which is a $(-q)$-analogue of the action of Harish-Chandra
induction and restriction on unipotent representations of $\mathrm{GL}_n(q)$. By definition, this
action is compatible with the decomposition into $\ell$-blocks of $\mathrm{GL}_n(q)$ and one can transfer this 
property to unitary groups using the correspondence between $\mathrm{GU}_n(q)$ and $\mathrm{GL}_n(-q)$.

\subsubsection{The Lie algebras $\frakg_e$ and $\frakg_{e,\circ}$}
\label{subsec:g-e}
Recall from Theorem \ref{thm:cat} that the eigenvalues of $X$ on $E$ and
$F$ are all powers of $-q$. If we denote again by ${q}$ the image of $q$ under the canonical map 
$\scrO \twoheadrightarrow \k$, then
the eigenvalues of $X$ on $\k E$ and $\k F$ 
belong to the finite set $(-{q})^\mathbb{Z} \subset \k^\times$. This set has exactly $e$ elements, where $e$ is the order of $-{q}$ in $\k^\times$. 

\begin{definition}
We define $\scrI_e$ to be the subset $(-{q})^\bbZ$ of $\k^\times$. We denote by 
$\scrI_e$ and $\scrI_{e,\circ}$ the finite quivers $\scrI_e({q}^2)$ and
$\scrI_e(-{q})$. 
\end{definition}

The quivers $\scrI_e$ and $\scrI_{e,\circ}$ have the same set of vertices, but the arrows in 
$\scrI_e$ are the composition of two consecutive arrows in $\scrI_{e,\circ}$. The quiver 
$\scrI_{e,\circ}$ is cyclic, whereas the quiver $\scrI_{e}$ is cyclic if $e$ is odd,
and is a union of two cyclic quivers if $e$ is even. Therefore the corresponding
Kac-Moody algebras are isomorphic to $\widehat{\fraks\frakl}_e$ or
$(\widehat{\fraks\frakl}_{e/2})^{\oplus2}$.

\smallskip

To avoid cumbersome notation, we will write  $(\bullet)_{e,\circ} = (\bullet)_{\scrI_{e,\circ}}$
and $(\bullet)_{e} = (\bullet)_{\scrI_{e}}$. 
We must introduce $\frakg_{e,\circ}$ and $\frakg_{e}$ such that 
$\frakg_{e,\circ}' = [\frakg_{e,\circ},\frakg_{e,\circ}]$ and $\frakg_{e}' = 
[\frakg_{e},\frakg_{e}]$. 
The Chevalley generators of $\frakg_{e,\circ}'$ and $\frakg_{e}' $ are $e_{i,\circ}, f_{i,\circ}$ and 
$e_i,f_i$ respectively, for $i \in (-q)^\bbZ$. It is easy to see that there exists a morphism of Lie algebras $\kappa : \frakg_e' \longrightarrow
\frakg_{e,\circ}'$ defined by 
 $$\kappa(e_i) = [e_{-qi,\circ},e_{i,\circ}] \quad \text{and} \quad
    \kappa(f_i)= [f_{-qi,\circ},f_{i,\circ}]. $$
It restricts to a map between the coroot lattices sending $\alpha_i^\vee$ to $\alpha_{i,\circ}^\vee
+\alpha_{-qi,\circ}^\vee$

\smallskip

We denote by $\frakg_{e,\circ}$ the Kac-Moody algebra associated with the lattices
$\X_{e,\circ} = \P_{e,\circ} \oplus \bbZ\delta_\circ$ and $\X_{e,\circ}^\vee = 
\Q_{e,\circ}^\vee \oplus \bbZ \partial_\circ$, where $\delta_\circ = \sum \alpha_{i,\circ}$, $\partial_\circ = \Lambda_{1,\circ}^\vee$
and the pairing $\X_{e,\circ}^\vee
\times \X_{e,\circ} \longrightarrow \bbZ$ is given by
 $$ \langle \alpha_{j,\circ}^\vee, \Lambda_{i,\circ} \rangle_{e,\circ} = \delta_{ij}, \quad
 \langle \partial_\circ,\Lambda_{i,\circ} \rangle_{e,\circ} = \langle \alpha_{j,\circ}^\vee,
 \delta_\circ \rangle_{e,\circ} = 0, \quad
  \langle \partial_\circ,\delta_\circ \rangle_{e,\circ} = 1.$$
Then $\frakg_{e,\circ}$, $\frakg_{e,\circ}'$ are isomorphic to
$\widehat{\fraks\frakl}_e$, $\widetilde{\fraks\frakl}_e$ (see Example \ref{ex:onegenerator}). 

\smallskip

Let  $\widetilde\frakg_e$ be the usual Kac-Moody algebra associated with $\scrI_e$.
Its derived Lie subalgebra is equal to $\frakg'_e$.
Let $\widetilde \X_e$ and $\widetilde \X_e^\vee$ be the lattices corresponding to $\widetilde\frakg_e$.
If $e$ is odd, then $\scrI_e$ is a cyclic
quiver and $\widetilde\frakg_e$ is isomorphic to $\widehat{\fraks\frakl}_e$. 
Let $\alpha_1$ be the affine root, then we have $\widetilde \X_e = \P_e \oplus \bbZ \tilde\delta$ and 
$\widetilde \X_e^\vee = \Q_e^\vee \oplus \bbZ \widetilde \partial$
with $\tilde\delta = \sum \alpha_i$ and $\widetilde\partial = \Lambda_1^\vee$. 
If $e$ is even, then $\scrI_e$ is the disjoint union
of two cyclic quivers and $\tilde\frakg_e$ is isomorphic
to $(\widehat{\fraks\frakl}_{e/2})^{\oplus2}$. 
Let $\alpha_1$ and $\alpha_{-q^{-1}}$
be the affine roots, then we have $\widetilde \X_e = \P_e \oplus
 \bbZ \delta_1 \oplus \bbZ\delta_2$ and $\widetilde \X_e^\vee = \Q_e^\vee \oplus \bbZ \partial_1
\oplus \bbZ \partial_2$ with $\delta_1 = \sum_{j \text{ odd}} \alpha_{-q^{j}}$, 
$\delta_2 = \sum_{j \text{ even}} \alpha_{q^{j}}$, $\partial_1 = \Lambda_{-q^{-1}}^\vee$ and
$\partial_2 = \Lambda_{1}^\vee$. We abbreviate $\widetilde\partial=\partial_1+\partial_2$ and $\widetilde\delta=\delta_1+\delta_2$.

\smallskip

The map
$\kappa : \frakg_e' \longrightarrow
\frakg_{e,\circ}'$ may not extend to a morphism of Lie algebra $\widetilde\frakg_e \longrightarrow \frakg_{e,\circ}$.
For this reason we'll define $\frakg_e$ to be a  Lie subalgebra of $\widetilde\frakg_e$ containing $\frakg'_e$. 
More precisely, we define $\frakg_e$ to be Lie subalgebra $\frakg'_e\oplus\bbC\partial$ of $\widetilde\frakg_e,$
where $\partial$ is the element given by $\partial = \Lambda_{1}^\vee+\Lambda_{-q^{-1}}^\vee$. 
If $e$ is odd, then $\partial=2\widetilde\partial+h$ for some coweight $h$ of $\fraks\frakl_e$, hence we have $\frakg_e=\widetilde\frakg_e$.
If $e$ is even, then $\partial=\widetilde\partial$ and $\frakg_e$ is strictly smaller than $\widetilde\frakg_e$.

\smallskip

In both cases we can view $\frakg_e$ as the Kac-Moody algebra associated with the lattice 
$\X_{e}^\vee =  \Q_{e}^\vee \oplus \bbZ \partial$ above and a lattice $\X_{e} = \P_{e} \oplus \bbZ\delta$ that we now define.
If $e$ is odd we set $\delta=\tilde\delta/2$, if $e$ is even we set $\X_e = \widetilde \X_e / (\delta_1-\delta_2)$ with $\delta=\widetilde\delta/2$.
The perfect pairing $\X_e^\vee\times\X_e\longrightarrow\bbZ$ 
is induced in the obvious way by the pairing $\widetilde\X_e^\vee\times\widetilde\X_e\longrightarrow\bbZ$. 
We have 
$$\langle \alpha_{j}^\vee,\Lambda_{i} \rangle_{e} = \delta_{ij},
\quad \langle \partial,\alpha_j  \rangle_e = \delta_{j1}+\delta_{j,-q^{-1}},
 \quad \langle \partial,\Lambda_1\rangle_{e} = \langle\alpha_j^\vee,\delta\rangle_e = 0,\quad
\langle\partial,\delta\rangle_e=1.$$
This definition of $\frakg_e$ ensures that $\kappa$ extends to a Lie algebra
homomorphism $\frakg_e \longrightarrow \frakg_{e,\circ}$.

\begin{lemma}\label{lem:kappa}
There is a well-defined morphism of Lie algebras
$\frakg_e \longrightarrow \frakg_{e,\circ}$ which extends $\kappa$ whose restriction to $\X_{e}^\vee$ 
is given by  
$$ \kappa(\alpha_i^\vee) = \alpha_{i,\circ}^\vee + \alpha_{-qi,\circ}^\vee, 
  \quad \kappa(\partial)= \partial_\circ.$$
The restriction $\kappa : \X^\vee_e\longrightarrow\X^\vee_{e,\circ}$ has an adjoint $\kappa^* : \X_{e,\circ} \longrightarrow \X_{e}$ such that
$$ \kappa^*(\Lambda_{i,\circ})\equiv\Lambda_{i} + \Lambda_{-q^{-1}i}\ \rm{mod}\ \delta, 
  \quad  \kappa^*(\delta_\circ)= \delta.$$
\end{lemma} 

\begin{proof} 
Recall that $i \in (-q) ^\bbZ$ and that we have already defined the Lie algebra homomorphism 
$\kappa : \frakg'_e \longrightarrow \frakg'_{e,\circ}$ sending $e_i$ to $\kappa(e_i) = [e_{-qi,\circ},e_{i,\circ}]$. 
We set $\kappa(\partial)= \partial_\circ$ and $\kappa(e_i)= [e_{-qi,\circ},e_{i,\circ}]$ in $\frakg_{e,\circ}$.
By definition of $\frakg_e$ we have 
$$[\partial,e_i] = \langle \partial, \alpha_i \rangle_{e}\, e_i = (\delta_{i1} + \delta_{i,-q^{-1}}) e_i.$$ 
Thus, we have
$$[\kappa(\partial),\kappa(e_i)] = \langle \partial_\circ, \alpha_{i,\circ} + 
\alpha_{-qi,\circ}\rangle_{e,\circ} \, \kappa(e_i) = (\delta_{i1} + \delta_{i,-q^{-1}}) \kappa(e_i)
=\kappa([\partial,e_i]),$$
because the weight of
$\kappa(e_i)$ in $\frakg_{e,\circ}$ is $\alpha_{i,\circ} +
\alpha_{-qi,\circ}$. 
The same holds for $f_i$ instead of $e_i$. 
The second claim is an easy computation using the relation $\langle \kappa(h),\alpha_\circ \rangle_{e,\circ}
 = \langle h,\kappa^*(\alpha_\circ)\rangle_e$ for all $h\in \X_e^\vee$ and $\alpha_\circ\in \X_{e,\circ}$. 
\end{proof}

\begin{remark} If $e$ is even then $\langle\partial,\Lambda_i\rangle_e=0$, hence $ \kappa^*(\Lambda_{i,\circ})=\Lambda_{i} + \Lambda_{-q^{-1}i}$ for all $i$.
\end{remark}

\subsubsection{Action of $\frakg'_e$ on $[\scrU_\k]$}
\label{subsec:derivedaction-uk}
The quotient map $ \scrO \twoheadrightarrow \k$ yields a morphism of quivers
$\sp : \scrI_\infty \longrightarrow \scrI_e$ and a surjective morphism of abelian groups
$\sp : \P_\infty \twoheadrightarrow \P_e$  sending $\Lambda_i$ to $\Lambda_{\sp(i)}$.
In addition any integrable representation $V$ of $\frakg_\infty$ can be ``restricted''
to an integrable representation of the derived algebra $\frakg'_e$, where
$e_i \in \frakg'_e$ (resp. $f_i \in \frakg'_e$) act as $\sum_{\sp(j)=i} e_j$ (resp.
$\sum_{\sp(j)=i} f_j$). From the definition of the action of $\frakg_\infty$ and 
$\frakg_e$ on Fock spaces, see \eqref{eq:EF}, we deduce that the map
$|\mu,Q_{t}\rangle_\infty\mapsto|\mu,Q_t\rangle_e$ induces the following isomorphism
of $\frakg'_e$-modules
 \begin{equation}\label{eq:fockinftyande} 
 \sp \, : \, \mathrm{Res}_{\frakg'_e}^{\frakg_\infty} \, \bfF(Q_t)_\infty
 \, \mathop{\longrightarrow}\limits^\sim\, \bfF(Q_t)_e.\end{equation}
We show that under the decomposition map, this isomorphism endows $[\scrU_\k]$ with
a structure of $\frakg'_e$-module which is compatible with the one coming
from the representation datum.

\begin{proposition}\label{prop:lreduction}
For each $i \in \scrI_e$, let $\k E_i$ (resp. $\k F_i$)
be the generalized $i$-eigenspace of $X$ on $\k E$ (resp. $\k F$). 
  \begin{itemize}[leftmargin=8mm]
    \item[$\mathrm{(a)}$] The operators $[\k E_i], [\k F_i]$ endow $[\scrU_\k]$ with a structure
    of $\frakg'_e$-module.
    \item[$\mathrm{(b)}$] The decomposition map $d_\scrU : [\scrU_K] \longrightarrow [\scrU_k]$
     is a $\frakg'_e$-module isomorphism 
     $$ \mathrm{Res}_{\frakg'_e}^{\frakg_\infty} \, [\scrU_K]\simto [\scrU_k].$$
  \end{itemize}
\end{proposition}

\begin{proof}
Recall from \S\ref{sec:rep-datum} that the endofunctor $F$ of $\scrO G\mod$ is represented by the
bimodule $\bigoplus_{n\in\bbN} \scrO G_{n+2}\,e_{{n+2,n}}.$
Under base change from $\scrO$ to $K$ and $\k$, it yields a functor
$KF$ on $KG\mod$ and a functor $\k F$ on $\k G\mod$. They are represented by the bimodules
$\bigoplus_{n\in\bbN} K G_{n+2}\,e_{{n+2,n}}$ and $\bigoplus_{n\in\bbN} \k G_{n+2}
\,e_{{n+2,n}}$ respectively. 
For any $\scrO G$-module $M$ we write $KM=K\otimes_\scrO M$ and $\k M=\k\otimes_\scrO M$.
The associativity of the tensor product implies that
$KF(KM)=K(FM)$ and $\k F(\k M)=\k(FM).$
\smallskip

Similarly, the endomorphism $X$ of $F$ which acts by right multiplication by the
element \eqref{X} yields an endomorphism of $KF$ and of $\k F$.
For each $i\in\scrI_\infty$ (resp. $i\in\scrI_e$), let $K F_i$ 
(resp. $\k F_i$) be the generalized $i$-eigenspace of $X$ acting on $K F$ (resp.
on $\k F$). We define $K E_i$ and $\k E_i$ in a similar way.
Finally, for $i \in \scrI_e$, we write $E_i=(\bigoplus_{\sp(j)=i}KE_j)\cap E$ 
and $F_i=(\bigoplus_{\sp(j)=i}KF_j)\cap F$. Notice that, although the sums above are
a priori infinite sums, they are well-defined as subfunctors of $E$ and $F$.
Then, we have
\begin{itemize}[leftmargin=8mm]
  \item[(i)] $E=\bigoplus_{i\in\scrI_e}E_i$ and $F=\bigoplus_{i\in\scrI_e}F_i$,

  \item[(ii)] for each $i\in\scrI_e$, the functors $\k E_i$, $\k F_i$ are isomorphic to
  the specialization of $E_i$, $F_i$ to $\k$.
\end{itemize}

\smallskip
We deduce that the decomposition map $d_\scrU : [\scrU_K] \longrightarrow 
[\scrU_k]$ is a $\bbC$-linear isomorphism which intertwines the operators
$\bigoplus_{\sp(j)=i} KE_j$, $\bigoplus_{\sp(j)=i} KF_j$ on $[\scrU_K]$ with the
operators $\k E_i$, $\k F_i$ on $[\scrU_k]$ for each $i\in \scrI_e$.
\end{proof}

We can fit the three isomorphisms written in Proposition \ref{prop:lreduction}(b), 
\eqref{eq:fockinftyande} and Theorem \ref{thm:char0}(b) into the following diagram:
$$\xymatrix{ \mathrm{Res}_{\frakg'_e}^{\frakg_\infty}\, [\scrU_K]\ar[d]^{\rotatebox{90}{$\sim$}}
\ar[r]_{\ \ d_\scrU}^{\ \ \sim} &[\scrU_\k] \ar@{-->}[d]\\
 \mathrm{Res}_{\frakg'_e}^{\frakg_\infty} \big(\displaystyle\bigoplus_{t\in\bbN}\bfF({Q_t})_{\infty}\big)\ar[r]_{\qquad \sp}^{\qquad \sim} &\displaystyle\bigoplus_{t\in\bbN}\bfF({Q_t})_{e}.
}$$
By composition, we get an explicit description of the $\frakg'_e$-module
structure on~$[\scrU_\k]$. 

\begin{corollary}\label{cor:derivedaction-uk}
The map $|\mu, Q_t\rangle_e \longmapsto [V_{\varpi_t(\mu)}]$ induces a
$\frakg'_e$-module isomorphism
$$\bigoplus_{t\in\bbN}\bfF({Q_t})_{e} \, \mathop{\longrightarrow}\limits^\sim \, [\scrU_\k].$$
\qed
\end{corollary}

\subsubsection{Action of $\frakg_e$ on $[\scrU_\k]$}\label{subsec:action-uk}
We now define an action of $\frakg_e$ on $[\scrU_\k]$ by extending the action from $\frakg'_e$
to $\frakg_e$ on $\bigoplus_{t\in\bbN}\bfF({Q_t})_{e}$. This amounts to defining the action of
$\partial$, or equivalently to extending the grading from $\P_e$ to $\X_e = \P_e \oplus 
\bbZ \delta$. To this end we shall use the (level 1) action of $\frakg_{e,\circ}$ on the
Fock space $\bfF(1)_{e,\circ}$. 

\smallskip

Recall from Lemma \ref{lem:kappa} that $\kappa : \frakg_e \longrightarrow \frakg_{e,\circ}$ is a 
Lie algebra homomorphism such that
$e_i  \longmapsto [e_{-qi,\circ},e_{i,\circ}],$
 $f_i  \longmapsto [f_{-qi,\circ},f_{i,\circ}]$ and $\partial \longmapsto  \partial_{\circ}.$
Any integrable $\frakg_{e,\circ}$-representation 
(resp. $\frakg'_{e,\circ}$-representation) can be ``restricted'' to an 
integrable $\frakg_e$-representation (resp. $\frakg'_{e}$-representation)
through $\kappa$. We denote by $\mathrm{Res}_{\frakg_{e}}^{\frakg_{e,\circ}}$ and 
$\mathrm{Res}_{\frakg'_{e}}^{\frakg'_{e,\circ}}$
the corresponding operations.

\begin{lemma}\label{lem:isofocks}
The map $|\mu,Q_t\rangle_e \longmapsto (-1)^{a(\varpi_t(\mu))}
|\varpi_t(\mu),1\rangle_{e,\circ}$
induces an isomorphism of $\frakg'_{e}$-modules
  $$ \bigoplus_{t \in \bbN} \bfF(Q_t)_e \ \mathop{\longrightarrow}
  \limits^\sim\ \mathrm{Res}_{\frakg'_{e}}^{\frakg'_{e,\circ}}
  \, \bfF(1)_{e,\circ}$$
where $a$ is Lusztig's $a$-function (see \cite[4.4.2]{Lu84}). 
\end{lemma}

\begin{proof} The map is clearly an isomorphism of vector spaces. 
We first show the compatibility of the action for the Lie algebras coming from the
quiver in characteristic zero. Recall that $\scrI_{\infty} = (-q)^\bbZ$ with the action
of $q^2$, where $q$ is seen as an element of $K$. Let $\scrI_{\infty,\circ}$
be the quiver with the same set of vertices, but with arrows given by multiplication by $-q$. Let
$\frakg_{\infty,\circ}$ be the corresponding derived Lie algebra. It is isomorphic
to $\fraks\frakl_\bbZ$ whereas $\frakg_{\infty}$ is isomorphic to $(\fraks\frakl_\bbZ)^{\oplus 2}$.
As before, we can embed $\frakg_{\infty}$ into $\frakg_{\infty,\circ}$ by sending
$e_i$ to $[e_{-qi,\circ},e_{i,\circ}]$ and $f_i$ to $[f_{-qi,\circ},f_{i,\circ}]$.
Now, the isomorphism in the lemma can be deduced from the isomorphism of $\frakg_\infty$-modules
\begin{equation}\label{eq:fockcharO} \bigoplus_{t \in \bbN} \bfF(Q_t)_\infty \ \mathop{\longrightarrow}
  \limits^\sim\ \mathrm{Res}_{\frakg_{\infty}}^{\frakg_{\infty,\circ}}
  \, \bfF(1)_{\infty,\circ} \end{equation}
using the transitivity of the restriction and \eqref{eq:fockinftyande}. 

\smallskip

We prove the isomorphism \eqref{eq:fockcharO} using the explicit description of 
the action of the Chevalley generators on Fock spaces in terms of $\beta$-sets.
Fix $t \geqslant 0$ and $\mu$ a bipartition.
Let $\lambda = \varpi_t(\mu)$, i.e., we have $\tau_{2}(\lambda,0) = (\mu,\sigma_t)$. 
Given $i  = (-q)^j$, we want to compare the action of $f_i$ on $|\mu,Q_t\rangle_\infty$ 
and $[f_{-qi,\circ},f_{i,\circ}]$ on $|\lambda,1\rangle_{\infty,\circ}$. Let $\beta$ be the
set of $\beta$-numbers of the charged partition $(\lambda,0)$. 

\smallskip

First, we consider the right hand side of \eqref{eq:fockcharO}. We distinguish four cases:
\smallskip

\noindent \emph{Case 1.} If $j \notin \beta$, then one can add an $i$-node neither on 
$|\lambda,1\rangle$ nor on $f_{-qi,\circ}|\lambda,1\rangle$ (if the later is $\neq 0$), therefore $[f_{-qi,\circ},f_{i,\circ}]|\lambda,1\rangle=0.$

\smallskip

\noindent \emph{Case 2.} If $j,j+2 \in \beta$, then one can add a $(-qi)$-node neither on 
$f_{i,\circ}|\lambda,1\rangle$ nor on $|\lambda,1\rangle$, therefore
$[f_{-qi,\circ},f_{i,\circ}]|\lambda,1\rangle=0.$

\smallskip

\noindent \emph{Case 3.} If $j, j+1 \in \beta$ and $j+2 \notin \beta$, then 
$[f_{-qi,\circ},f_{i,\circ}] |\lambda,1\rangle = -f_{i,\circ} f_{-iq,\circ} |\lambda,1\rangle
= -|\lambda',1\rangle$ where $\lambda'$ is obtain by adding first an $(-qi)$-node, then
an $i$-node. The charged $\beta$-set of $(\lambda',0)$ is obtained from $\beta$ by changing
$j$ to $j+2$. 
 
\smallskip

\noindent \emph{Case 4.} If $j\in \beta$ and $j+1,j+2 \notin \beta$, then 
$[f_{-qi,\circ},f_{i,\circ}] |\lambda,1\rangle = f_{-qi,\circ} f_{i,\circ} |\lambda,1\rangle
= |\lambda'',1\rangle$ where $\lambda''$ is obtain by adding first an $i$-node, then
an $(-qi)$-node. The charged $\beta$-set of $(\lambda'',0)$ is again obtained from $\beta$ by changing
$j$ to $j+2$. 

\smallskip

Now, we consider the left hand side of \eqref{eq:fockcharO}. 
Let $\sigma_t = (\sigma_1,\sigma_2)$ be the $2$-core of $\lambda$. 
By definition, the bipartition $\mu$ is the unique bipartition whose charged $\beta$-sets
satisfy
 $$ \beta = \beta_0(\lambda) = (-1 + 2\beta_{\sigma_1}(\mu^1)) \sqcup 2\beta_{\sigma_2}(\mu^2).$$
By Proposition \ref{prop:tensorfock}, the Fock space $\bfF(Q_t)_\infty$ is 
identified with the tensor product of the level $1$ Fock spaces 
$\bfF(-q^{-1}( q^2)^{\sigma_1})_\infty \otimes \bfF((q^{2})^{\sigma_2})_\infty$. If $j$ is odd, then 
$i = (-q)^j = -q^{-1}(q^2)^{(j+1)/2}$ and the action of $f_i$ on
$|\mu_1,-q^{-1}( q^2)^{\sigma_1}\rangle_\infty$
corresponds to changing $(j+1)/2$ to $(j+1)/2+1$ in $\beta_{\sigma_1}(\mu_1)$. If $j$ is
even, then $i = (q^2)^{j/2}$ and $f_i$ acts on  $|\mu_2,(q^{2})^{\sigma_2}\rangle_\infty$ by changing $j/2$
to $j/2+1$. Using the previous equality of $\beta$-sets, this amount to changing 
$j$ to $j+2$ in $\beta$. This proves that the action of $f_i$ on $|\mu,Q_t\rangle_\infty$ 
and $[f_{-qi,\circ},f_{i,\circ}]$ on $|\lambda,1\rangle_{\infty,\circ}$ coincide up 
to a sign. 

\smallskip

Finally, it remains to see that in the case (3) the difference 
$a(\lambda')-a(\lambda)$ is odd, whereas in the case (4) the difference $a(\lambda'')-a(\lambda)$ is even.
This is a staightforward computation using the formula for the $a$-function given in \cite[4.4.2]{Lu84}.
\end{proof}

We can now define the action of $\partial$ on 
$\bigoplus_{t\in\bbN}\bfF(Q_t)_e$ and $[\scrU_\k]$.
To do that, consider the Fock space $\bfF(1)_{e,\circ}$ as a charged Fock space for the charge $s = 0$.
This endows $\bfF(1)_{e,\circ}$ with an integrable representation of $\frakg_{e,\circ}$ as in
\S \ref{subsec:xgrading}. 
Consequently, by Corollary \ref{cor:derivedaction-uk} and Lemma \ref{lem:isofocks}
we can endow $\bigoplus_{t\in\bbN} \bfF(Q_t)_e$, and therefore $[\scrU_\k]$, with
an integrable representation of $\frakg_e$ such that the maps
$$|\mu,Q_t\rangle_e\mapsto [V_\lambda] \longmapsto (-1)^{a(\lambda)}|\lambda,1\rangle_{e,\circ}$$
with $\lambda=\varpi_t(\mu)$ induce
$\frakg_e$-module isomorphisms
\begin{equation}\label{eq:xgrading}
\bigoplus_{t\in\bbN}\bfF(Q_t)_e\ \mathop{\longrightarrow}
  \limits^\sim\  [\scrU_\k] \ \mathop{\longrightarrow}
  \limits^\sim\ \mathrm{Res}_{\frakg_{e}}^{\frakg_{e,\circ}}
  \, \bfF(1)_{e,\circ}.\end{equation}

\begin{remark} The map $[V_\lambda] \longmapsto (-1)^{a(\lambda)}|\lambda,1\rangle_{e,\circ}$ also
endows $[\scrU_\k]$ with a structure of $\frakg_{e,\circ}$-module, but which does not
come from Harish-Chandra induction and restriction. On the other hand, this would be the case if we
were working with finite linear groups instead of unitary groups. In that case, the 
Fock space one should consider would be $\bfF(1)$ with the level $1$ action of
the Kac-Moody algebra associated with the quiver $q^\bbZ(q)$. Using the analogy between
$\mathrm{GL}_n(-q)$ and $\mathrm{GU}_n(q)$ which interwines (at the level of characters)
Harish-Chandra induction with the $2$-Harish-Chandra induction (see \cite[sec. 3]{BMM}), one could 
expect that the action of $\frakg_{e,\circ}$ on $[\scrU_\k]$ comes from some the truncation
of induction and restriction functors coming from Deligne-Lusztig varieties. Note that these
functors are no longer exact but only triangulated which explains the appearance of signs in the 
formulae for the action of $e_{i,\circ}$ and $f_{i,\circ}$ on the standard basis elements
$|\lambda,1\rangle_{e,\circ}$.
\end{remark}

\subsubsection{The $\frakg_e$-action on $\scrU_\k$}
From Corollary \ref{cor:derivedaction-uk} and \eqref{eq:xgrading}, we know that the
representation datum on $\scrU_\k$  yields an integrable representation of
$\frakg_e$ on $[\scrU_\k]$. In order to show that it endows $\scrU_\k$ with a structure
of categorical $\frakg_e$-representation,
it only remains to prove that there is a decomposition of the category into
weight spaces. To this end, we show that unipotent characters lying in the same $\ell$-block
have the same weight.

\begin{lemma}\label{lem:sameblocksameweight}
Let $\lambda$ and $\mu$ be partitions of $n$. If $V_\mu$,
$V_\lambda$ belong to the same block of $\k G_n$ then $[V_\lambda]$, 
$[V_\mu]$ have the same weight for the action of $\frakg_e$.
\end{lemma}

\begin{proof}
Recall that $[\scrU_\k]$ is isomorphic to $\mathrm{Res}_{\frakg_e}^{\frakg_{e,\circ}}\, \bfF(1)_{e,\circ}$ 
as a $\frakg_e$-module by \eqref{eq:xgrading}.
Here, the restriction is obtained through the map $\kappa : \frakg_e \longrightarrow 
\frakg_{e,\circ}$ defined in Lemma \ref{lem:kappa}. In particular, we have the following
equality of weights in $\X_e$
\begin{equation}\label{eq:equalityweights}
 \mathrm{wt}([V_\lambda]) = \kappa^* (\mathrm{wt}(|\lambda,1\rangle_{e,\circ})).
\end{equation}
The indecomposable modules $V_\lambda$, $V_\mu$ lie in the same block
of $\k G_n$ if and only if the unipotent character $E_\lambda$, $E_\mu$
lie in the same $\ell$-block. Now, by 
Propositions \ref{prop:Delta},  \ref{prop:lblocks} this is equivalent to the 
weights of $|\lambda,1\rangle_{e,\circ}$ and $|\mu,1\rangle_{e,\circ}$
to be equal.
The lemma follows from \eqref{eq:equalityweights}.
\end{proof}

From the lemma we deduce that the classes of the simple modules $[D_\lambda]$ are
also weight vectors. Indeed, they are linear combination of $[V_\mu]$'s in the 
same block, therefore with the same weight. Given $\omega \in \X_e$, we define
$\scrU_{\k,\omega}$ to be the Serre subcategory of $\scrU_\k$ generated by the
simple modules $D_\lambda$ such that $[D_\mu]$ has weight $\omega$. From Lemma
\ref{lem:sameblocksameweight} we deduce that 
  $$ \scrU_\k = \bigoplus_{\omega \in \X_e} \scrU_{\k,\omega}$$
and we obtain our second and main categorification theorem for the unipotent representations
of finite unitary groups.

\begin{theorem}\label{thm:charl}
The representation datum $(E,F,T,X)$ associated with Harish-Chandra induction and restriction and the decomposition $\scrU_\k = \bigoplus_{\omega \in \X_e} \scrU_{\k,\omega}$
yields a categorical representation of $\frakg_e$ on $\scrU_\k$. Furthermore, the map
$ |\mu,Q_t\rangle_e \longmapsto [V_{\varpi_t(\mu)}]$ induces a $\frakg'_e$-module
isomorphism $\bigoplus_{t \in \bbN} \bfF(Q_t)_e \simto [\scrU_\k]$.
\qed
\end{theorem}

\subsection{Derived equivalences of blocks of $\scrU_\k$}\label{sec:broue}

In this section we apply the categorical techniques of \cite{CR} to produce some derived
equivalences between blocks of $\scrU_\k$. In the linear prime case (when $e$ is even),
we use the existence of \emph{good blocks} \cite{Li12} to deduce Brou\'e's abelian defect group
conjecture for finite unitary groups.

\smallskip

Recall that $d$, $e$ and $f$ denote respectively the order of $q^2$, $-q$ and $q$ modulo 
$\ell$.

\subsubsection{Characterization of the blocks of $\scrU_\k$}
\label{subsec:block=weight}
Recall that we proved in Lemma \ref{lem:sameblocksameweight} that each weight space of
$\scrU_\k$ is a union of
blocks of $\scrU_\k$ (or equivalently unipotent blocks of $\k G\mod$). Here we investigate
which block can occur in a given weight space. More precisely, we show that each weight space
(resp. each weight space on a Harish-Chandra series) is indecomposable when $e$ is odd
(resp. when $e$ is even).

\smallskip

Throughout this section, we will denote by $\omega_{\lambda,\circ}$ the weight of 
$|\lambda,1\rangle_{e,\circ}$ for the action of $\frakg_{e,\circ}$ on 
$\bfF(1)_{e,\circ}$ and by $\omega_{\lambda}$ the weight of $[V_\lambda]$
for the action of $\frakg_e$ on $[\scrU_\k]$. Recall from \eqref{formA}, \eqref{formB} that
\begin{equation}\label{eq:weightcirc}
\begin{aligned}
 \omega_{\lambda,\circ} = &\, t_{\pi_{\lambda_{[e]}}}(\Lambda_{1,\circ}) - w_e(\lambda) \delta_\circ \\
 	= &\, \Lambda_{1,\circ} + \pi_{\lambda_{[e]}} - 
	\Delta(\lambda_{[e]},1)\delta_{\circ} - w_e(\lambda) \delta_\circ
\end{aligned}
\end{equation}
where $\lambda_{[e]}$ is the $e$-core of $\lambda$, and $w_e(\lambda) =
| \lambda^{[e]}|$ is its $e$-weight. Using Proposition \ref{prop:lblocks}, this shows that
$\omega_{\lambda,\circ}$ characterizes the $\ell$-block in which $V_\lambda$ lies. In other
words, $\ell$-blocks correspond to weight spaces for the action of $\frakg_{e,\circ}$
on $\bfF(1)_{e,\circ}$. Since the equality \eqref{eq:equalityweights} gives
$$\omega_\lambda = \kappa^*(\omega_{\lambda,\circ}),$$
we are left with computing the different weights that can appear in
$(\kappa^*)^{-1}(\omega_\lambda)$. We shall start with the case where $e$ is odd.

\begin{lemma}\label{lem:weighteodd}
Assume $e$ is odd. Let $\lambda$, $\nu$ be partitions. Then
$\omega_{\lambda,\circ} = \omega_{\nu,\circ}$ if and only if
$\omega_{\lambda} = \omega_{\nu}$.
\end{lemma}

\begin{proof}
Since $\omega_\lambda = \kappa^*(\omega_{\lambda,\circ})$, it is enough to
show that $\kappa^* : \X_{e,\circ} \longrightarrow \X_e$ is injective.
With $e$ being odd, we have $2\,\alpha^\vee_{i,\circ} = \sum_j (-1)^j\kappa(\alpha^\vee_{(-q)^j i})$ for each $i$,
hence $2\Q^\vee_{e,\circ}\subseteq\kappa(\Q^\vee_e).$
If $\kappa^*(\alpha_\circ)=0$ then
$0=\langle\kappa(\X_e^\vee),\alpha_\circ\rangle_{e,\circ}=\langle\kappa(\Q_e^\vee)+\bbZ\partial_\circ,\alpha_\circ\rangle_{e,\circ}.$
We deduce that $\langle\X_{e,\circ}^\vee,\alpha_\circ\rangle_{e,\circ}=0,$ hence $\alpha_\circ=0$.
\end{proof}

We now assume that $e$ is even. Then $\kappa^*$ is no longer injective, therefore
weight spaces of $\scrU_\k$ might contain several blocks in general. However, 
one can show that $\kappa^*$ is injective on weights coming from partitions
with the same $2$-core.

\begin{lemma}\label{lem:weighteeven}
Assume $e$ is even. Let $\lambda$, $\mu$ be partitions. Then 
$\omega_{\lambda,\circ} = \omega_{\nu,\circ}$  if and only if 
$\omega_{\lambda} = \omega_{\nu}$ and $\lambda_{[2]} = \nu_{[2]}$.
\end{lemma}

\begin{proof} If $\omega_{\lambda,\circ} = \omega_{\nu,\circ}$ then 
the partitions $\lambda$ and $\nu$ have the same $e$-core. With
$e$ being even, they also have the same $2$-core.

\smallskip

We now prove the converse. 
An easy computation shows that the kernel of $\kappa^*$ is spanned by the
weight $\sum_{j=0}^{e-1} (-1)^j \Lambda_{(-q)^j,\circ}$. In particular, if
$\omega_\lambda = \omega_\nu$ then $\omega_{\lambda,\circ} - 
\omega_{\nu,\circ}$ lies in $\P_{e,\circ}$. Therefore by \eqref{eq:weightcirc}, this difference
equals $\pi_s = \sum_{i\in \scrI_e} (s_i-s_{-qi})\Lambda_{i,\circ}$ where
$s = \lambda_{[e]} - \nu_{[e]}$. 
Since $\kappa^*(\pi_s)=0$, we deduce from the formulas in Lemma \ref{lem:kappa} that 
$s_i = s_{q^2i}$ for all $i \in\scrI_e$. On the other hand, since 
$\lambda_{[2]} = \nu_{[2]}$ then by Lemma \ref{lem:ecore2core} below we must have
$\sum_{j \text{ odd}} s_{-q^j} = 0$ and $\sum_{j \text{ even}} s_{q^j} 
=0$. This proves that $s =0$, thus $\pi_s = 0$ and therefore 
$\omega_{\lambda,\circ} = \omega_{\nu,\circ}$. 
\end{proof}

Assume $e$ is even. Given $t\in\bbN$ and $\omega \in \X_e$, we define
$\scrU_{\k,t}$ (resp. $\scrU_{\k,t,\omega}$) to be the Serre
subcategory of $\scrU_{\k}$ generated by simple modules $D_\lambda$
where  $\lambda$ has a $2$-core equal to $\Delta_t = (t,t-1,\ldots,1)$ (resp. with 
in addition $\omega_\lambda = \omega$). With $e$ being even, any pair of
partitions with the same $e$-core have the same $2$-core, therefore
$\scrU_{\k,t}$ is a direct summand of $\scrU_\k$. 
From the previous lemmas and the characterization of blocks of 
$\scrU_\k$ (see Proposition \ref{prop:lblocks}) we get the following more precise result.

\begin{proposition} \label{prop:main-blocks}
Let $\omega\in\X_e$.
\begin{itemize}[leftmargin=8mm]
  \item[$\mathrm{(a)}$] If $e$ is odd, then the category $\scrU_{\k,\omega}$ is an indecomposable
  summand of $\scrU_\k$.
  \item[$\mathrm{(b)}$] If $e$ is even, then the category $\scrU_{\k,t,\omega}$ is an 
  indecomposable summand of $\scrU_{\k,t}$ and $\scrU_\k$ for all $t\in\bbN$.
  \qed
\end{itemize}
\end{proposition}

\begin{remark}
The indecomposable categories in {\rm (a)} and {\rm (b)} are both equal to the unique
block $B_{\nu,w}$ of $\scrU_\k$ such that $\omega=\omega_{\nu,\circ}-w\,\delta$.
We will call $\omega$ the \emph{degree} of the block $B_{\nu,w}$.
\end{remark}

When $e$ is even $[\scrU_{\k,t}]$ is spanned by the classes $[V_\lambda]$ where
the $2$-core of $\lambda$ is the triangular partition $\Delta_t$. In particular,
$[\scrU_{\k,t}]$ is stable by the action
of $[E_i]$ and $[F_i]$, and since $\scrU_{\k,t}$ is a direct summand of $\scrU_\k$, 
the category $\scrU_{\k,t}$ itself is stable by $E_i$ and $F_i$. This allows
to refine Theorem \ref{thm:charl} into the following result when $e$ is even.

\begin{corollary}\label{cor:cat-even}
Assume $e$ is even, and let $t \in \bbN$. 
The Harish-Chandra induction and restriction functor yield a representation of
$\frakg_e$ on $\scrU_{\k,t}$ which categorifies $\bfF(Q_t)_e$.
\qed
\end{corollary}

\subsubsection{Action of the affine Weyl group}
\label{sec:Weyl-action}
In this section we study the action of the affine Weyl group of $\frakg_e$ on
the weight spaces of $\scrU_{\k}$ and $\scrU_{\k,t}$ with a view  to
understanding the action on unipotent blocks. As above, the results will
depend on the parity of $e$.

\smallskip

We start with the case where $e$ is even. To avoid cumbersome notation we fix
$t$ and we write $Q = Q_t = (Q_1,Q_2)$. 
Recall from \eqref{sigma} that $\sigma  = (\Delta_t)_{[2]} = (\sigma_1,\sigma_2)$ is given by
 $$\sigma=\begin{cases}\big(-t/2,\,t/2\big)&\ \text{if}\ t\ \text{is\ even},\\
 \big((1+t)/2,-(1+t)/2\big)&\ \text{if}\ t\ \text{is\ odd},\end{cases}$$
so that we have $Q = (-q^{-1} q^{2\sigma_1},q^{2\sigma_2})$ by \eqref{Q}.

\smallskip

In \S\ref{subsec:g-e} we defined $\frakg_e$ as a subalgebra of $\widetilde\frakg_e:=\frakg_{e,1} \oplus \frakg_{e,2}$
where each $\frakg_{e,p}$ is isomorphic to $\widehat{\fraks\frakl}_{e/2}$. The derived
Lie algebras of $\frakg_e$, $\widetilde\frakg_e$ coincide, but the derivation $\partial$ of $\frakg_e$ is equal to $(\partial_1,\partial_2)
 = (\Lambda_{-q^{-1}}^\vee, \Lambda_1^\vee)$. Recall also that the Lie
algebra $\frakg_{e,1}$ is generated by $\partial_1=\Lambda_{-q^{-1}}^\vee$ and $e_i,f_i$ for $i \in (-q)^{-1+2\bbZ}$, 
whereas $\frakg_{e,2}$ is generated by $\partial_2=\Lambda_{1}^\vee$ and
$e_i,f_i$ for $i \in q^{2\bbZ}$. For each $p = 1,2$, the Fock space $\bfF(Q_p)_{e,p}$
can be endowed with an action of $\frakg_{e,p}$. To do so, in order to define the action of $\partial_p$, we must fix a charge,
see \S\ref{subsec:xgrading}. First, observe that the bijection $\phi:(-q)^{p-2+2\bbZ}\to q^{2\bbZ}$ such that $\phi(i)= i(-q)^{2-p}$
identifies the quiver associated with $\frakg_{e,p}$ with the quiver $q^{2\bbZ}(q^2)$, and that $\phi(Q_p)=q^{2\sigma_p}.$
Thus, we can fix the charge to be $\sigma_p$.
Consequently, there is a well-defined action of $\widetilde\frakg_{e}$ on $\bfF(Q_1)_{e,1} \otimes \bfF(Q_2)_{e,2}$. 

\begin{proposition}\label{prop:tensoreven}
Assume $e$ is even. Then the linear map such that
$|\mu,Q\rangle_{e} \longmapsto |\mu^1,Q_1\rangle_{e,1} \otimes 
|\mu^2,Q_2\rangle_{e,2}$ is an isomorphism of $\frakg_e$-modules
 $$ \bfF(Q)_e \ \mathop{\longrightarrow}\limits^\sim \ 
 \mathrm{Res}_{\frakg_e}^{\widetilde\frakg_{e}} \big( \bfF(Q_1)_{e,1}
 \otimes \bfF(Q_2)_{e,2}\big).$$  
\end{proposition}

\begin{proof}
The map is an isomorphism of $\frakg'_e$-modules 
by Proposition \ref{prop:tensorfock} (recall that the derived algebras of $\frakg_e$ and $\widetilde\frakg_{e}$ coincide).
Therefore we are left with comparing the action of the derivation $\partial$ on the elements
$|\mu,Q\rangle_e$ and $|\mu^1,Q_1\rangle_{e,1} \otimes |\mu^2,Q_2\rangle_{e,2}.$

\smallskip

First, we compute the action of $\partial$ on $|\mu,Q\rangle_e$.
Set $\lambda = \varpi_t(\mu)$. By \eqref{eq:xgrading},
the weight of $|\mu,Q\rangle_e$ is $\omega_\lambda$.
By Lemma \ref{lem:kappa} and \eqref{eq:weightlevel1}, we get
$$ \langle \partial, \omega_{\lambda} \rangle_e =
 \langle \partial,\kappa^*(\omega_{\lambda,\circ}) \rangle_e = 
 \langle \kappa(\partial),\omega_{\lambda,\circ} \rangle_{e,\circ} 
= \langle \partial_{\circ} ,\omega_{\lambda,\circ}\rangle_{e,\circ} = -n_{1}(\lambda,1).$$

\smallskip

Now, we compute the action of $\partial$ on $|\mu^1,Q_1\rangle_{e,1} \otimes |\mu^2,Q_2\rangle_{e,2}.$
The embedding $\frakg_e \hookrightarrow\widetilde \frakg_{e}$ sends $\partial$ to $(\partial_1,\partial_2)$,
therefore the $\delta$-part of the weight of $|\mu^1,Q_1\rangle_{e,1} \otimes 
|\mu^2,Q_2\rangle_{e,2}$ is 
$$-n_{-q^{-1}}(\mu^1,Q_1) - \Delta(\sigma_1,e/2) - n_{1}(\mu^2,Q_2) - 
\Delta(\sigma_2,e/2). $$  
Since $e$ is even, removing an $e$-hook to $\lambda$ amounts to removing an $e/2$-hook
to $\mu^1$ or to $\mu^2$. Since $(-q)$ has order $e$ and $q^2$ has order $e/2$, this 
has the effect of substracting $1$ to $n_{1}(\lambda,1)$ and to 
$n_{-q^{-1}}(\mu^1,Q_1)+ n_{1}(\mu^2,Q_2)$. Therefore to prove the equality between
$\delta$-weights we can assume that $\lambda$ is an $e$-core. In that case both
$\mu^1$ and $\mu^2$ are $e/2$-cores and we use the following technical lemma together
with Lemma \ref{lem:Delta} to conclude.\end{proof}

\begin{lemma}\label{lem:ecore2core}
Let $\mu = (\mu^1,\mu^2)$ be a bipartition and $\lambda = \varpi_t(\mu)$.
Let $s = (s_{i})_{i\in \scrI_e} = \lambda_{[e]}$, $s^{1} = (s_{-q^j})_{j\,\mathrm{odd}}$ and $s^{2} = (s_{q^j})_{j\, \mathrm{even}}$.
Then, we have
 \begin{itemize}[leftmargin=8mm]
  \item[$\mathrm{(a)}$] $\tau_{e/2}(\mu^p,\sigma_p) \in \scrP^{e/2}\times \{s^{p}\}$ for $p=1,2$,
  \item[$\mathrm{(b)}$]  $w_e(\lambda)=w_{e/2}(\mu^1)+w_{e/2}(\mu^2)$.
\end{itemize}
\qed
\end{lemma}

From Proposition \ref{prop:tensoreven} we deduce that the action of the affine Weyl group $W_e$ of
$\frakg_e$ on $\bfF(Q)_e$ can be read off from the action of $W_{e,1}\times W_{e,2}$
on $\bfF(Q_1)_{e,1} \otimes \bfF(Q_2)_{e,2}$. We have the following partition
of weights into orbits under $W_e$.

\begin{proposition}\label{prop:orbiteven}
Assume $e$ is even. Let $\lambda$, $\nu$ be partitions. The weights $\omega_\lambda$ and
$\omega_\nu$ are conjugate under $W_e$ if and only if $\lambda_{[2]} = \mu_{[2]}$ and
$w_{e}(\lambda) = w_e(\mu)$.
\end{proposition}

\begin{proof} Since each $[\scrU_{\k,t}]$ is stable by $\frakg_e$, the partitions
$\lambda$ and $\mu$ must have the same 2-core in order to be conjugate under $W_e$.  
We will therefore assume that $\lambda_{[2]} = \mu_{[2]} = \sigma_t$ and we denote
by $\mu$ and $\gamma$ the bipartitions such that $\lambda = \varpi_t(\mu)$ 
and $\nu = \varpi_t(\gamma)$.

\smallskip

We write $\X_{e,1} = \P_{e,1} \oplus \,\bbZ \,\delta_1$ and $\X_{e,2} = \P_{e,2} \oplus \,\bbZ \,\delta_2$.
Recall that $\X_e$ is the quotient of the lattice $\widetilde\X_e=\X_{e,1}\oplus \X_{e,2}$ by
the line spanned by $\delta_1 - \delta_2$. To avoid cumbersome notation, we
will write $\widetilde\omega_\mu$, $\widetilde\omega_{\gamma}$ for the weights of 
$|\mu^1,Q_1\rangle_{e,1} \otimes |\mu^2,Q_2\rangle_{e,2}$ and
$|\gamma^1,Q_1\rangle_{e,1} \otimes |\gamma^2,Q_2\rangle_{e,2}$ in
$\widetilde\X_{e}$ and $\omega_\mu$, $\omega_{\gamma}$
for their image in $\X_e$. 
Notice that, by Proposition \ref{prop:tensoreven}, the weights
$\omega_{\lambda}$, $\omega_{\nu}$ are conjugate under $W_e$ 
if and only if $\omega_\mu$, $\omega_{\gamma}$ are conjugate 
under $W_{e,1}\times W_{e,2}$. 

\smallskip

Now, by Proposition \ref{prop:Delta} we have
\begin{align*}
&\widetilde\omega_\mu\in W_{e,1}(\Lambda_{Q_1})+W_{e,2}(\Lambda_{Q_2})+w_{e/2}(\mu^1)\delta_1
+w_{e/2}(\mu^2)\delta_2,\\
&\widetilde\omega_\gamma\in W_{e,1}(\Lambda_{Q_1})+W_{e,2}(\Lambda_{Q_2})+w_{e/2}(\gamma^1)\delta_1
+w_{e/2}(\gamma^2)\delta_2.
\end{align*}
In particular, the weights $\widetilde\omega_\mu$, $\widetilde\omega_{\gamma}$ are conjugate 
under $W_{e,1}\times W_{e,2}$ if and only if $$w_{e/2}(\mu^1)\delta_1
+w_{e/2}(\mu^2)\delta_2 =  w_{e/2}(\gamma^1)\delta_1+w_{e/2}(\gamma^2)\delta_2$$
in $\widetilde\X_e$.
Now, if $w_e(\lambda) = w_e(\nu)$ then $w_{e/2}(\mu^1)+w_{e/2}(\mu^2) = w_{e/2}(\gamma^1)+
w_{e/2}(\gamma^2)$  by Lemma \ref{lem:ecore2core}, hence there is an integer $m$ such that
the weights $\widetilde\omega_\mu$, $\widetilde\omega_{\gamma}+m(\delta_1-\delta_2)$ are conjugate 
under $W_{e,1}\times W_{e,2}.$
We deduce that the weights $\omega_\mu$, $\omega_{\gamma}$ are conjugate 
under $W_{e,1}\times W_{e,2}$, hence $\omega_{\lambda}$, $\omega_{\nu}$ are conjugate under $W_e$.

\smallskip

Conversely, assume that $\omega_{\lambda}$, $\omega_{\nu}$ are conjugate under $W_e$ or, equivalently,
that $\omega_\mu$, $\omega_{\gamma}$
are conjugate under $W_{e,1}\times W_{e,2}$. This means that there exist
$w \in W_{e,1}\times W_{e,2}$ and $m \in \bbZ$ such that $w(\widetilde\omega_\mu) = 
\widetilde\omega_\gamma + m(\delta_1-\delta_2)$. If $m$ is negative, we can
add $-m$ $e/2$-hooks to $\mu^1$ and to $\gamma^2$ to obtain bipartitions
$\mu'$ and $\gamma'$ such that $\widetilde\omega_{\mu'} = \widetilde\omega_{\mu} -m\delta_1$ 
and $\widetilde\omega_{\gamma'} = \widetilde\omega_{\gamma} -m\delta_2$. Since $w$ acts
trivially on $\delta_1$ and $\delta_2$, we now have $w(\widetilde\omega_{\mu'}) = 
\widetilde\omega_{\gamma'}$. By Proposition \ref{prop:Delta}, we deduce that 
$w_{e/2}((\mu')^1) = w_{e/2}((\gamma')^1)$ and $w_{e/2}((\mu')^2) = w_{e/2}((\gamma')^2)$. In particular, we have
$w_{e/2}((\mu')^1) + w_{e/2}((\mu')^2) = w_{e/2}((\gamma')^1) + 
w_{e/2}((\gamma')^2)$, and the same equality holds for $\mu$ and $\gamma$.
This proves that $w_e(\lambda)= w_e(\nu)$ by Lemma \ref{lem:ecore2core}. The same argument applies when
$m$ is nonnegative by adding hooks to $\mu^2$ and $\gamma_1$.
\end{proof}

We now turn to the case where $e$ is odd. Let $\rho_\circ = \sum \Lambda_{i,\circ}^\cl
\in  \Q_{e,\circ}^\cl$ be the sum of the classical fundamental weights.
For each $\alpha_\circ\in\X_{e,\circ}$ and $w\in W_{e,\circ}$, we abbreviate
$w\bullet \alpha_\circ=w(\alpha_\circ+\rho_\circ)-\rho_\circ$.

\begin{proposition}\label{prop:orbitodd}
Assume $e$ is odd. Let $\lambda$, $\nu$ be partitions. The weights $\omega_\lambda$
and $\omega_\nu$ are conjugate under $W_e$ if and only if  $\lambda_{[e]}\in
\frakS_{\scrI_e}\bullet\nu_{[e]}+2\bbZ^e$ and $w_e(\lambda)=w_e(\nu)$.
\end{proposition}

\begin{proof}
Let $W_{e,\circ}^{(2)}$ be the subgroup of $W_{e,\circ}$ which is equal to
$\frakS_{I_e}\ltimes 2\Q^\cl_{e,\circ}$
under the isomorphism $W_{e,\circ}=\frakS_{I_e}\ltimes\Q^\cl_{e,\circ}$.
Similarly, let $\frakg^{(2)}_{e,\circ}$ be the Lie subalgebra of $\frakg_{e,\circ}$ which is equal to
$ \fraks\frakl_e(\bbC) \otimes \bbC[t^2,t^{-2}] \oplus
  \bbC c \oplus \bbC \partial$ under the isomorphism
$\frakg_{e,\circ}=\widehat{\fraks\frakl}_e$ and
let $\pi$ be the Lie algebra automorphism of $\frakg_{e,\circ}$ which takes $x\otimes t^n$ to $x\otimes t^{n-\height(\alpha)}$
for each $\alpha$-root vector $x\in\fraks\frakl_e(\bbC)$. 
The embedding $\kappa \,:\,\frakg_e\to\frakg_{e,\circ}$ defined in Lemma \ref{lem:kappa}
yields an embedding $\kappa \,:\,W_e\to W_{e,\circ}$.

\begin{lemma} We have
 \begin{itemize}[leftmargin=8mm]
  \item[$\mathrm{(a)}$]
  $\kappa(\frakg_e)=\pi(\frakg^{(2)}_{e,\circ}),$
  \item[$\mathrm{(b)}$]
$\kappa(W_e)=(1,-\rho_\circ)\cdot W_{e,\circ}^{(2)}\cdot (1,\rho_\circ).$
\end{itemize} 
\end{lemma}

By Lemma \ref{lem:weighteodd}, we must characterize the pairs $\{\lambda,\nu\}$
such that $\omega_{\lambda,\circ}$ and $ \omega_{\nu,\circ}$ are $W_e$-conjugate in
$\X_{e,\circ}$. By \eqref{eq:weightcirc} and Proposition \ref{prop:Delta}, we deduce that
$\omega_{\lambda,\circ}$ and $\omega_{\nu,\circ}$ are $W_e$-conjugate 
if and only if $t_{\pi_s+\rho_\circ}(\Lambda_{1,\circ})-w_e(\lambda)\,\delta_\circ$ and
$t_{\pi_u+\rho_\circ}(\Lambda_{1,\circ})-w_e(\nu)\,\delta_\circ$ are
$W_{e,\circ}^{(2)}$-conjugate, where we abbreviate $s=\lambda_{[e]}$ and $u=\nu_{[e]}$.
We have $w_e(\lambda)=w_e(\nu)$ if and only if $t_{\pi_s+\rho_\circ}(\Lambda_{1,\circ})-
w_e(\lambda)\,\delta_\circ$ and $t_{\pi_u+\rho_\circ}(\Lambda_{1,\circ})-w_e(\nu)\,
\delta_\circ$ are $W_{e,\circ}$-conjugate. Hence $\omega_{\lambda,\circ}$ and
$\omega_{\nu,\circ}$ are $W_e$-conjugate if and only if $w_e(\lambda)=w_e(\nu)$ and 
$\pi_s\in \frakS_{I_e}\bullet\pi_u+2\Q_{e,\circ}^\cl$. The proposition is proved.
\end{proof}

\begin{proof}[Proof of the lemma]
It is enough to prove claim (a).
Let $\widetilde\kappa=\pi^{-1}\circ\kappa$.
We must prove that $\widetilde\kappa(\frakg_e)=\frakg^{(2)}_{e,\circ}.$
The inclusion $\widetilde\kappa(\frakg_e)\subseteq\frakg^{(2)}_{e,\circ}$ is obvious.
Let us concentrate on the reverse inclusion.

\smallskip

Let $\frakg^{(0)}_{e,\circ}$ be the Lie subalgebra of $\frakg_{e,\circ}$ which is equal to
$ \fraks\frakl_e(\bbC)$ under the isomorphism
$\frakg_{e,\circ}=\widehat{\fraks\frakl}_e.$ 
First, we consider the root vectors $\widetilde\kappa(e_i)$ of $\frakg_{e,\circ}$ as $i$ runs over the set $I_e$.
If $i=1$ or $-q^{-1}$ then $\widetilde\kappa(e_i)$ is of the form $x\otimes t^{3-e}$ where $x$
is a root vector of $\frakg^{(0)}_{e,\circ}$ whose weight has height $2-e$,
else it is of the form $y\otimes t^2$ where 
$y$ is a root vector of $\frakg^{(0)}_{e,\circ}$ whose weight has height $2$.
Now, it is easy to see that all root vectors of $\frakg^{(0)}_{e,\circ}$ associated with the negative simple roots
can be decomposed as the Lie bracket of one vector of the first type and $(e-3)/2$ vectors of the second type.
A similar result holds for the root vectors of $\frakg^{(0)}_{e,\circ}$ associated with the (positive) simple roots,
by considering the $f_i$'s instead of the $e_i$'s. This implies that $\frakg^{(0)}_{e,\circ}$
is contained into $\widetilde\kappa(\frakg_e)$. 

\smallskip

To conclude, observe that the Lie algebra
$\frakg^{(2)}_{e,\circ}$ is generated by $\frakg^{(0)}_{e,\circ}$ and any vector $\widetilde\kappa(e_i)$ of the second type
mentioned above.
\end{proof}

\subsubsection{Derived equivalences of blocks of $\scrU_\k$} \label{subsec:broue}
Now that we know the orbits of the affine Weyl group on the weight spaces of
$[\scrU_\k]$ (hence on the blocks of $\scrU_\k$), we can apply Proposition
\ref{prop:main-blocks} and the work of Chuang and Rouquier \cite{CR} to produce
derived equivalences between blocks of $\scrU_\k$ in the same $W_e$-orbit.

\smallskip

The following definition is taken from \cite{Li12} for $e$ even and from \cite{Li14}
for $e$ odd.

\begin{definition}
Let $\nu = (s_{(-q)^i}) \in \bbZ^{\scrI_e}$ be an $e$-core and $w \geqslant 0$. The unipotent block 
$B_{\nu, w}$ is a \emph{good block} if for every $i = 0, \ldots, e-2$ we have
$ s_{(-q)^i} \leqslant s_{(-q)^{i+2}} + w-1.$
\end{definition}

We will say that two unipotent blocks are \emph{$W_e$-conjugate} (or conjugate under
$W_e$) if the corresponding weight spaces in $\scrU_\k$ (or $\scrU_{\k,t}$ when $e$
is even) are $W_e$-conjugate.

\begin{lemma}
Any unipotent block is $W_e$-conjugate to a good block.
\end{lemma}

\begin{proof}
Let $\lambda$ be a partition and $w = w_e(\lambda)$. Assume first that $e$ is even. Let $\nu$ be any other partition such that $w_{e}(\lambda) = w_{e}(\nu)$. 
Write $\lambda_{[e]} = (s_{(-q)^i})$ and $\nu_{[e]} = (t_{(-q)^i})$. 
By Lemma \ref{lem:ecore2core}, the partitions $\lambda$, $\mu$
have the same $2$-core if and only if $\sum s_{q^{2i}} -\sum s_{-q^{2i+1}} = \sum t_{q^{2i}}-
\sum t_{-q^{2i+1}}$.  In such case,  Proposition \ref{prop:orbiteven} ensures that the degrees of the
blocks $B_{\lambda_{[e]}, w}$ and $B_{\nu_{[e]}, w}$ are $W_e$-conjugate. If one chooses for
$\nu$ the partition with $e$-core $t_{(-q)^i} = s_{(-q)^i} + \lfloor i/2 \rfloor (w-1)$, then
 $B_{\nu_{[e]}, w}$ is a good block.

\smallskip

If $e$ is odd, let $\nu$ be any partition with $w_e(\nu) = w_e(\lambda)$. Let $\varepsilon_i = 0$
if $s_{q^{2i+2}}-s_{q^{2i}}$ is even, and $1$ otherwise. We choose for $\nu$ the partition with
$e$-core $\nu_{[e]} = (t_{q^i})$ such that $t_{q^{2i}} = s_{1}+2i (w-1) + \varepsilon_0 +
\cdots + \varepsilon_{i-1}$. Then $B_{\nu_{[e]},w}$ is a good block, and $t_{q^{2i}}-s_{q^{2i}}
= 2i(w-1)+ s_{1} - s_{q^2} + s_{q^2} - \cdots - s_{q^{2i}} + \varepsilon_0 + \cdots +
\varepsilon_{i-1}$ is even. We deduce from Proposition \ref{prop:orbitodd} that the the blocks
$B_{\lambda_{[e]}, w}$ and $B_{\nu_{[e]}, w}$ are $W_e$-conjugate.
\end{proof}

Together with the results of Livesey \cite{Li12} on the structure of good blocks in the linear prime
case, we deduce that Brou\'e's abelian defect group conjecture holds for unipotent blocks
when $e$ is even.

\begin{theorem}\label{thm:broue}
Assume $e$ is even. Let $B$ be a unipotent block of $G_n$ over $\k$ or $\scrO$,
and $D$ a defect group of $B$. If $D$ is abelian, then $B$ is derived equivalent to the Brauer
 correspondent of $B$ in $N_{G_n}(D)$.
\end{theorem}

\begin{proof} 
Let $B'$ be a good block which is $W_e$-conjugate to $B$. By \cite[thm.~6.4]{CR}, the blocks $B$ and
$B'$ are derived equivalent. Now it follows from \cite[Theorem 7.1]{Li12}
that any good block is derived equivalent to its Brauer correspondent.
\end{proof}

\begin{remark} Let $B$ be any block of a finite general unitary group. Then, by the Jordan
decomposition \cite[Th\'eor\`eme B]{BoRo03}, the block $B$ is  Morita equivalent to a direct
product of unipotent blocks of unitary groups  $\GU_m(q^r)$ with $r$ odd and linear groups
$\GL_n(q^s)$ with $s$ even. However, it is not yet proven that this Morita equivalence preserves
the local structure of the blocks. If it were the case, then Theorem \ref{thm:broue}
would hold for any block of $G_n$.
\end{remark}

\subsection{The crystals of $\scrU_K$ and $\scrU_\k$}\label{sec:isocrystals}
In this section we show how to compare the crystal of the categorical representations
$\scrU_K$ and $\scrU_\k$ (which are related to Harish-Chandra induction and restriction)
with the crystals of the Fock spaces related to $[\scrU_K]$ and $[\scrU_\k]$. 
This solves the main conjecture of Gerber-Hiss-Jacon \cite{GHJ} and gives a combinatorial way to
compute the (weak) Harish-Chandra branching graph and the Hecke algebras associated
to the weakly cuspidal unipotent modules.

\subsubsection{Crystals and Harish-Chandra series}
Recall that to any categorical representation one can associate a perfect basis, and
hence and abstract crystal (see Proposition \ref{prop:PBfromcategorification}).
In the previous section we constructed a categorical action on the categories
of unipotent representations over $K$ (denoted by $\scrU_K$) and over $\k$ (denoted
by $\scrU_k$). From these two categorical representations we get:
\begin{itemize}[leftmargin=8mm]
  \item[(a)] an abstract crystal $B(\scrU_K)=\big(\Irr(\scrU_K),\widetilde E_i,\widetilde F_i\big)$
  with the canonical labeling $\Irr(\scrU_K)=\{[E_\lambda]\,|\,\lambda\in\scrP\}$.
  Since the representation of $\frakg_\infty$ on $\scrU_K$ is the direct sum of its
  representation on the subcategories $\scrU_{K,t}$ for each $t\in\bbN$, we deduce that
  $B(\scrU_K)$ is the direct sum of the abstract crystals $B(KG,E_t)=\big(\Irr(KG, E_t),
  \widetilde E_i,\widetilde F_i\big)$,
  
  \item[(b)] an abstract crystal $B(\scrU_\k)=\big(\Irr(\scrU_\k),\widetilde E_i,\widetilde F_i\big)$
  with the canonical labeling $\Irr(\scrU_\k)=\{[D_\lambda]\,|\,\lambda\in\scrP\}.$
\end{itemize}
By contruction of $E$ and $F$, these crystals are related to the (weak) Harish-Chandra 
series, as stated in the following proposition. 

\begin{proposition}\label{prop:wHC} Let $R = K$ or $\k$, and 
  $\scrI = (-q)^{\bbZ} \subset R^\times$.
  \begin{itemize}[leftmargin=8mm]
    \item[$\mathrm{(a)}$]
   $D \in \Irr(\scrU_R)$ is weakly cuspidal if and only if $\widetilde E_i D =0$
    for all $i \in \scrI$.
    \item[$\mathrm{(b)}$] If $M,N \in \Irr(\scrU_R)$, then $N$ appears in the head of $F(M)$ if and only if
    there exists $i \in \scrI$ such that $N \simeq \widetilde F_i M$.
    \item[$\mathrm{(c)}$] If $D \in \Irr(\scrU_R)$ is weakly cuspidal, then 
    $$\wIrr(RG,D)= \{\widetilde F_{i_1}\cdots\widetilde F_{i_m}(D)\,\mid\,m\in\bbN,\,
    i_1,\dots,i_m\in \scrI\}.$$
  \end{itemize}
  In other words, the (uncolored) crystal graph associated with $B(\scrU_R)$ coincides with the weak
  Harish-Chandra branching graph, and its connected components with the weak Harish-Chandra series.
  \end{proposition}
\begin{proof}
Assertions (a) and (b) follow from the definition of $\widetilde E_i$ and $\widetilde F_i$ 
and the fact that $E= \bigoplus E_i$ and $F= \bigoplus F_i$. 

\smallskip
Let $D$ be a cuspidal unipotent $RG_r$-module. Since the $F_i$'s are exact functors, we have $F_{i_1}\cdots F_{i_m}(D) \twoheadrightarrow  \widetilde F_{i_1}\cdots\widetilde F_{i_m}(D)$ which shows that  $\widetilde F_{i_1}\cdots\widetilde F_{i_m}(D) \in \wIrr(\k G_n,D)$. Therefore to prove (c) it is enough to
show that given $n = r+2m$ we have
 $$\wIrr(\k G_n,D) \subset \{\widetilde F_{i_1}\cdots\widetilde F_{i_m}(D)\,\mid\,i_1,
 \dots,i_m\in \scrI\}.$$
We argue by induction on $m$. Let $M\in\wIrr(\k G_{n+2},D)$. Then $F^{m+1}(D)$ maps onto $M$.
Thus there is a vertex $i\in \scrI$ such that $F_iF^{m}(D)$ maps onto $M$. Since $F_i$ is exact
we deduce that there is a composition factor $N$ of $F^m(D)$ such that $F_iN$ maps onto $M$, \emph{i.e.},
 such that $\widetilde F_iN=M$. The module $N$ lies in a weak Harish-Chandra series 
 $\wIrr(RG_{n},E)$ for some weakly cuspidal module $E$. Since $F$ and therefore
 $F_i$ preserves the series, we must have $E \simeq D$ (using Proposition \ref{prop:partitionwcusp}).
 By induction on $m$, one can write 
 $N = \widetilde F_{i_1}\cdots\widetilde F_{i_m}(D)$, hence $M = \widetilde F_i N = 
 \widetilde F_i \widetilde F_{i_1}\cdots\widetilde F_{i_m}(D)$.
\end{proof}

Note that by \cite{L77}, the ordinary Harish-Chandra series and weak Harish-Chandra series
on unipotent modules coincide when $R = K$.

\subsubsection{Comparison of the crystals}
Throughout this section we will assume that $e$, the order of $(-q)$ modulo $\ell$ is odd.
Recall from \S\ref{subsec:derivedaction-uk} that $[\scrU_\k]$ is isomorphic to a direct
sum of Fock spaces $\bfF(Q_t)_e$ where $t$ runs over $\bbN$ and $Q_t = (Q_1,Q_2)$ is defined by
$$Q_t=\begin{cases}\big((-q)^{-1-t}\,,\,(-q)^t\big)&\ \text{if}\ t\ \text{is\ even},\\
  \big((-q)^t\,,\,(-q)^{-1-t}\big)&\ \text{if}\ t\ \text{is\ odd.}\end{cases}$$
To any charged Fock space one can associate an abstract crystal, see \S \ref{subsec:crystalfock}. 
We now show how to choose the charge for each $\bfF(Q_t)_e$ so that the crystal
will coincide with the Harish-Chandra branching graph. 
Recall from \eqref{sigma} that 
 $$\sigma_t=\begin{cases}\big(-t/2,\,t/2\big)&\ \text{if}\ t\ \text{is\ even},\\
 \big((1+t)/2,-(1+t)/2\big)&\ \text{if}\ t\ \text{is\ odd}.\end{cases}$$
Now, we define 
\begin{equation}\label{eq:s}
s_t=(s_1,s_2)=-{1\over 2}(e+1,0)+\sigma_t,
\end{equation}
so that with the assumption
on $e$ we have $Q_p=q^{2s_p}$ for each $p=1,2$. In other words $s_t$ is a charge for
$\bfF(Q_t)_e$ with respect to $q^2$. We denote by $B(s_t)_e$ the corresponding
abstract crystal of $\bfF(Q_t)_e$, with the canonical labeling
$B(s_t)_e=\{b(\mu,s_t)\,|\,\mu\in\scrP^2,\,t\in\bbN\}.$
Finally, we set $B_e=\bigsqcup_{\,t\in\bbN}B(s_t)_e$.

\smallskip 

We can now prove our main theorem, which compares the crystal $B_e$ with the
abstract crystal coming from the categorical action of $\widehat \fraksl_e$ on $\scrU_\k$.

\begin{theorem}\label{thm:wHC}
Assume that $e$ is odd. The map $b(\mu,s_{t})\mapsto [D_{\varpi_t(\mu)}]$ is a
crystal isomorphism $B_e\simto B(\scrU_\k)$.
\end{theorem}

\begin{proof}
Recall from Proposition \ref{prop:PBfromcategorification} that the set
$B(\scrU_\k)=\{[D_\lambda]\,|\,\lambda\in\scrP\}$ is a perfect basis of $[\scrU_\k]$.
Further, the discussion in \S \ref{subsec:crystalfock} implies that the family
$\bfB^\vee_{e}:=\bigoplus_{t\in\bbN}\bfB^\vee(s_t)$ is a perfect basis of $\bigoplus_{t\in\bbN}\bfF({Q_t})_e$ which comes with a canonical labeling
$\bfB^\vee_{e}=\{\bfb^\vee(\mu,s_t)\,|\,\mu\in\scrP^2,\,t\in\bbN\}$.
Finally, by Proposition \ref{prop:Uglov} the map $\bfb^\vee(\mu,s_t)\mapsto b(\mu,s_t)$ is a crystal isomorphism $\bfB^\vee_e\simto B_{e}$.

\smallskip

Hence, we must prove that the bijection $\varphi:B(\scrU_\k)\to \bfB^\vee_e$ such that $[D_{\varpi_t(\mu)}]\mapsto \bfb^\vee(\mu,s_{t})$ is a
crystal isomorphism. Since both bases are perfect,
by Proposition \ref{prop:perfect-bases}, it is enough to prove that there is a partial order $\succ$ on the set of partitions such that,
under the $\bbC$-linear isomorphism $\bigoplus_{t\in\bbN}\bfF(Q_t)_e\simto[\scrU_\k]$ in
\eqref{eq:xgrading} which maps $|\mu,Q_t\rangle_e$ to $[V_{\varpi_t(\mu)}]$ for each $\mu,$ $t$, we have
\begin{align*}
\bfb^\vee(\mu,s_{t})=\varphi([D_{\varpi_t(\mu)}])\in [D_{\varpi_t(\mu)}]+\sum_{\varpi_u(\gamma)\succ\varpi_t(\mu)}\bbC\,[D_{\varpi_u(\gamma)}].
\end{align*}

\smallskip

To do so, we use the main result of \cite{RSVV} to show that one can take $\succ$ to
be the dominance order. Given a positive integer $d$, and elements $\kappa, s_1,\dots, s_l\in\bbC$,
let $\scrO^{s,\kappa}\{d\}$ be the category O of the rational double affine Hecke algebra of type $G(l,1,n)$ with
parameters $\kappa$ and $s=(s_1,\dots,s_l)$. The parameters are normalized as in \cite[sec.~6.2.1]{RSVV}.
The category $\scrO^{s,\kappa}\{d\}$ is a highest weight category.
We will abbreviate $\scrO^{s,\kappa}=\bigoplus_{d\in\bbN}\scrO^{s,\kappa}\{d\}$, where
$\scrO^{s,\kappa}\{0\}$ is the category of finite dimensional vector spaces.
Let $\Delta(\mu)^{s,\kappa}$ and $T(\mu)^{s,\kappa}$ be the standard and the tilting module with lowest weight $\mu\in\scrP^l$.
To characterize the highest weight structure on $\scrO^{s,\kappa}$ we must fix a partial order $\geqslant_{s,\kappa}$ on the set of $l$-partitions such that
$$(T(\mu)^{s,\kappa}:\Delta(\gamma)^{s,\kappa})\neq 0\Rightarrow \mu\geqslant_{s,\kappa}\gamma.$$
There are several choices for this order.
From now on we will assume that $s\in\bbZ^l$ and $\kappa=-e$ with $e$ a positive integer. 
Then, we choose the same order as in \cite[sec.~6.2.2]{RSVV}, following an idea of \cite{DG}.
More precisely,  recall the combinatorics of $l$-partitions in \S \ref{sec:combinatorics}.
If $A$, $B$ are boxes of Young diagrams of $l$-partitions, then we write $A>_sB$ if we have
$\ct^s(A)<\ct^s(B)$ or if $\ct^s(A)=\ct^s(B)$ and $p(A)>p(B)$.
We'll abbreviate $\widetilde\ct(A) = \ct^s(A) -ep(A)/l.$
Then, we define the partial order on $\scrP^l$ by setting
$\mu\geqslant_{s,-e}\gamma$ if and only if there are orderings $Y(\mu)=\{A_n\}$ and
$Y(\gamma)=\{B_n\}$ such that $A_n\geqslant_sB_n$ for all $n$.  
Note that, since $-e < 0$, if $\mu\geqslant_{s,-e}\gamma$
and $Y(\mu)=\{A_n\}$ and $Y(\gamma)=\{B_n\}$ are as above, then we have
$\widetilde\ct(A_n) \leqslant \widetilde\ct(B_n)$ for all $n$.

\smallskip

For each $p=1,2,\dots,l$ we set $Q_p=u^s$
where $u=\exp(2\sqrt{-1}\pi/e)$. We consider the charged Fock space $\bfF(s)=(\bfF(Q),s)$ of $\widehat{\fraks\frakl}_e$.
It is equipped with the Uglov canonical basis $\bfB^+(s)=\{\bfb^+(\mu,s)\,|\,\mu\in\scrP^l\}$.
Now, set $s^\star=(-s_l,\cdots,-s_2,-s_1)$ and
$\mu^\star=({}^t\mu^l,\dots,{}^t\mu^2,{}^t\mu^1)$ where  ${}^t\lambda$ denotes the transposed partition of $\lambda$.
By \cite[thm.~7.3]{RSVV} there is a $\bbC$-linear isomorphism $[\scrO^{s^\star,-e}]\simto\bfF(Q)$ such that
$$[\Delta(\gamma^\star)^{s^\star,-e}]\mapsto|\gamma,Q\rangle,\quad [T(\mu^\star)^{s^\star,-e}]\mapsto \bfb^+(\mu,s)$$
from which we deduce that
$$\bfb^+(\mu,s)\in |\mu,Q\rangle+\sum_{\mu^\star>_{s^\star,-e}\gamma^\star}\bbC\,|\gamma,Q\rangle.
$$
Using the transpose and the fact that $\gamma^\star>_{s^\star,-e}\mu^\star$
if and only if $\mu>_{s,-e}\gamma$ we obtain
\begin{equation}\label{eq:unitriangular}
\bfb^\vee(\mu,s)\in |\mu,Q\rangle+\sum_{\mu>_{s,-e}\gamma}\bbC\,|\gamma,Q\rangle.
\end{equation}

Now, we consider the particular case where $e$ is odd and $s=s_t$ for some $t\in\bbN$
as in \eqref{eq:s}. The following holds.

\begin{lemma} Assume that the positive integer $e$ is odd. For all bipartitions $\mu$,
$\gamma$ we have $\mu>_{s_t,-e}\gamma\Rightarrow \varpi_t(\gamma)>\varpi_t(\mu),$
where $>$ is the dominance order.
\end{lemma}

\begin{proof} Let $\lambda=\varpi_t(\mu)$ and $\nu = \varpi_t(\gamma)$.
By construction, the charged $\beta$-sets of  $\lambda$ and $\mu$ are related by 
$$\beta_0(\lambda)=\{2u-1\,|\,u\in\beta_{\sigma_1}(\mu^1)\}\sqcup\{2u\,|\,u\in\beta_{\sigma_2}(\mu^2)\}.$$
In order to compare $\lambda$ and $\nu$ we will compare the numbers of elements
in the sets $\{A\in Y(\lambda)\,|\,\ct(A)\geqslant j\}$ and 
$\{B\in Y(\nu)\,|\,\ct(B)\geqslant j\}$ and invoke \cite[lem. 4.2]{DG}. 
To this end we proceed as in \cite[p.~814]{DG} by computing, for all integer $j$
\begin{align*}
\sharp\{A\in Y(\lambda)\,|\,\ct(A)\geqslant j\} 
&=\sum_{\substack{k\in\bbZ\\ k\geqslant j}}\sharp\{u\in\beta_0(\lambda)\,|\,u\geqslant k+1\}+a_j\\
&=\sum_{\substack{k\in\bbZ\\ 2k\geqslant j}}\big(\sharp\{u\in\beta_{\sigma_1}(\mu^1)\,|\,u\geqslant k+1\}
+\sharp\{u\in\beta_{\sigma_2}(\mu^2)\,|\,u\geqslant k+1\}\big)\\
&+\sum_{\substack{k\in\bbZ\\ 2k+1\geqslant j}}\big(\sharp\{u\in\beta_{\sigma_1}(\mu^1)\,|\,u\geqslant k+2\}
+\sharp\{u\in\beta_{\sigma_2}(\mu^2)\,|\,u\geqslant k+1\}\big)+a_j\\
&=\sum_{\substack{k\in\bbZ\\ 2k\geqslant j}}\sharp\{A\in Y(\mu^1)\,|\,\ct(A)=k-\sigma_1\}\\
&+\sum_{\substack{k\in\bbZ\\ 2k\geqslant j}}\sharp\{A\in Y(\mu^2)\,|\,\ct(A)=k-\sigma_2\}\\
&+\sum_{\substack{k\in\bbZ\\ 2k+1\geqslant j}}
\sharp\{A\in Y(\mu^1)\,|\,\ct(A)=k-\sigma_1+1\}\\
&+\sum_{\substack{k\in\bbZ\\ 2k+1\geqslant j}}\sharp\{A\in Y(\mu^2)\,|\,\ct(A)=k-\sigma_2\}+a_{j,t},
\end{align*}
where the integers $a_j$ and $a_{j,t}$ are constants which depend only on $j$ and $t$.
Next, assume that $e$ is odd. Recall from \eqref{eq:s} that $s_1=\sigma_1+d$ and $s_2=\sigma_2$
with $d=-(e+1)/2$. Now, we identify $Y(\mu)$ with $Y(\mu^1)\sqcup Y(\mu^2)$. We deduce that
\begin{align*}
\sharp\{A\in Y(\lambda)\,|\,\ct(A)\geqslant j\} 
&=\sharp\{A\in Y(\mu^1)\,|\,\ct^{s_1}(A)\geqslant \lceil j/2\rceil+d\}\\
&+\sharp\{A\in Y(\mu^2)\,|\,\ct^{s_2}(A)\geqslant \lceil j/2\rceil\}\\
&+\sharp\{A\in Y(\mu^1)\,|\,\ct^{s_1}(A)\geqslant \lceil (j-1)/2\rceil+d+1\}\\
&+\sharp\{A\in Y(\mu^2)\,|\,\ct^{s_2}(A)\geqslant \lceil (j-1)/2\rceil\}+a_{j,t}.
\end{align*}
In particular,  we get
\begin{align*}
\sharp\{A\in Y(\lambda)\,|\,\ct(A)\geqslant j\} 
&=\sharp\{A\in Y(\mu^1)\,|\,\widetilde\ct(A)\geqslant \lceil j/2\rceil-e-1/2\}\\
&+\sharp\{A\in Y(\mu^2)\,|\,\widetilde\ct(A)\geqslant \lceil j/2\rceil-e\}\\
&+\sharp\{A\in Y(\mu^1)\,|\,\widetilde\ct(A)\geqslant \lceil (j-1)/2\rceil-e+1/2\}\\
&+\sharp\{A\in Y(\mu^2)\,|\,\widetilde\ct(A)\geqslant \lceil (j-1)/2\rceil-e\}+a_{j,t}.
\end{align*}
Finally, since $e$ is odd and $\ct^{s_t}(A)$ is an  integer we can remove the symbols $\lceil\bullet\rceil$ and we obtain
\begin{align*}
\sharp\{A\in Y(\lambda)\,|\,\ct(A)\geqslant j\} 
&=\sharp\{A\in Y(\mu)\,|\,\widetilde\ct(A)\geqslant j/2-e-1/2\}\\
&+\sharp\{A\in Y(\mu)\,|\,\widetilde\ct(A)\geqslant  j/2-e\}+a_{j,t}.
\end{align*}
Now, assume that $\mu\geqslant_{s_t,-e}\gamma$. Then we can choose orderings $Y(\mu)=\{A_n\}$
and $Y(\gamma)=\{B_n\}$ such that $\widetilde\ct(A_n)\leqslant\widetilde\ct(B_n)$ for all $n$.
We deduce that for all integer $j$ we have
$$\sharp\{A\in Y(\nu)\,|\,\ct(A)\geqslant j\} 
\geqslant\sharp\{A\in Y(\lambda)\,|\,\ct(A)\geqslant j\},$$
from which we deduce that $\nu\geqslant\lambda$ (for the dominance order) by 
\cite[lem~4.2]{DG}. \end{proof}

\smallskip

\noindent \emph{Proof of Theorem \ref{thm:wHC}. Continued.}
From the lemma above together with \eqref{eq:unitriangular} we deduce that
\begin{align*}
\bfb^\vee(\mu,s_t)\in |\mu,Q_t\rangle_e+\sum_{\varpi_t(\gamma)>\varpi_t(\mu)}
\bbC\,|\gamma,Q_t\rangle_e.
\end{align*}
Hence, under the identification of $[\scrU_\k]$ with $\bigoplus_{t\in\bbN}\bfF({Q_t})_e$ we
obtain
\begin{align*}
\bfb^\vee(\mu,s_t)\in [V_{\varpi_t(\mu)}]+\sum_{\varpi_t(\gamma)>\varpi_t(\mu)}
\bbC\,[V_{\varpi_u(\gamma)}].
\end{align*}
On the other hand, by Proposition \ref{prop:unitriangularunip} we have
\begin{align*}[D_\lambda]\in [V_\lambda]+\sum_{\nu>\lambda}\bbC\,[V_\nu].\end{align*}
Therefore one can use Proposition \ref{prop:perfect-bases}, to deduce that the map
$[D_{\varpi_t(\mu)}]\mapsto \bfb^\vee(\mu,s_{t})$ is a
crystal isomorphism $B(\scrU_\k)\simto \bfB^\vee_e$. \end{proof}

As a direct consequence of this theorem and Proposition \ref{prop:wHC}, the unipotent
module $D_{\varpi_t(\mu)}$ in $\scrU_\k$ is weakly cuspidal if and only if
we have $\widetilde E_i(b(\mu,Q_{t}))=0$ for all $i\in \scrI_e$.

\begin{remark}
Our choice for the charge $s_t$ differs (by a sign) from the choice used by Gerber-Hiss-Jacon in \cite{GHJ}.
\end{remark}

\begin{remark}
When $e$ is even and $t\geqslant 0$, the category $\scrU_{\k,t}$ categorifies the level $2$ 
Fock space $\bfF(Q_t)_e$ by Corollary \ref{cor:cat-even}. It has a tensor product decomposition 
into level $1$ Fock spaces $\bfF(-q^{-1}q^{2\sigma_1})_e\otimes\bfF(q^{2\sigma_2})_e$. 
Here, we view $q$ as an element of $\k$.
One can therefore construct an abstract crystal of 
$\bfF(Q_t)_e$ out of the crystals of the charged level 1 Fock spaces $(\bfF(q^{2\sigma_1}),\sigma_1),$
$(\bfF(q^{2\sigma_2}),\sigma_2)$ under the identification 
$\bfF(Q_t)_e=\bfF(q^{2\sigma_1})\otimes\bfF(q^{2\sigma_2})$ in \S\ref{sec:Weyl-action}.
This crystal coincides with the
crystal $B(\scrU_{\k,t})$ coming from the categorification. 
\end{remark}

\begin{remark}
By Theorem \ref{thm:char0}, the category $\scrU_{K,t}$ categorifies the level $2$ 
Fock space $\bfF(Q_t)_\infty.$
It has a tensor product decomposition 
into level $1$ Fock spaces $\bfF(-q^{-1}q^{2\sigma_1})_\infty\otimes\bfF(q^{2\sigma_2})_\infty$.
Here, we view $q$ as an element of $K$.
One can therefore construct an abstract crystal of 
$\bfF(Q_t)_\infty$ out of the crystals of the charged level 1 Fock spaces
$(\bfF(q^{2\sigma_1}),\sigma_1),$
$(\bfF(q^{2\sigma_2}),\sigma_2)$
under the identification 
$\bfF(Q_t)_\infty=\bfF(q^{2\sigma_1})\otimes\bfF(q^{2\sigma_2})$ associated with the map
$(-q)^{p-2+2\bbZ}\to q^{2\bbZ}$ such that $i\mapsto i(-q)^{2-p}$ for each $p=1,2$.
This crystal coincides with the
crystal $B(\scrU_{K,t})$ coming from the categorification. 
\end{remark}

\subsubsection{Corollaries}
As a byproduct of our main theorem, we obtain a proof of the conjectures of
Gerber-Hiss-Jacon stated in \cite{GHJ}.

\begin{corollary}\label{cor:1}
Assume that $e$ is odd. The modules $D_\lambda$ and $D_\nu$ lie in the same weak
Harish-Chandra series if and only if the corresponding vertices of the abstract crystal $B_e$ belong to the same connected component.
In particular, if this holds then $\lambda$ and $\nu$ have the same $2$-core.
\end{corollary}

\begin{proof}
Let $\mu,\gamma$ be bipartitions and $t,u$ be non-negative integers such that
$\lambda=\varpi_t(\mu)$ and $\nu=\varpi_u(\gamma)$. Then by Proposition \ref{prop:wHC} and
Theorem \ref{thm:wHC}, the vertices $b(\mu,s_t)$ and $b(\gamma,s_u)$ of the abstract
crystal $B_e$ lie in the same connected component. This implies $t=u$.
\end{proof}

\begin{corollary}\label{cor:3}
Assume that $e$ is odd and that $\lambda$ is a partition of $n$.
If $D_\lambda$ is weakly cuspidal, then its $\ell$-block contains a cuspidal
simple $KG_n$-module (not necessarily unipotent). 
In particular, if this holds then
$\lambda$ is a $2$-core. 
\end{corollary}

\begin{proof}
It follows from the combination of \cite[sec.~5.5, thm.~7.6]{GHJ} and the crystal isomorphism given 
in Theorem \ref{thm:wHC}.
\end{proof}

Let $r > 0$, $m\geqslant 0$ and $n=r+2m$. Recall from \S\ref{subsec:hecke} that there is
$\k$-algebra homomorphism
$$\phi_{\k,m}\, :\, \bfH^{q^2}_{\k,m}\longrightarrow\End_{\k G_n}(F^m)^\op.$$
As in Theorem \ref{thm:cat}, one can show that the evaluation at a weakly cuspidal
module is a Hecke algebra of type $B_m$. 

\begin{corollary}\label{cor:2} 
Let $\lambda$ be a partition of $r>0$ such that $D_\lambda$ is weakly cuspidal.
Let $(t,t-1,\ldots,1)$ be the $e$-core of $\lambda$. Then for $m \geqslant 0$ and $n = r+2m$,
the map $\phi_{\k,m}$ factors to an algebra isomorphism
  $$\bfH^{Q_t,\,q^2}_{\k,m}\mathop{\longrightarrow}\limits^\sim \scrH(\k G_n,D_\lambda).$$
\end{corollary}

\begin{proof}
We first compute the weight of $[D_\lambda]$. By Corollary \ref{cor:3}, the $e$-core
of $\lambda$ is a 2-core $\Delta_t = (t,t-1,\ldots,1)$ for some $t \geqslant 0$. Then 
from Proposition \ref{prop:Delta} and \eqref{eq:equalityweights} we deduce that the projection
to $\P_e$ of the weights $\omega_\lambda$ and $\omega_{\Delta_t}$ are equal.
By Lemma \ref{lem:isofocks} and  \eqref{eq:xgrading}
we deduce that this projection equals $\Lambda_{Q_t}$. We can now invoke
Proposition \ref{prop:unicite} for $D_\lambda$: we have $E D_\lambda = 0$ since $D_\lambda$
is weakly cuspidal, and $\End_{\scrU_\k}(D_\lambda)=\k$ since $D_\lambda$ is simple
and $\k$ is a splitting field for $G_r$. This shows that $\phi_{\k,m}$ induces
the expected isomorphism.
\end{proof}

\begin{remark} Corollary \ref{cor:1} proves Conjecture 5.4 in \cite[sec.~5.4]{GHJ}.
Corollary \ref{cor:3} proves Conjecture 5.5 in \cite[sec.~5.4]{GHJ}.
Corollary \ref{cor:2} generalizes the argument in \cite[prop.~5.21]{G01},
in the particular case of the unitary group, and it proves the conjecture in 
\cite[\S 5]{GH}.
\end{remark}

\section{The representation of the Heisenberg algebra on $\scrU_\k$}\label{sec:heisenberg}
Every cuspidal $\k G_n$-module is weakly cuspidal. Therefore, every Harish-Chandra
series in $\Irr(\scrU_\k)$ is partitioned into weak Harish-Chandra series. Proposition
\ref{prop:wHC} and Theorem \ref{thm:wHC} yield a complete (combinatorial) description of
the partition of $\Irr(\scrU_\k)$ in weak Harish-Chandra series, which coincides with
the decomposition of $[\scrU_\k]$ for the action of $\frakg_e$. 
In this section we propose a conjectural generalization to Harish-Chandra series 
by endowing $\scrU_\k$ with an extra action of a Heisenberg algebra.

\smallskip

Throughout this section we will assume that $e$, the order of $-q$ modulo $\ell$ 
is odd, and hence equal to $d$, the order of $q^2$ modulo $\ell$.

\subsection{The $q$-Schur algebra}
Fix some arbitrary integers $m\geqslant n>0$. Assume that $R$ is one of the fields $K$ or $\k$,
where $(K,\scrO,\k)$ is an $\ell$-modular system as in \S\ref{sec:basics}. In particular
$q \in R^\times$. Let $v$ be a formal variable and $A=R[v,v^{-1}]$. Let $\bfU_A(\frakg\frakl_m)$
be Lusztig's divided power version over $A$ of the quantized enveloping algebra of
$\frakg\frakl_m$.

\begin{definition}
The \emph{$q$-Schur algebra} $\bfS^q_{R,m,n}$ over $R$ is the quotient algebra of 
$\bfU_A(\frakg\frakl_m)$ by the ideal generated by the element $v-q$ and the annihilator ideal of the $n$-fold tensor power of the natural module
of $\bfU_A(\frakg\frakl_m)$.
\end{definition}

Let $\Lambda^+(m)\subseteq\bbZ^m$ be the set of dominant weights. Since $m\geqslant n$ we have
a canonical inclusion of the set of partitions $\scrP_n$ of $n$ into $\Lambda^+(m)$.
For each dominant weight $\lambda$, the irreducible $\bfU_v(\frakg\frakl_m)$-module of highest 
weight $\lambda$ has an integral form $\bfL_A(\lambda)$ which is a module for
$\bfU_A(\frakg\frakl_m)$. If $\lambda$ is a partition of $n$, this module factors through the
quotient $\bfS^q_{R,m,n}$. This yields a complete set of non-isomorphic irreducible
$\bfS^q_{R,m,n}$-modules
$$\Irr(\bfS^q_{R,m,n})=\{L_R(\lambda)\,|\,\lambda\in\scrP_n\}.$$

There is a natural notion of tensor product of two $\bfU_A(\frakg\frakl_m)$-modules coming from 
the comultiplication of the Hopf algebra $\bfU_A(\frakg\frakl_m)$. For $n_1+n_2=n$, the tensor
factors through a bifunctor
$$\bullet\otimes\bullet\,:\,\bfS^q_{R,m,n_1}\mod\times\,\bfS^q_{R,m,n_2}\mod
\longrightarrow \bfS^q_{R,m,n}\mod.$$

\subsection{The cuspidal algebra} Recall that $GL_n$ denotes the finite linear
group $\mathrm{GL}_n(q^2)$ over the field with $q^2$ elements and that $R=K$ or $\k$.
Let $P_{R,n} := R GL_n/B_n$ be the $RGL_n$-module arising from the permutation representation
of $GL_n$ on the cosets of a Borel subgroup $B_n$. 

\begin{definition}
The \emph{unipotent cuspidal algebra} is the quotient algebra $\bfC_{R,n}$ of $RGL_n$ by
the annihilator ideal $I_{R,n}$ of $P_{R,n}$.
\end{definition}

The $R$-algebra $\bfC_{R,n}$ is actually a quotient algebra of the \emph{unipotent block} of
$RGL_n$, by which we mean the sum of blocks corresponding to the irreducible constituents of
the module $P_{R,n}$. Therefore, the pull-back of modules by the quotient map $RGL_n
\twoheadrightarrow \bfC_{R,n}$ gives a fully faithful functor $h_n\,:\,\bfC_{R,n}\mod
\to RGL_n\umod.$ We will view the category $\bfC_{R,n}\mod$ as a subcategory of $RGL_n\umod$.
Then, we have a canonical identification 
$$\Irr(\bfC_{R,n})=\Unip(RGL_n).$$

\smallskip

First, set $R=K$. The set of
unipotent characters of $KGL_n$ is $$\Unip(KGL_n)=\{\rho_\lambda\,|\,\lambda\in\scrP_n\},$$
where $\rho_\lambda$ is given by the following formula
$$\rho_\lambda= \frac{1}{|\frakS_n|} \sum_{w \in \frakS_n}
  \phi_\lambda(w) R_{\bfT_w}^{\bfG\bfL_n}(1)$$
and $\bfT_w$ runs over the set of the $GL_n$-conjugacy classes of $F_{q^2}$-stable maximal
tori in $\bfG\bfL_n$. Note that unlike the case of finite unitary groups, the class
function $\rho_\lambda$ is a true character (compare with the virtual character
$\chi_\lambda$ defined in \S \ref{subsec:uK}).

\smallskip

Now, consider the case $R=\k$. Let $d_{\scrO GL_n}:[K GL_n\mod]\to[\k GL_n\mod]$ be the decomposition map. 
Dipper showed in \cite{D1,D2} that there is a unique labeling 
$$\Unip(\k GL_n)=\{S_\lambda\,|\,\lambda\in\scrP_n\}$$
such that $d_{\scrO GL_n}(\rho_\lambda)=[S_\lambda]$ modulo $\bigoplus_{\mu>\lambda}\bbZ\,[S_\mu]$.

\smallskip

Finally, the Harish-Chandra induction relative to the subgroup $GL_{n_1}\times GL_{n_2}$
of $GL_n$ yields a bifunctor
$$\bullet\odot\bullet\,:\,\bfC_{R,n_1}\mod\times\bfC_{R,n_2}\mod
\longrightarrow \bfC_{R,n}\mod.$$
We have the following result, due to Takeuchi \cite{T96}. 
See \cite[sec.~3.5, thm.~4.2a]{BDK} for a formulation closer to ours.

\begin{proposition}[\cite{T96}]\label{prop:schur}
For each $m\geqslant n,$ there is an equivalence of abelian categories $\beta\,:\,
\bfS^{q^2}_{R,m,n}\mod\simto\bfC_{R,n}\mod$ such that
 \begin{itemize}[leftmargin=8mm]
  \item[$\mathrm{(a)}$] $\beta$ intertwines the bifunctors $\otimes$ and $\odot$,
  \item[$\mathrm{(b)}$] $\beta(L_K({}^t\lambda))=\rho_\lambda$ and $\beta(L_\k({}^t\lambda))=
  S_\lambda$ for each $\lambda\in\scrP_n$.
\qed
\end{itemize}
\end{proposition}

\subsection{Categorification of the Heisenberg operators}
Recall that $R=K$ or $\k$. Let $n=r+2m$ with $r,m>0$, and $L_{r,m} \simeq G_r \times GL_m$
be the corresponding standard Levi  as defined in \S \ref{sec:def-unitary}. For any module
$X$ in $R GL_{m}\umod$ we consider the following functors
\begin{align*}
&B_X\,:\,RG_r\umod\to RG_{n}\umod,\ M\mapsto R_{L_{r,m}}^{G_{n}}(M\otimes X),\\
&B^*_X\,:\,RG_{n}\umod\to RG_{r}\umod,\ M\mapsto 
\Hom_{R  GL_m}\big(X,{}^*\!R_{L_{r,m}}^{G_{n}}(M)\big).
\end{align*}
The functor $B^*_X$ is right adjoint to the functor $B_X$.
In the particular case where $X =R$ is the trivial module, we recover the functors
$F$ and $E$ as defined in \S\ref{sec:rep-datum}.

\smallskip

Now, let us consider the particular case $R=\k$. Assume that $I_{\k,m} X = 0$, so that $X$
can be viewed as an objet of $\bfC_{\k, m}\mod$. The inclusion functor $h_m\,:\,\bfC_{\k, m}\mod
\to \k GL_m\umod$ has a (left exact) right adjoint $h^!_m\,:\, \k GL_m\umod\to\bfC_{\k,m}\mod$
which takes a $\k GL_m$-module $M$ to the annihilator of the ideal $I_{\k,m}$ in $M$.
We have
$$B^*_X=\Hom_{\bfC_{\k,m}}\big(X,{}^*\!R_{L_{r,m}}^{G_{n}}(\bullet)\big):=\Hom_{\bfC_{\k,m}}\big(X,(\id\otimes h^!_m)\,{}^*\!R_{L_{r,m}}^{G_{n}}(\bullet)\big).$$
Since the algebra $\bfC_{\k,m}$ has a finite global dimension by Proposition \ref{prop:schur}, 
we may define the right derived functor
\begin{align*}\R \!B^*_X\,:\,D^b(\k G_{n}\umod)\to D^b(\k G_{r}\umod),\ M\mapsto 
\RHom_{\bfC_{\k,m}}\big(X,{}^*\!R_{L_{r,m}}^{G_{n}}(M)\big).
\end{align*}
Note that a module in $\k GL_m\umod$ may not have a finite global dimension,
and therefore $\R\! B_X^*$ might not exist if $X$ is not assumed to be annihilated
by $I_{\k,m}$. Using standard arguments, like tensor-hom adjunction, we can prove the following.

\begin{proposition}
For each module $X$ in $\bfC_{\k,m}\mod$ the following hold
 \begin{itemize}[leftmargin=8mm]
  \item[$\mathrm{(a)}$] $B_X$ is exact and extends to an exact functor 
$D^b(\k G_r\umod)\to D^b(\k G_{n}\umod)$, 
\item[$\mathrm{(b)}$] $(B_X,\R\! B^*_X)$ is an adjoint pair of triangulated functors,
\item[$\mathrm{(c)}$] if $X_1\in\bfC_{\k,m_1}\mod$, $X_2\in\bfC_{\k,m_2}\mod$, then there
are isomorphisms of functors\\
$B_{X_1}B_{X_2}\simeq B_{X_1\odot X_2}\simeq B_{X_2}B_{X_1}$ and
$\R\! B^*_{X_1}\R\! B^*_{X_2}\simeq \R\! B^*_{X_1\odot X_2}\simeq \R\! B^*_{X_2}\R\! B^*_{X_1}.$
\qed
\end{itemize}
\end{proposition}

\subsection{The modular Steinberg character}
Recall that $\rho_{(1^m)}$ is the Steinberg character of $KGL_m$ while $\rho_{(m)}$ is the
trivial one. We deduce that $S_{(m)}$ is the trivial $\k GL_m$-module, while $S_{(1^m)}$ is
the top of a modular reduction of the Steinberg character $\chi_{(1^m)}$ called the
\emph{modular Steinberg character} of $\k GL_m$. We will write $St_m=S_{(1^m)}$. The
$\k GL_m$-module $St_m$ is cuspidal if and only if $m=1$ or $m=e\ell^j$ for some $j\in\bbN$
If this holds, it is the only cuspidal module in $\Unip(\k GL_m)$. See \cite{GHM1} for details.

\smallskip

For any partition $\lambda$ of $m$ of the form $\lambda=(1^{(m_{-1})},e^{(m_0)},(e\ell)^{(m_1)},(e\ell^2)^{(m_2)},\dots)$ 
in exponential notation, we set 
$$GL_\lambda=(GL_1)^{m_{-1}}\times\prod_{j\geqslant 0}(GL_{e\ell^j})^{m_j},\qquad
St_\lambda=(St_1)^{\otimes m_{-1}}\otimes\bigotimes_{j\geqslant 0}(St_{e\ell^j})^{\otimes m_j}.$$
Then the unipotent cuspidal pairs of $\k GL_m$ are the pairs $(GL_\lambda,St_\lambda)$ where
$\lambda$ runs over the set of all partitions of $m$ of the form above. For any partition $\lambda\in\scrP_m$, we set
$$X_\lambda=R^{GL_m}_{GL_\lambda}(St_\lambda)\in\Unip(\k GL_m).$$
Note that $X_\lambda$ can be seen as an object
of the category $\bfC_{\k, m}\mod$ since it is annihilated by $I_{\k,m}$.
In particular, for each $j\in\bbN$ we have
$$S_{(1^{e\ell^j})}=St_{e\ell^j}=St_{(e\ell^j)}=X_{(e\ell^j)}.$$

\smallskip

The functors $F,$ $E$ in the representation data on $\scrU_\k$ as defined in 
\S\ref{sec:rep-datum} are given by $F=B_{X_{(1)}}$ and $E=B^*_{X_{(1)}}$. The following
proposition is easy to prove.

\begin{proposition} For a given unipotent $\k G_n$-module $M$, the following conditions
are equivalent
 \begin{itemize}[leftmargin=8mm]
  \item[$\mathrm{(a)}$] $M$ is cuspidal,
  \item[$\mathrm{(b)}$] $B^*_{X}(M)=0$ for all $X\in\bfC_{\k,m}\mod$ and all $m>0$,
  \item[$\mathrm{(c)}$] $E_i(M)=B^*_{S_{(1^{e\ell^j})}}(M)=0$ for all $i\in I_e$, $j\in\bbN$,
  $j\neq 0$.
\qed
\end{itemize}
\end{proposition}

\subsection{The Heisenberg algebra and the Fock space}\label{sec:Heisenberg}
The \emph{Heisenberg algebra} is the Lie $\bbC$-algebra $\frakH$ generated by elements
${\bf1}, b_n,b^*_n$ with $n>0$ and the defining relations
$$[b_n,b_m]=[b^*_n, b^*_m]=0,\quad [b_n,b^*_m]=-n\delta_{n,m}\,{\bf1},\quad n,m>0.$$
The value of $\bf1$ on a given representation is called the \emph{level}.

\smallskip

Let $n$ be a non-negative integer and $\nu$ be a partition of $n$. For any permutation
$x$ in $\frakS_n$ of \emph{cycle-type} $\nu$ let $c_\nu=c_x$ be the conjugacy class
of $x$. We interpret $\frakS_0$ as the trivial group. For any integer $e>0$ we abbreviate 
$$c_n=c_{(n)},
\quad e\nu=(e\nu_1,e\nu_2,\dots).$$

\smallskip

Set $\pmb\Lambda=\bigoplus_{m\in\bbN}\bbC\Irr(\frakS_m)$ and $\langle\bullet,
\bullet\rangle_{\pmb\Lambda}=\bigoplus_{m\in\bbN}\langle\bullet,\bullet\rangle_{\frakS_m}$.
The induction and restriction yield a pair of linear maps $\Ind_{n,m}^{n+m}$,
$\Res_{n,m}^{n+m}$ between the $\bbC$-vector spaces $\pmb\Lambda$ and
$\pmb\Lambda\otimes\pmb\Lambda$. They are adjoint relatively to the scalar product
$\langle\bullet,\bullet\rangle_{\pmb\Lambda}$.
There is a unique representation of $\frakH$ of level $e$ on $\pmb\Lambda$ such that,
for each $\phi\in\bbC\Irr(\frakS_m)$ and $\psi\in\bbC\Irr(\frakS_{m+en})$, we have
\begin{align}\label{Heisenberg}
b_{n}(\phi)=\Ind_{m,en}^{m+en}(\phi\times c_{en}),\quad
b^*_{n}(\psi)=\langle\Res_{m,en}^{m+en}(\psi),c_{en}\rangle_{\frakS_{en}}.
\end{align}
Write $b_\nu=\prod_ib_{\nu_i}$ and 
$b^*_{\nu}=\prod_ib^*_{\nu_i}$.
For each $f$ in $\pmb\Lambda,$ we define $b_f,b_f^*\in U(\frakH)$ by
$$\gathered
b_f=\sum_{\nu\in\scrP}\langle c_\nu,f\rangle_{\pmb\Lambda}\,b_{\nu},\quad
b^*_f=\sum_{\nu\in\scrP}\langle c_\nu,f\rangle_{\pmb\Lambda}\,b^*_{\nu}.
\endgathered$$
If $f=\phi_\nu$, then we will abbreviate 
\begin{align}\label{Heisenberg2} 
a_\nu=b_{\phi_\nu},\quad a^*_\nu=b^*_{\phi_\nu}.
\end{align}

\smallskip

We can now define the representation of the Heisenberg algebra on Fock spaces.
We'll assume that $I=A^{(1)}_{e-1}$.
Hence, since $e$ is odd, we have $\frakg'\simeq\frakg'_e\simeq\widetilde{\fraks\frakl}_e$ and $\frakg\simeq\frakg_e\simeq\widehat{\fraks\frakl}_e$.

\smallskip

First, assume that $l=1$. Then, the representation of $\frakg'$ on
$\bfF(Q)$ admits a commuting action of $\frakH$ of level $e$ such that the representation
of $\frakH\times\frakg'$ is irreducible. The $\bbC$-linear isomorphism $\bfF(Q)\simto\pmb\Lambda$
such that $|\nu,Q\rangle\mapsto \phi_\nu$ identifies the pairing $\langle\bullet,
\bullet\rangle$ with $\langle\bullet,\bullet\rangle_{\pmb\Lambda}$ and the $\frakH$-action 
on $\bfF(Q)$ with \eqref{Heisenberg}.

\smallskip

Now, assume that $l\geqslant 1$. 
By \eqref{eq:EF}, we have a $\frakg'$-module
isomorphism $\bfF(Q)\simto\bigotimes_{p=1}^l\bfF({Q_p})$. 
The representation of $\frakg'$ on  $\bfF(Q)$ admits a commuting action of $\frakH$ of level
$el$ which is given by the tensor product of the representations of $\frakH$ on $\bfF(Q_1),
\bfF(Q_2),\dots,\bfF(Q_l)$. This representation coincides with the specialization of the
representation of the quantized Heisenberg algebra considered in \cite[sec.~4.3]{U}.
See \cite[prop.~4.6]{SV} for more details.

\subsection{The main conjecture}
Assume that $m_{-1}=0$, hence $m$ is a multiple of $e$ and the partition $\lambda$ of $m$  is of the form
$\lambda=e\,\nu$ with $\nu$ a partition of $m/e$.
We abbreviate
$A_{\nu}=B_{S_{{}^t\lambda}}$ and $A^*_{\nu}=\R\! B^*_{S_{{}^t\lambda}}.$

The functors $A_\nu$, $A^*_\nu$ yield linear endomorphisms of the vector space
$[\scrU_\k]$. Using the isomorphism $[\scrU_\k]\simeq\bigoplus_{t\in\bbN}\bfF(Q_t)_e$ we 
can endow $[\scrU_\k]$ with a structure of a representation of level $2e$ of the Heisenberg
algebra $\frakH$ which commutes with the action of $\frakg_e$. By Proposition \ref{prop:wHC},
the set $\{[D]\,$\,|\,$D\in\Irr(\scrU_\k)$ is weakly cuspidal$\}$ is a basis of the subspace
$[\scrU_\k]^{\leqslant 0}=\{x\in\scrU_\k\,|\,E_i(x)=0,\,\forall i\in I_e\}$. We define
\begin{align*}
[\scrU_\k]^\hw
&=\{x\in[\scrU_\k]^{\leqslant 0}\,|\,b_n^*(x)=0,\,\forall n\geqslant 1\},\\
&=\{x\in[\scrU_\k]\,|\,b_n^*(x)=E_i(x)=0,\,\forall n\geqslant 1,\,\forall i\in I_e\}.
\end{align*}
Recall the elements $a_\nu$, $a^*_\nu$ in $U(\frakH)$ introduced in \eqref{Heisenberg2}.

\begin{conjecture}\hfill
 \begin{itemize}[leftmargin=8mm]
  \item[$\mathrm{(a)}$] $\{[D]$\,|\,$D\in\Irr(\scrU_\k)$ is cuspidal$\,\}$ is a basis of
  $[\scrU_\k]^\hw$,
  \item[$\mathrm{(b)}$] $a_\nu=[A_\nu]$ and $a^*_\nu=[A^*_\nu]$ in 
  $\End([\scrU_\k])$.
\end{itemize}
\end{conjecture}

\begin{remark}
Part (a) of the conjecture implies that $\{[D]$\,|\,$D\in\Irr(\scrU_\k)$ is cuspidal$\,\}=
\Irr(\scrU_\k)\cap [\scrU_\k]^\hw$, because $\Irr(\scrU_\k)$ is a basis of $[\scrU_\k]$.
In particular, the subset $\{[D]$\,|\,$D\in\Irr(\scrU_\k)$ is cuspidal$\,\}$ of
$\bigoplus_{t\in\bbN}\bfF(Q_t)_e$ depends only on the integer $e$, and not on the characteristic
$\ell$ of $\k$.
\end{remark}

\end{document}